\documentclass[12pt]{article}
\usepackage[top=30truemm,bottom=30truemm,left=25truemm,right=25truemm]{geometry}
\usepackage{amsmath, amssymb, amsthm}
\usepackage{amscd}
\usepackage{tikz}
		\usetikzlibrary{calc, patterns, intersections}
\usepackage{tikz-cd}
\usepackage{enumerate}
\usepackage{url}
\usepackage{hyperref}
\usepackage{empheq}

\usepackage[all]{xy}
\newtheorem{theorem}{Theorem}[section]

\theoremstyle{definition} %%%
\newtheorem{df}[theorem]{Definition}
\newtheorem{ex}[theorem]{Example}
\newtheorem{rem}[theorem]{Remark}

\theoremstyle{remark}

\theoremstyle{theorem} %%%
\newtheorem{thm}[theorem]{Theorem}
\newtheorem*{inthm}{Theorem}
\newtheorem{lem}[theorem]{Lemma}
\newtheorem{prop}[theorem]{Proposition}
\newtheorem{cor}[theorem]{Corollary}

\newtheoremstyle{yproof}%
	{}{}{\normalfont}{}{\itshape}{\!\!.}{ }%
	{\thmname{#1} \thmnote{of #3\,\,}}
\theoremstyle{yproof}
\newtheorem*{pf}{Proof}
\numberwithin{equation}{section}
\allowdisplaybreaks

\newcommand\alert{\textit}
\renewcommand\qed{\hfill $\Box$}
\newcommand\fin{\hfill $\Diamond$}

\newcommand\Z{\mathbb{Z}}		% the set of integers
\newcommand\lm{\lambda}		% $\lambda$
\newcommand\al{\alpha}			% $\alpha$ root
\newcommand\cmin{b}	% the minimal element of c^{-1}(i)

\newcommand\per[1]{{\hat{#1}}}
\newcommand\sem[1]{{\boldsymbol{#1}}}
\newcommand\cyl[1]{{\mathring{#1}}}

\newcommand\cdiag{\theta}		% cylindric diagram
\newcommand\pdiag{\Theta}		% periodic diagram
\newcommand\peri{\omega}		% period

\newcommand\cyli{\mathcal{C}}	% cylinder
\newcommand\ptn{\mathcal{P}}	% the set of partitions
\newcommand\LE{{\mathrm{ST}}}	% the set of standard tableaux
\newcommand\lex{{\mathfrak{t}}}	% a standard tableau
\newcommand\h{\mathfrak{h}}		% Cartan subalgebra
\newcommand\fdwt{\varpi}			% fundamental domain
\newcommand\con{\mathbf{c}}		% content
\newcommand\predom{\zeta}	% predominant weight corresp to $\cdiag$
\newcommand\bra{\langle}		% the natural pairing
\newcommand\ket{\rangle}		%

\newcommand{\precdot}{
	\mathrel{\prec\!\!\!\cdot}}

\newcommand\rble{\prec}
\newcommand\rbleq{\preceq}

\newcommand\btm{\Gamma}			% the bottom set

\newcommand\hk{{\mathbf{h}}}
\newcommand\hkr{{{\gamma}_{\lex}}}
\newcommand\arm{{\mathrm{Arm}}}
\newcommand\leg{{\mathrm{Leg}}}

\newcommand\I{{\mathcal J}}
\newcommand\ideal{{\xi}}

\newcommand\nmult{{N}}

\newcommand\Supp{{\mathrm{Supp}}}
\newcommand\btmmax{{\btm_{\mathrm{max}}}}
\newcommand\btmmin{{\btm_{\mathrm{min}}}}

\newcommand{\oordneq}{{<^{\mathrm{or}}}}
\newcommand{\hpord}{\leq^{\mathrm{hp}}}

\newcommand\oord{{\leq^{\mathrm{or}}}}

\newcommand\tcord{{{\trianglelefteq^{\mathrm{tc}}}}}

\newcommand\neword{\unlhd}%{{ \triangleleft}}
\newcommand\newordl{\unrhd}%{{\triangleright}}
\newcommand\cylord{\trianglelefteq}

\newcommand\gam{{\gamma}}

\newcommand\im{{\mathrm{Im}}}

\newcommand\DS{{D}}

\newcommand\xs{{x^S}}
\newcommand\xe{{x^E}}
\newcommand\xse{{x^{SE}}}
\begin{document}

%%%%%%%%%%%%%%%%%%%%%%%%%%%%%%%%%%%%%%%%%%%
\title{Poset Structure concerning Cylindric Diagrams}
\author{Kento Nakada\footnote{Graduate School of Education, Okayama University, Japan e-mail address: \texttt{nakada@okayama-u.ac.jp}},
Takeshi Suzuki\footnote{Graduate School of Natural Science and Technology, Okayama University, Japan e-mail address: \texttt{suzuki@math.okayama-u.ac.jp}} and 
Yoshitaka Toyosawa\footnote{e-mail address: \texttt{prkr5rq9@s.okayama-u.ac.jp}}
}
\date{}
\maketitle

\begin{abstract}
The purpose of the present paper is to 
give a realization of a cylindric diagram
as a subset of root systems of type $A_{\kappa-1}^{(1)}$
and several characterization of its poset structure.
Furthermore, the set of order ideals of a cylindric diagram is 
described as a weak Bruhat interval of the Weyl group.
%%%%%%%%%%%%%%%%%%%%%%%%%%%%%%%%%%%%%%%%%%%%%%%%%%%%%%%%%%%%%%%%%%%%%%%%%%%%%%%%%%%%%%%%%%%%%%%
\end{abstract}

\section*{Introduction}

\begin{figure}[h!]
\begin{center}
\begin{tikzpicture}[scale=.32, rotate=270]

%%%%%%%%%%%%%%%%%%%%%%%%%%%%%%%%%%%%%%%%%%%%%%%%%%%%%%%%%%%%%%%%%%%%%%%%%%%%
\fill[pink] (0,-10) |- (1,5) |- (3,3) |- (4,1) -- (4,-10) -- cycle;

\draw[xshift=-8cm, yshift=10cm, ultra thick] (2,3) -- (3,3) |- (4,1) -- (4,-5);
\draw[xshift=-4cm, yshift=5cm, ultra thick] (0,-15) |- (1,5) |- (3,3) |- (4,1) -- (4,-5);
\draw[ultra thick] (0,-10) |- (1,5) |- (3,3) |- (4,1) -- (4,-5);
\draw[xshift=4cm, yshift=-5cm, ultra thick] (0,-5) |- (1,5) |- (3,3) |- (4,1) -- (4,-5);
\draw[xshift=8cm, yshift=-10cm, ultra thick] (0,0) |- (1,5) |- (2,3);

%%%%%%%%%%%%%%%%%%%%%%%%%%%%%%%%%%%%%%%%%%%%%%%%%%%%%%%%%%%%%%%%%%%%%%%%%%%%

\begin{scope}[xshift=-8cm, yshift=10cm]
\draw[dotted, ultra thick] (1,-17) -- +(0,17);
\draw (2,-20) grid (4,1);
\draw (2,1) grid (3,3);
%\draw (0,3) grid (1,5);
\end{scope}

\begin{scope}[xshift=-4cm, yshift=5cm]
\draw[dotted, ultra thick] (2,-16) -- +(0,-3); 
\draw (0,-15) grid (4,1);
\draw (0,1) grid (3,3);
\draw (0,3) grid (1,5);
\end{scope}

\draw[dotted, ultra thick] (2,-11) -- +(0,-3);
\draw (0,-10) grid (4,1);
\draw (0,1) grid (3,3);
\draw (0,3) grid (1,5);

\begin{scope}[xshift=4cm, yshift=-5cm]
\draw[dotted, ultra thick] (2,-6) -- +(0,-3);
\draw (0,-5) grid (4,1);
\draw (0,1) grid (3,3);
\draw (0,3) grid (1,5);
\end{scope}

\begin{scope}[xshift=8cm, yshift=-10cm]
%\draw (0,-5) grid (4,1);
\draw (0,0) grid (2,3);
\draw (0,3) grid (1,5);
\draw[dotted, ultra thick] (3,0) -- +(0,3);
\end{scope}

\draw[dashed] (0,5) -- +(4,0) -- +(4,-5);
%\draw (0,5) to [out=60, in=120] +(4,0) to [out=-60, in=60] +(0,-5);
\draw (0,5) sin (2,6) cos  (4,5);
\node at (2,6) [fill=white,inner sep=5pt] {4};
\draw (4,5) cos (5,2.5) sin (4,0);
\node at (5,2.5) [fill=white,inner sep=5pt] {5};

\draw[thick, ->] (3,4) -- +(4,-5);
\node at (7,1.5) {$(4,-5)$};

\end{tikzpicture}
\end{center}

%\caption{Contents on cylindric diagram associated with $\lm=(5,4,4,2) \in \ptn_{4,4}$.}\label{fig1}
\caption{
  A periodic diagram of period $\peri=(4,-5)$.} %with respect to $\lm=(9,7,7,5)$.}
  \label{periodic diagram}
\end{figure}
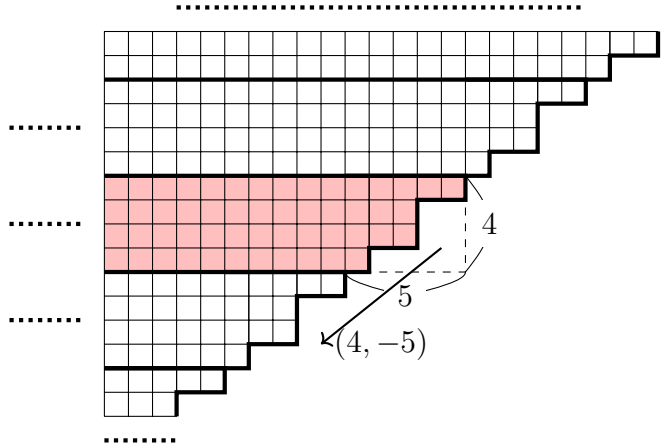

%Let $\peri\in \Z_{\geqq 1}\times\Z_{\leqq -1}$.
A {\it periodic $($Young$)$ diagram} %of period $\peri\in \Z_{\geqq 1}\times\Z_{\leqq -1}$ 
%skew diagram of period $\peri$ 
is a Young diagram consisting of infinitely 
many cells in $\Z^2$ which is invariant under parallel 
translations generated by a certain vector $\peri\in \Z^2$ called the period
(see Figure \ref{periodic diagram}).
The image 
of a periodic diagram under 
the natural projection 
%$\pi:\Z^2\to \Z^2/\Z\peri$ 
onto the cylinder $\Z^2/\Z\peri$
is called a \alert{cylindric diagram}.
Diagrams given as a set-difference of
two cylindric diagrams are called  cylindric skew diagrams.

%The set of standard tableaux on a periodic diagram can be identified 
%with the set of  standard tableaux on the corresponding cylindric 
%skew diagram.
We note that cylindric skew diagrams have been known to 
parameterize a
certain class of irreducible modules over 
the Cherednik algebras (double affine Hecke algebras)
(\cite{Suz2005,SuzVaz2005})
and the (degenerate) affine Hecke algebras (\cite{Kle1996, Ruff2006}) of type $A$,
where standard tableaux on those diagrams also appear.
%We also note that %the authors presented in \cite{SuzToy2022} the conjecture of the hook formula giving the number of cylindric standard tableaux on cylindric skew diagrams.
%the conjecture of hook formulas for cylindric skew diagrams was given in \cite{SuzToy2022}.

Let $\peri=(m,-\ell)\in\Z_{\geqq 1}\times\Z_{\leqq -1}$ and
let $\cdiag$ be a cylindric diagram in $\Z^2/\Z\peri$.
%of period $\peri$.
The lattice $\Z^2$ admits a partial order $\leq$ defined by
$$(a,b)\leq (c,d) \iff a\geqq c\text{ and }\ b\geqq d,$$
which induces %and it induces 
a poset structure  on  $\Z^2/\Z\peri$
and also on $\cdiag$.
Together with the content map $\con:\cdiag\to \Z/\kappa\Z$, where
$\con(a,b)=b-a\mod \kappa$ and  $\kappa=\ell+m$,
the cylindric digram $\cdiag$ is 
a locally finite $\Z/\kappa\Z$-colored
$d$-complete poset in the sense of \cite{Str2020,Str2021}.

The purpose of the present paper is to investigate
the poset $(\cdiag,\leq)$ as well as the poset
$(\I(\cdiag),\subset)$, 
where $\I(\cdiag)$ denotes the set of cylindric skew diagrams (or proper order ideals) included in $\cdiag$.
%in other words, %the set of 
%proper order ideals of $\cdiag$.

%\if0
%%%%%%%%%%%%%%%%%%%%%%%%%%%%%%%%%%%%%%%%%%%%%%%%%%%%%%%%%%%%%%%%
We briefly review a description in  the classical case.
Let $\lm\subset\Z^2$ be a finite Young diagram. 
The associated  Grassmannian permutation $w_\lm$
 is an element of the Weyl group of the root system $R$ of type $A_{n}$
 where  $n=\sharp\{\con(x)\mid x\in \lm\}$.
It is known %has known been  shown  by Stembridge \cite{Stem2001}
that
%It holds that
 the poset  $(\lm,\leq)$ is dually isomorphic to the poset 
$(R(w_\lm^{-1}),\oord)$, where 
$R(w_\lm^{-1}):=R_+\cap w_\lm^{-1} R_-$
%is the inversion set of $w$
% where   $R_+$ (resp. $R_-$) is the set of positive (resp. negative) roots
%Moreover the map $\hk$ gives an isomorphism $(\lm,\leq)\cong (R(w_\lm),\oord)$ as posets,
and $\oord$ is the ordinary order (or the standard order) 
defined by
%as the transitive closure of the relation
$$\al\oord \beta \iff  \beta-\al\in \sum_{i \in [1,n]}\Z_{\geqq 0}\al_i$$
for $\al,\beta\in R(w_\lm)$
%where $\{\al_1,\al_2,\dots,\al_{n-1}\}$
%is  the set of simple roots.
with $\Pi$ being the set of simple roots (\cite{Stem2001}).

Let $\cdiag$ be a cylindric diagram
in $\Z^2/\Z\peri$. %of period $\peri=(m,-\ell)$.
We would like to describe the poset $(\cdiag,\leq)$
in terms of the root system of type $A_{\kappa-1}^{(1)}$ with $\kappa=\ell+m$.

A key ingredient in our approach is 
the {\it colored hook length} (\cite{Nak2008, Proc2014}), %which have appeared in \cite{Proc2014, Nak2008}, 
given by
$$\hk(x)=\sum_{y\in H(x)}\al_{\con(y)}\ \ (x\in\cdiag),$$
where $H(x)$ denotes the hook at $x$ and $\al_i$  are  simple roots.
(See Section \ref{ss:hook} for precise definitions.)
We will show that the map $\hk$ embeds
the cylindric diagram $\cdiag$
into the set $R_+$ of positive (real) roots,
and that the image $\hk(\cdiag)$ is given by the inversion set $R(w_\cdiag)$ 
associated with a semi-infinite word $w_\cdiag$, which can be thought as
an 
%semi-infinite
 analogue of the Grassmannian permutation.
%%%%%%%%
Moreover, we show that the image $\hk(\cdiag)$ is also characterized as
the subset %$\DS(\predom_\cdiag)$ 
of $R_+$ consisting of those elements satisfying
$$\bra\predom_\cdiag,\al^\vee\ket=-1,$$
where $\predom_\cdiag$ is 
 a predominant integral weight determined by $\cdiag$ (see Section \ref{sec:predom} and \ref{sec:inversion} for details).

%%%%%%%%%%%%%%%%%%%%%%%%%%%%%%%%%%%%%%%%%%%%%%%%%%%%%%%%%%%%%%%%%%%%%%%%%%%%%

%We will give  a natural generalization of the colored hook length $\hk$ 
%and the inversion set $R(w_\lm)$ and then 
%obtain a bijection
%$$\hk:\cdiag\to R(w_\cdiag).$$
%\fi

Unlike the classical case, 
the ordinary order in $R(w_\cdiag)$ does not lead a poset isomorphism,
and we need to introduce a modified order $\neword$ in $R(w_\cdiag)$ by
$$\al\oord \beta \iff
\beta-\al\in \sum_{\gamma\in \Pi_\cdiag}\Z_{\geqq 0}\gamma,$$
to obtain a poset isomorphism $(\cdiag,\leq)\cong
(R(w_\cdiag),\neword)$, where $\Pi_\cdiag$ is a certain subset of the affine root system (see Section \ref{sec:inversion set}).

Another description of the poset $\cdiag$ is given by 
a linear extension or (reverse) standard tableau $\lex$
on $\cdiag$, which is by definition a bijective order preserving map $\cdiag\to\Z_{\geqq1}$.
%A standard tableau 
A linear extension $\lex:\cdiag\to\Z_{\geqq1}$ brings a poset structure 
to $\Z_{\geqq1}$ and the resulting poset is an infinite analogue of the heap,
which is originally introduced by Stembridge \cite{Stem2001}.
%
%
\if0
%For a cylindric diagram $\cdiag$ and a standard tableau $\lex$,
%the following posets are associated:
\begin{itemize}
\item $(\cdiag,\le)$, 
\item $(R(w_\cdiag),\neword)=(D(\predom_\cdiag),\neword)$,
\item $(\Z_{\geqq1},\hpord_\lex)$,
\end{itemize}
where the partial order $\neword$ is a new order to be introduced in Section \ref{sec:1stresult}
and the partial order $\hpord_\lex$ is an analogy to cylindric diagrams of the heap ordering.
The set $D(\predom_\cdiag)$ is a subset of positive roots determined by some integral weight $\predom_\cdiag$ (see Section \ref{sec:hooks}).
 %and a standard tableau $\lex$ on $\cdiag$,
%some integral weight $\predom_\cdiag$ and the inversion set $R(w_\cdiag)$ determine.
%The integral weight $\predom_\cdiag$ is pre-dominant, namely,
%$\bra \predom_\cdiag, \al^\vee \ket \geqq -1$ for all $\al\in R_+$.
The colored hook length $\hk$ can extend to cylindric diagrams.
To describe of $\im\hk$,
we will show that the map $\hk$ gives a bijection
from $\cdiag$ to $D(\predom_\cdiag)$ (Theorem \ref{thm:hkinR+}).
The set $R(w_\cdiag)$ is a subset of positive roots associated a ``semi-infinite word'' $w_\cdiag$.
Moreover, the set $D(\predom_\cdiag)$ coincides with the inversion set $R(w_\cdiag)$ (see Section \ref{sec:hooks}).
\fi
%
In summary,  we  have the following:
\begin{inthm}[Theorem \ref{th:cylord} and Proposition \ref{prop:cdiagandheap}]
The followings are poset isomorphisms:
\[ (\Z_{\geqq1},\hpord_\lex) \stackrel{\lex}{\leftarrow} (\cdiag,\le) \stackrel{\hk}{\rightarrow} (R(w_\cdiag),\neword). \]
\end{inthm}

Another goal of this paper is to
describe the poset structure %of the set 
$\I(\cdiag)$.
%of cylindric skew diagrams included in $\cdiag$
%or equivalently, the set of 
%(or proper order ideals of $\cdiag$)
%with the partial order given by the inclusion relations.
For a finite Young diagram $\lm$,
%it is determined some dominant integral weight $\zeta_\lm$.
%Let $W^\lm$ denote the set of the minimal length coset representatives of $W/W_\lm$,
%where $W_\lm=\{w \in W \mid w \zeta_\lm=\lm\}$.
it is known that the set $\I(\lm)$ of order ideals
 of $\lm$  is isomorphic to
the interval $[e,w_\lm]=\{u \in W \mid e \preceq u \preceq w_\lm\}$ 
with weak right Bruhat order (\cite[Proposition I]{Proc2014}).
For a cylindric diagram $\cdiag$,
we define a ``semi-infinite Bruhat interval'' $[e,w_{\cdiag})$,
and we have the following:
%For a cylindric diagram $\cdiag$,
%some integral weight $\predom_\cdiag$ is determined.
%In particular, it is predominant, namely, $\bra \predom_\cdiag, \al^\vee \ket \geqq -1$ for all $\al \in R_+$.
%As the second result,
%we will see the following:
\begin{inthm}[Theorem \ref{th:ideal}]
The map 
$$\Phi: (\I(\cdiag),\subset) \to ([e,w_{\cdiag}),\rbleq)$$
given by $\Phi(\ideal) = w_{\ideal}$ is a poset isomorphism.
\if0
%The set $\I(\cdiag)$ of order ideals of $\cdiag$ is poset isomorphic to the semi-infinite Bruhat interval $[e,w_\cdiag)$ with weak right Bruhat order:
%\[ (\I(\cdiag),\subset) \cong ([e,w_\cdiag),\rbleq). \]
The following is a poset isomorphism:
\[ (\I(\cdiag),\subset) \stackrel{\Phi}{\rightarrow} 
([e,w_{\cdiag}),\preceq)
% \stackrel{\Psi}{\rightarrow} (\I(R(w_\cdiag)),\subset), 
\]
where 
$\Phi(\ideal)=w_{\ideal}$. %and $\Psi(w)=R(w)$.
\fi
\end{inthm}
%We will see that the poset isomorphism from the set $\I(\cdiag)$ to
%the subset of Weyl group such that $\bra \predom_\cdiag, \al^\vee \ket = -1$
%for all $\al \in R(w_\cdiag)$,
%which is an opposite concept of minuscule element, with Bruhat ordering
%(Theorem \ref{thm1}).

%We consider the poset structure of $\I(\cdiag)$.
%For an integral weight $\predom_\cdiag$,
%For a cylindric diagram $\cdiag$,
%some integral weight $\predom_\cdiag$ is determined.
%In particular, it is predominant, namely, $\bra \predom_\cdiag, \al^\vee \ket \geqq -1$ for all $\al \in R_+$.
%We will describe that the $\predom_\cdiag$-``pluscule'' elements which is an opposite concept of minuscule elements (see Section 4)
%are one-to-one correspondence with proper order ideals of $\cdiag$.
%An element $w$ of Weyl group satisfies $\bra \predom_\cdiag, \al^\vee \ket = -1$
%for all $\al \in R(w)$ is called a $\predom_\cdiag$-``pluscule'' element.
%Note that it is an opposite concept of minuscule element.
%We will see that the poset isomorphism from the set $\I(\cdiag)$ to
%the subset of Weyl group such that $\bra \predom_\cdiag, \al^\vee \ket = -1$
%for all $\al \in R(w_\cdiag)$,
%which is an opposite concept of minuscule element, with Bruhat ordering.

%%%%%%%%%%%%%%%%%%%%%%%%%%%%%%%%%%%%%%%%%%%
\section{Cylindric diagrams}
%%%%%%%%%%%%%%%%%%%%%%%%%%%%%%%%%%%%%%%%

%%%%%%%%%%%%%%%%%%%%%%%%%%%%%%%%%%%%%%%%%%%%%%%%%%
\subsection{Cylindric diagrams as posets}\label{SS:cylindric}
%%%%%%%%%%%%%%%%%%%%%%%%%%%%%%%%%%%%%%
Let $(P,\le)$ be a poset.
For $x,y\in P$, define an \alert{interval} $[x,y]$ by
$$[x,y] = \{ z \in P \mid x \leq z \leq y\}.$$
We say that $y$ \alert{covers} $x$ if $[x,y]=\{x,y\}$.
%An interval is called a \alert{chain} if it is a totally ordered set.
%
%In the rest,
%we assume that a poset $P$ is \textit{locally finite}, namely,
%$\#[x,y]<\infty$ for any $x,y \in P$.

\begin{df}
Let $(P,\le)$ be a poset.
A subset $J$ of $P$ is called an \alert{order filter} (resp. \alert{order ideal})
if the following condition holds:
$$x \in J,\ x \le y \implies y \in J
\qquad (\text{resp. } x \in J,\ x \ge y \implies y \in J).$$

\if0
A subset $F$ of $P$ is called an \alert{order filter}
if the following condition holds:
$$x \in F,\ x \leqq y \implies y \in F.$$
%An order filter $F$ is said to be \alert{non-trivial} if $F\neq \emptyset$ nor $F\neq P$.

A subset $J$ of $P$ is called an \alert{order ideal}
if the following condition holds:
$$x \in J,\ x \geqq y \implies y \in J.$$
\fi

An order filter (resp. order ideal) $J$ is said to be \alert{proper} if
$J \neq P$,
and it is said to be \alert{non-trivial} if $J \neq P$ nor $J \neq \emptyset$.

%An order filter (resp. ideal $J$) is said to be \alert{non-trivial} if $F\neq \emptyset$ (resp. $J\neq \emptyset$) nor $F\neq P$ (resp. $J\neq P$).
%
%For $x\in P$, put
%For $x_1,x_2,\dots\in P$, the set 
%$\bigcup_iF_{x_i}$, where $F_x:=\{y\in P\mid x\geqq y\}$ 
%is an order filter of $P$ and it is called the order filter
%generated by $x_1,x_2,\dots$.
%The condition (Y3) above means that   
%$\sidiag$ is an order filter of $\Z^2$,
%and a semi-infinite periodic diagram
%is a $\Z\omega$-orbit of an order filter generated by finitely many elements.
\end{df}

For $\peri \in \Z_{\geqq 1} \times \Z_{\leqq -1}$,
we let $\Z \peri$ denote the subgroup of (the additive group) $\Z^2$ generated by $\peri$,
and define the cylinder $\cyli_\peri$ by
$$
\cyli_\peri = \Z^2 / \Z \peri.
$$
Let $\pi:\Z^2 \to \cyli_\peri$ be the natural projection.
The cylinder $\cyli_\peri$ inherits a $\Z^2$-module structure via $\pi$.
%The group $\Z^2$ acts on $\cyli_\peri$ by
%\begin{align*}
%\Z^2\times \cyli_\peri &\to \cyli_\peri\\
%(u,x)\ \ &\mapsto x+u:=\pi(\tilde x+u),
%\end{align*}
%where $\tilde x$ is any element of $\pi^{-1}(x)$.
%
%
%We regardas

Define a poset structure on   $\Z^2$ 
by 
%the following partial order
$$
(a,b) \leq (a', b') \iff \hbox{$a \geqq a'$ and $b \geqq b'$ 
as integers}.
$$
For $x,y \in \cyli_\peri$,
write $x \le y$ if there exists $\tilde{x}, \tilde{y} \in \Z^2$ such that
$\pi(\tilde{x})=x$, $\pi(\tilde{y})=y$ and $\tilde{x} \le \tilde{y}$.
It is not difficult to see the following:
%%%%%%%%%%%%%%%%%%%%%%%%%%%%%%%%%%%%%%%%%
\begin{lem}
Let $\peri \in \Z_{\geqq 1} \times \Z_{\leqq -1}$. Then
the relation $\leq$ on $\mathcal{C}_\omega$ is a partial order,
and the projection $\pi$ is order preserving.
\end{lem}
%%%%%%%%%%%%%%%%%%%%%%%%%%%%%%%%%%%%%%%%%%%

In the rest of this section, we fix  $\omega \in \mathbb{Z}_{\geqq1} \times \mathbb{Z}_{\leqq -1}$.

%%%%%%%%%%%%%%%%%%%%%%%%%%%%%%%%%%%%%%%%%%%%%%%%%%%%%%%
\begin{df}

\noindent
(1) A non-trivial order filter of $\mathcal{C}_\omega$ is called a \alert{cylindric diagram}.

%\noindent
%(2) A proper order ideal of a cylindric diagram is called a \alert{cylindric skew diagram}.

\noindent
(2)
A non-trivial order filter $\Theta$ of $\Z^2$ is called a \alert{periodic diagram of period $\peri$} if $\Theta + \omega = \Theta$.
\end{df}
%%%%%%%%%%%%%%%%%%%%%%%%%%%%%%%%%%%%%%%%%%%%%%%%%%%%%%%%%%%%%

%%%%%%%%%%%%%%%%%%%%%%%%%%%%%%%%%%%%%%%%%%%%%%
\begin{lem}\label{lem;CYD}
%Let $\peri\in\Z_{\geqq1} \times \Z_{\leqq -1}$.

\noindent
{\rm (1)} For a cylindric diagram $\cdiag$ in $\cyli_\peri$,
the inverse image $\pi^{-1}(\cdiag)$ is a periodic diagram of period $\peri$.

\noindent
{\rm (2)} For a periodic diagram $\Theta$ of period $\peri$,
the image $\pi(\Theta)$ is a cylindric diagram in $\cyli_\peri$. 
\end{lem}
%%%%%%%%%%%%%%%%%%%%%%%%%%%%%%%%%%%%%%%%%%%%%%%%%%%%%%%%%%%%%%%

Figure \ref{periodic diagram} indicates a periodic diagram of period $\peri=(4,-5)$.
The set consisting of colored cells is a fundamental domain with respect to the action of $\Z\peri$,
and it is in one to one correspondence with the associated cylindric diagram.

\begin{df}
%Let $m,\ell\in\Z_{\geqq1}$.
Let $m,\ell \in \Z_{\geqq1}$.
A non-increasing sequence $\lambda=(\lambda_1,\dots,\lambda_m)$ of (possibly negative) integers is called a \alert{generalized partition of length $m$}.
%and it
%is said to be \alert{$\ell$-restricted} % generalized partition of length $m$}
%if it satisfies %the following conditions:
For $\peri=(m,-\ell) \in \Z_{\geqq1} \times \Z_{\leqq-1}$,
we 
denote by $\ptn_{\peri}$ the set of %$\ell$-restricted 
generalized partitions of length $m$ satisfying 
$$
%\lambda_1 \geqq \cdots \geqq \lambda_m,\quad
\lambda_1-\lambda_m \leqq \ell.
$$
\end{df}
%\def\sinf{\frac{\infty}{2}}

%We call $\lm$ the \alert{fundamental domain} o $\perid{\lm}$.

For %a generalized partition  
$\lambda=(\lambda_1,\dots,\lambda_m) \in \ptn_\peri$,
% \in \P_{m,\ell}$.
we define
\begin{align*}
\sem\lambda &= \{ (a,b) \in \mathbb{Z}^2 \mid 1 \leqq a \leqq m,\ b \leqq \lambda_a \}, \\
\per{\lambda} &
%=\per{\lm}_{(m,-\ell)}
= \sem\lambda + \mathbb{Z} \peri,\\
\cyl\lm &
%=\cyl\lm_{(m,-\ell)}
=\pi(\per\lm).
\end{align*}
Note that  
$\sem\lambda = \per{\lambda} \cap ([1,m] \times \mathbb{Z})$ 
and $\sem\lm$ is a fundamental domain of $\per\lambda$ with respect to the action of
$\mathbb{Z} (m,-\ell)$. %(see Figure \ref{fig1}).

\if0
It is easy to see that $\per{\lm}$ is a periodic diagram of period  $\peri=(m,-\ell)$, and 
that any diagram of period $(m,-\ell)$ is of the form $\per\lm$ for some 
 $\lm \in \ptn_{\peri}$.
\fi

%It is easy to see the fllowing:

%\begin{prop}
%$\mathrm{(1)}$ 
If  $\lm\in \ptn_{\peri}$ then $\per{\lambda}$ is a periodic diagram of period  $\peri$
and $\cyl\lm$ is a cylindric diagram.
Moreover,
 any periodic (resp. cylindric)
  diagram of period $\peri$ is of the form $\per\lambda$
  (resp. $\cyl\lm$) for some 
 $\lambda \in \ptn_{\peri}$.

For a poset $P$ and its order filter $J$,
we denote the set-difference $P\setminus J$ also by $P/J$.
It is easy to see the following:
%
%%%%%%%%%%%%%%%%%%%%%%%%%%%%%%%%%%%%%%%%%%%%%%%%%%%%%%%%%%%%%%%%%%
\begin{prop}\label{prop:skew}

For a  subset $\xi$ of $\cyli_\peri$, the following conditions are equivalent $:$
%\begin{enumerate}[\rm(i)]
%\item

\noindent
$\mathrm{(i)}$ $\xi$ is a proper order ideal of a cylindric diagram in $\cyli_\peri$.

\noindent
$\mathrm{(ii)}$ $\xi$ is a set-difference $\cdiag/\eta$ of
two cylindric diagrams $\cdiag,\eta$ in $\cyli_\peri$ with $\cdiag\supset \eta$.

\noindent
$\mathrm{(iii)}$ $\xi$ is an intersection of a proper order ideal and 
a proper order filter of $\cyli_\peri$.

\noindent
%\item
$\mathrm{(iv)}$ $\xi$ is a finite subset of $\cyli_\peri$
and  satisfies the following condition:
$$x,y \in \xi \implies [x,y] \subset \xi.$$

\noindent
%\item 
$\mathrm{(v)}$ $\xi$ is a finite subset of $\cyli_\peri$
and satisfies the following condition:
$$x,x+(1,1) \in \xi \implies
x+(0,1),x+(1,0) \in \xi\ \ (\text{the skew property})$$
%\end{enumerate}
\end{prop}
%%%%%%%%%%%%%%%%%%%%%%%%%%%%%%%%%%%%%%%%%%%%%%%%%%%%%%%%%%%%%%%%%%%%%%

\begin{df}
A  subset $\xi$ of $\cyli_\peri$ is called a \alert{cylindric skew diagram}
if it satisfies one of the conditions (i)--(v) in Proposition \ref{prop:skew}.
\end{df}

\begin{figure}[h]
\centering
\begin{tikzpicture}[scale=.5, rotate=270]

\fill[pink] (0,5) |- (1,2) |- (2,1) |- (3,0) |- (4,-2) |- (3,1) |- (1,3) |- (0,5) -- cycle;

%%%%%%%%%%%%%%%%%%%%%%%%%%%%%%%%%%%%%%%%%%%%%%%%%%%%%%%%%%%%%%%%%%%%%%%%%%%%
% \lambda
\draw[xshift=-4cm, yshift=5cm, ultra thick, black!60] (2,3) -- (3,3) |- (4,1) -- (4,-2);
%\draw[ultra thick] (0,-7) |- (1,5) |- (3,3) |- (4,1) -- (4,-7);
\draw[xshift=4cm, yshift=-5cm, ultra thick, black!60]  (0,2) |- (1,5) |- (2,3);

\draw[ultra thick] (0,2) |- (1,5) |- (3,3) |- (4,1) -- (4,-2);

% \mu
\draw[xshift=-4cm, yshift=5cm, ultra thick, black!60] (2,1) |- (3,0) |- (4,-2);
\draw[ultra thick] (0,2) -- (1,2) |- (2,1) |- (3,0) |- (4,-2);
\draw[xshift=4cm, yshift=-5cm, ultra thick, black!60] (0,2) -- (1,2) |- (2,1);
%\fill[green] (0,2) |- (1,5) |- (3,3) |- (4,1) |- (3,-2) |- (2,0) |- (1,1) |- (0,2) -- cycle;

%%%%%%%%%%%%%%%%%%%%%%%%%%%%%%%%%%%%%%%%%%%%%%%%%%%%%%%%%%%%%%%%%%%%%%%%%%%%

\if0
\begin{scope}[xshift=-4cm, yshift=5cm]
\draw[dotted, ultra thick] (1,-10) -- +(0,10);
\draw[black!60] (2,-12) grid (4,1);
\draw[black!60] (2,1) grid (3,3);
%\draw (0,3) grid (1,5);
\end{scope}
\fi

\draw[xshift=-4cm, yshift=5cm, black!60] (2,0) grid (3,3);
\draw[xshift=-4cm, yshift=5cm, black!60] (3,-2) grid (4,1);

%\draw[dotted, ultra thick] (2,-8) -- +(0,-3);
%\draw (0,-7) grid (4,1);
\draw (0,2) grid (1,5);
\draw (1,1) grid (2,3);
\draw (2,0) grid (3,3);
\draw (3,-2) grid (4,1);

\draw[xshift=4cm, yshift=-5cm] (0,2) grid (1,5);
\draw[xshift=4cm, yshift=-5cm] (1,1) grid (2,3);

\if0
\begin{scope}[xshift=4cm, yshift=-5cm]
%\draw (0,-5) grid (4,1);
\draw[black!60] (0,-2) grid (2,3);
\draw[black!60] (0,3) grid (1,5);
\draw[dotted, ultra thick] (3,-1) -- +(0,3);
\end{scope}
\fi
%%%%%%%%%%%%%%%%%%%%%%%%%%%%%%%%%%%%%%%%%%%%%%%%%%%%%%%%%%%%%%%%%%%%%%%%%%%%

\draw[xshift=-4cm, yshift=5cm, dotted, ultra thick] (1,2) -- +(0,3);
\if0
\coordinate (q) at (.5,.5);

\foreach \x/\y/\n in {
               3/1/8, 3/2/2, 3/3/1,
4/-1/9, 4/0/7, 4/1/4
}
\node[black!60] at ($(\x,\y)-(q)+(-4,5)$) {$\mathbf{\n}$};

\foreach \x/\y/\n in {
                              1/3/11, 1/4/6, 1/5/3,
                      2/2/10, 2/3/5,
               3/1/8, 3/2/2, 3/3/1,
4/-1/9, 4/0/7, 4/1/4
}
\node at ($(\x,\y)-(q)$) {$\mathbf{\n}$};

\foreach \x/\y/\n in {
                              1/3/11, 1/4/6, 1/5/3,
                      2/2/10, 2/3/5,
}
\node[black!60] at ($(\x,\y)-(q)+(4,-5)$) {$\mathbf{\n}$};
\fi

\draw[xshift=4cm, yshift=-5cm, dotted, ultra thick] (3,-1) -- +(0,3);

\end{tikzpicture}

\caption{A cylindric skew diagram.}\label{fig:skew diag}
\end{figure}
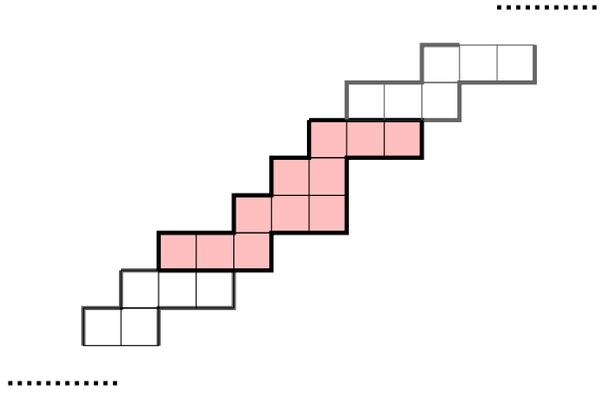

We  denote the set of proper order ideals of $\cdiag$
by $\I(\cdiag)$ and regard it
as a poset with the inclusion relation.
%We will describe the poset structure of $(\I(\cdiag),\subseteq)$
%in terms of the semi-infinite word $w_\cdiag$.
Note that any $\ideal\in \I(\cdiag)$ is a finite set
and thus $\I(\cdiag)=\bigsqcup_{n=0}^\infty \I_n(\cdiag)$,
where we put 
$$\I_n(\cdiag)=\{\ideal\in\I(\cdiag)\mid |\ideal|=n\}.$$
%
%
%
%%%%%%%%%%%%%%%%%%%%%%%%%%%%%%%%%%%%%%%%%%%%%%%%%%%%%%%%%%
\subsection{Standard tableaux} 
%%%%%%%%%%%%%%%%%%%%%%%%%%%%%%%%%%%%%%%%%%%%%%%%%%%%%%%%%

%Let $\omega\in \Z_{\geqq 1}\times \Z_{\leqq -1}$ and 
In the rest of present section, fix a cylindric  diagram $\cdiag$ in $\cyli_\peri$.
%For a cylindric diagram of period $\omega$, its bottom set $\Gamma$ consists
%of $\kappa$ elements (see Figure \ref{content})
%and 
%Note that the restriction of the content map
%gives a bijection $\con:\Gamma\to \Z/\kappa\Z$ (see Figure \ref{content}). 

%%%%%%%%%%%%%%%%%%%%%%%%%%%%%%%%%%%%%%%%%%%%%%%
\begin{df}
\noindent
(1)
For a cylindric diagram $\cdiag$,
a \alert{standard tableau} (or {\it linear extension}) 
of $\cdiag$ is a bijection $\lex:\cdiag \to \Z_{\geqq1}$ satisfying
$$
x<y \implies \lex(x)<\lex(y).
$$
We denote by $\LE (\cdiag)$ the set of standard tableaux of $\cdiag$.

\noindent
(2)
For a finite poset $P$ with $|P|=n$,
a standard tableau
of $P$ is a bijection 
$\lex:P \to [1,n]$ satisfying
$$
x<y \implies \lex(x)<\lex(y).
$$
We denote by $\LE (P)$ the set of standard tableaux of $P$.
\end{df}
%%%%%%%%%%%%%%%%%%%%%%%%%%%%%%%%%%%%%%%%%%%%%%%

\begin{figure}[h]
\centering
\begin{tikzpicture}[rotate=-90, scale=.5]
\newcommand\leftmost{-3}

\begin{scope}[black!70,xshift=-2cm,yshift=3cm]
%\draw[ultra thick] (1,\leftmost-3) |- (2,5) |- (3,4) -- (3,\leftmost-3);
\draw (1,2) grid +(1,3);
\draw (2,1) grid +(1,3);
\coordinate (c) at (.5,.5);
\foreach \n/\a/\b in {1/1/4, 2/2/3, 3/2/2, 4/1/3, 5/2/1, 6/1/2}
\node at ($(\a,\b)+(c)$) {$\n$};
\end{scope}

\begin{scope}[black!70,xshift=2cm,yshift=-3cm]
%\draw[ultra thick] (1,\leftmost+3) |- (2,5) |- (3,4) -- (3,\leftmost+3);
\draw (1,2) grid +(1,3);
\draw (2,1) grid +(1,3);
\coordinate (c) at (.5,.5);
\foreach \n/\a/\b in {1/1/4, 2/2/3, 3/2/2, 4/1/3, 5/2/1, 6/1/2}
\node at ($(\a,\b)+(c)$) {$\n$};
\end{scope}

\fill[red!30] (1,2) rectangle +(1,3);
\fill[red!30] (2,1) rectangle +(1,3);
%\draw[ultra thick] (1,\leftmost) |- (2,5) |- (3,4) -- (3,\leftmost);
\draw (1,2) grid +(1,3);
\draw (2,1) grid +(1,3);
%\draw[dotted, ultra thick] (2,0) -- +(0,-2);
\coordinate (c) at (.5,.5);
\foreach \n/\a/\b in {1/1/4, 2/2/3, 3/2/2, 4/1/3, 5/2/1, 6/1/2}{
%\fill[white!0] ($(\a,\b)+(c)$) circle (.4);
\node at ($(\a,\b)+(c)$) {$\n$};
}
\node at (6,3) {Standard tableau};

\begin{scope}[yshift=12cm]
\begin{scope}[black!70,xshift=-2cm,yshift=3cm]
%\draw[ultra thick] (1,\leftmost-3) |- (2,5) |- (3,4) -- (3,\leftmost-3);
\draw (1,2) grid +(1,3);
\draw (2,1) grid +(1,3);
\coordinate (c) at (.5,.5);
\foreach \n/\a/\b in {4/1/4, 1/2/3, 2/2/2, 5/1/3, 3/2/1, 6/1/2}
\node at ($(\a,\b)+(c)$) {$\n$};
\end{scope}

\begin{scope}[black!70,xshift=2cm,yshift=-3cm]
%\draw[ultra thick] (1,\leftmost+3) |- (2,5) |- (3,4) -- (3,\leftmost+3);
\draw (1,2) grid +(1,3);
\draw (2,1) grid +(1,3);
\coordinate (c) at (.5,.5);
\foreach \n/\a/\b in {4/1/4, 1/2/3, 2/2/2, 5/1/3, 3/2/1, 6/1/2}
\node at ($(\a,\b)+(c)$) {$\n$};
\end{scope}

\fill[red!30] (1,2) rectangle +(1,3);
\fill[red!30] (2,1) rectangle +(1,3);
%\draw[ultra thick] (1,\leftmost) |- (2,5) |- (3,4) -- (3,\leftmost);
\draw (1,2) grid +(1,3);
\draw (2,1) grid +(1,3);
%\draw[dotted, ultra thick] (2,0) -- +(0,-2);
\coordinate (c) at (.5,.5);
\foreach \n/\a/\b in {4/1/4, 1/2/3, 2/2/2, 5/1/3, 3/2/1, 6/1/2}{
%\fill[white!0] ($(\a,\b)+(c)$) circle (.4);
\node at ($(\a,\b)+(c)$) {$\n$};
}
\node at (6,3) {NOT standard tableau};
\end{scope}
\end{tikzpicture}
\caption{}
\end{figure}

\begin{rem}
Our standard tableaux are usually referred to as reverse standard tableaux.
\fin
\if0
\smallskip\noindent
%We also note that %
(2) T.Suzuki and Y.Toyosawa proposed a conjectural  hook formula 
concerning the number of standard tableaux on cylindric skew diagrams in
 \cite{SuzToy2022} 
%the conjecture of hook formulas for cylindric skew diagrams was given in \cite{SuzToy2022}.
\fi
\end{rem}

Let $\lex\in \LE(\cdiag)$. 
It is easy to see that 
the subset  $\lex^{-1}([1,n])$ of $\cdiag$ is
a proper order ideal,
% of $\cdiag$
and moreover the restriction 
$\lex|_{\lex^{-1}([1,n])}$ is a standard tableau on $\lex^{-1}([1,n])$.
%%%%%%%%%%%%%%%%%%%%%%%%%%%%%%%%%%%%%%%%%%%%%%%%%%%%%%%%%%%%%%%%%%%%%%%%%%%
%
Conversely, for $\ideal\in\I_n(\cdiag)$, any standard tableau
on $\ideal$ can be extended to a standard tableau on $\cdiag$.
In summary,  we have the following:
% as before.
%%%%%%%%%%%%%%%%%%%%%%%%%%%%%%%%%%%%%%%%%%%%%%%%%%%%
\begin{lem}\label{lem:standideal}
Let $n\in\Z_{\geqq 0}$.
The correspondence $\lex\mapsto \lex^{-1}([1,n])$ gives a surjective map
$$\LE(\cdiag)\to \I_n(\cdiag).$$
Moreover, for each $\lex\in\LE(\cdiag)$,
the restriction $\lex\mapsto \lex|_{\lex^{-1}([1,n])}$
gives a surjective map
$$\LE(\cdiag)\to \LE(\lex^{-1}([1,n])).$$
%gives a surjective map 
\end{lem}
%%%%%%%%%%%%%%%%%%%%%%%%%%%%%%%%%%%%%%%%%%

%%%%%%%%%%%%%%%%%%%%%%%%%%%%%%%%%%%%%%%%%%%%%%%%%%%%%%%%%%
\subsection{Content map and bottom set} 
%%%%%%%%%%%%%%%%%%%%%%%%%%%%%%%%%%%%%%%%%%%%%%%%%%%%%%%%%
Let $\pdiag$ be a periodic diagram of period $\omega$.
Define the \alert{content map} %(or \alert{coloring map}) %(or \alert{residue map})
$$\con:\pdiag \to \Z$$ by $\con(a,b) = b-a.$
Put $\kappa=|\con(\omega)|$.
Let $\cdiag=\pi(\pdiag)$.
Since $\con(x + \omega) = \con(x)-\kappa$, 
the content map $\con$ induces the map
$$\cdiag\to \Z/\kappa\Z,$$
which we denote by the same symbol $\con$.
%(See Figure \ref{content}, where the number in each cell is the content with modulo 8.)
%(the $\con$ plays a role of ``d-complete coloring''.)
It is easy to show the following:
%%%%%%%%%%%%%%%%%%%%%%%%%%%%%%%%%%%%%%%%%%%%%%%%%%%%%%%
\begin{prop}\label{prop:content}
%Let $\cdiag$ be a cylindric diagram.
For $x,y\in \cdiag$, the followings hold:
%\begin{enumerate}[\rm (1)]
%
%\item

\noindent
$\mathrm{(1)}$  If  $\con(x) - \con(y) \equiv 0, \pm1 \mod \kappa$,
then $x$ and $y$ are comparable.
%
%\item Let $[x,y]$ be a chain.
%If $\#[x,y] \leqq \kappa$ then $\con(z)$'s $(z \in [x,y])$ are distinct.
%If $\#[x,y] = \kappa+1$  then $\con(x)-\con(y)=\pm\kappa$.
%
%\item If $[x,y]$ is a $d_k$-interval, then $\con(x)=\con(y)$.
%
%\item 

\noindent
$\mathrm{(2)}$ If $x$ is covered by $y$, then $\con(x)-\con(y)\equiv\pm 1 \mod \kappa$.
\end{prop}
%%%%%%%%%%%%%%%%%%%%%%%%%%%%%%%%%%%%%%%%%%%%%%%%%%%%%%%%%%%

%%%%%%%%%%%%%%%%%%%%%%%%%%%%%
\begin{rem}
By Proposition \ref{prop:skew}  and Proposition \ref{prop:content},
cylindric diagrams are infinite (locally finite)
``$\Z/\kappa\Z$-colored $d$-complete posets'' in the sense of \cite{Str2020, Str2021}.
\fin
\end{rem}
%%%%%%%%%%%%%%%%%%%%%%%%%%%%%%%%%%%%%%%%%%%

\if0
%%%%%%%%%%%%%%%%%%%%%%%%%%
%\begin{df}
For a connected poset $P$ with minimal elements,
the \alert{bottom set} of $P$ is by definition the subset of $P$
consisting of all those elements $y$ satisfying
$$%\Gamma=\{y\in P\mid
x\in P\ \text{is minimal} \implies  [x,y] \text{ is a chain or empty}.$$
%\end{df}
%%%%%%%%%%%%%%%%%%%%%%%%%%%
\fi

Let $i \in \Z/\kappa\Z$.
By Proposition \ref{prop:content} (1),
the inverse image $\con^{-1}(i)$ is non-empty totally ordered subset of $\cdiag$.
Let $\cmin_i$ denote the minimum element in $\con^{-1}(i)$.

\begin{df}
%For a cylindric diagram $\cdiag$,
Define the \alert{bottom set} $\Gamma$ of $\cdiag$ by
$$
\Gamma = \{ \cmin_i \mid i \in \Z/\kappa\Z\}.
$$
%An example of a bottom set appears in Figure \ref{content}.
%(The set of yellowed cells represent the bottom set.)
\end{df}

\if0
%\noindent
%$\mathrm{(2)}$
%Let $
If  $\lm\in \ptn_{\peri}$, then $\per{\lambda}$ is a periodic diagram of period  $\peri$
and $\cyl\lm$ is a cylindric diagram.
Moreover,
 any periodic (resp. cylindric)
  diagram of period $\peri$ is of the form $\per\lambda$
  (resp. $\pi(\per{\lambda})$) for some 
 $\lambda \in \ptn_{\peri}$.
\fi

Figure \ref{content} indicates the periodic diagram $\per\lm$ with $\lm=(5,4,4,2)\in\ptn_{(4,-5)}$.
The number in each cell is the content with modulo 9.
%
%Colored cells form the fundamental domain $\sem\lm$.
%
Yellowed cells %correspond to cells in 
forms the bottom set of $\cyl\lm=\pi(\per\lm)$.
 
\if0
For generalized partitions $\lm=(\lm_1,\dots,\lm_m)$ and $\mu=(\mu_1,\dots,\mu_m)$,
%\in\P_{m,\ell}$,
we write $\lm\supset\mu$ if $\lm_i \geqq \mu_i$ for all $i$ with $1\leqq i \leqq m$.
For $\lm\supset\mu$, 
the set difference $\sem\lm/\sem\mu$ is an ordinary skew diagram associated with $\lm$ and $\mu$,
which is denoted also by $\lm/\mu$.
\fi

\begin{figure}[h!]
\begin{small}
\begin{center}

\begin{tikzpicture}[scale=.45, rotate=270]

%%%%%%%%%%%%%%%%%%%%%%%%%%%%%%%%%%%%%%%%%%%%%%%%%%%%%%%%%%%%%%%%%%%%%%%%%%%%

\fill[yellow!60, xshift=-4cm, yshift=5cm] (3,3) |- (4,1) |- (3,-1) |- (2,0) |- (3,3);
\fill[yellow] (0,2) |- (1,5) |- (3,3) |- (4,1) |- (3,-1) |- (2,0) |- (0,2) --cycle;
\fill[yellow!60, xshift=4cm, yshift=-5cm] (0,2) |- (1,5) |- (2,3) |- (0,2);

\draw[xshift=-4cm, yshift=5cm, ultra thick] (2,3) -- (3,3) |- (4,1) -- (4,-5);
\draw[ultra thick] (0,-7) |- (1,5) |- (3,3) |- (4,1) -- (4,-7);
\draw[xshift=4cm, yshift=-5cm, ultra thick]  (0,0) |- (1,5) |- (2,3);

%%%%%%%%%%%%%%%%%%%%%%%%%%%%%%%%%%%%%%%%%%%%%%%%%%%%%%%%%%%%%%%%%%%%%%%%%%%%

\begin{scope}[xshift=-4cm, yshift=5cm]
\draw[dotted, ultra thick] (1,-10) -- +(0,10);
\draw[black!60] (2,-12) grid (4,1);
\draw[black!60] (2,1) grid (3,3);
%\draw (0,3) grid (1,5);
\end{scope}

\draw[dotted, ultra thick] (2,-8) -- +(0,-3);
\draw (0,-7) grid (4,1);
\draw (0,1) grid (3,3);
\draw (0,3) grid (1,5);

\begin{scope}[xshift=4cm, yshift=-5cm]
%\draw (0,-5) grid (4,1);
\draw[black!60] (0,-2) grid (2,3);
\draw[black!60] (0,3) grid (1,5);
\draw[dotted, ultra thick] (3,-1) -- +(0,3);
\end{scope}

%%%%%%%%%%%%%%%%%%%%%%%%%%%%%%%%%%%%%%%%%%%%%%%%%%%%%%%%%%%%%%%%%%%%%%%%%%%%

\coordinate (q) at (.5,.5);
\foreach \c/\y in {8/-14, 0/-13, 1/-12, 2/-11, 3/-10, 4/-9, 5/-8, 6/-7, 7/-6, 8/-5, 0/-4, 1/-3, 2/-2, 3/-1, 4/0}
\node[black!60] at ($(2,\y)+(q)+(-4,7)$) {$\c$};

\foreach \c/\y in {7/-15, 8/-14, 0/-13, 1/-12, 2/-11, 3/-10, 4/-9, 5/-8, 6/-7, 7/-6, 8/-5, 0/-4, 1/-3}
\node[black!60] at ($(3,\y)+(q)+(-4,8)$) {$\c$};

\foreach \c/\y in {6/-7, 7/-6, 8/-5, 0/-4, 1/-3, 2/-2, 3/-1, 4/0, 5/1, 6/2, 7/3, 8/4}
\node at ($(0,\y)+(q)$) {$\c$};

\foreach \c/\y in {5/-8, 6/-7, 7/-6, 8/-5, 0/-4, 1/-3, 2/-2, 3/-1, 4/0, 5/1}
\node at ($(1,\y)+(q)+(0,1)$) {$\c$};

\foreach \c/\y in {4/-9, 5/-8, 6/-7, 7/-6, 8/-5, 0/-4, 1/-3, 2/-2, 3/-1, 4/0}
\node at ($(2,\y)+(q)+(0,2)$) {$\c$};

\foreach \c/\y in {3/-10, 4/-9, 5/-8, 6/-7, 7/-6, 8/-5, 0/-4, 1/-3}
\node at ($(3,\y)+(q)+(0,3)$) {$\c$};

\foreach \c/\y in {2/-2, 3/-1, 4/0, 5/1, 6/2, 7/3, 8/4}
\node[black!60] at ($(0,\y)+(q)+(4,-5)$) {$\c$};

\foreach \c/\y in {1/-3, 2/-2, 3/-1, 4/0, 5/1}
\node[black!60] at ($(1,\y)+(q)+(4,-4)$) {$\c$};

\end{tikzpicture}

\end{center}
\end{small}
%\caption{Contents on cylindric diagram associated with $\lm=(5,4,4,2) \in \ptn_{4,4}$.}\label{fig1}
\caption{%Contents on periodic diagram associated with
  %$\lm=(5,4,4,2) \in \ptn_{4,4}$.
  %Yellowed cells represent the bottom set.
  %The number in each cell is the content with modulo 8.
  }
  \label{content}
\end{figure}

%%%%%%%%%%%%%%%%%%%%%%%%%%%%%%%%%%%
%\section{Roots, words and weights}

%%%%%%%%%%%%%%%%%%%%%%%%%%%%%%%%%%%%%%%%%%%%%%%%%%%%%%%%%%%%%%%%%%%%%%%%%%%%%%%%
\subsection{Root systems and affine Weyl groups of type $A_{\kappa-1}^{(1)}$}
\label{sec:rootsystem}
%%%%%%%%%%%%%%%%%%%%%%%%%%%%%%%%%%%%%%%%%%%%%%%%%%%%%%%%%%%%%%%%%%%%%%%%%%%%%%%%

Let $\kappa \in \Z_{\geqq2}$. 
In the rest, we often identify $\Z/\kappa\Z$ with $\{0,1,\dots,\kappa-1\}$.
Let $\h$ be a $(\kappa+1)$-dimensional vector space
and choose elements $\al^\vee_i\ (i\in\Z/\kappa\Z)$ and $d$
of $\h$ so that
$$\{\al^\vee_0,\al^\vee_1,\dots,\al^\vee_{\kappa-1},d\}$$ 
forms a basis for $\h$.
Let $\h^*$ be the dual space of $\h$.
Define elements 
%$\al_0,\al_1,\dots,\al_{\kappa-1}$ 
$\al_j\ (j\in\Z/\kappa\Z)$
and $\fdwt_0$ of $\h^*$ by
\begin{align*}
&\bra \al_j,\al^\vee_i \ket = a_{ij},\quad
\bra \fdwt_0,\al^\vee_i\ket = \delta_{i0}
\quad (i,j\in\Z/\kappa\Z)
%(i=0,1,\dots,\kappa-1)
,\\
%\quad
&\bra \al_j,d \ket = \delta_{j0},\quad
\bra \fdwt_0,d \ket =0,
\end{align*}
where $\bra\cdot,\cdot\ket:\h^*\times\h\to\Z$ is the natural pairing
and the integer $a_{ij}$ is defined by
\[
a_{ij}=
\begin{cases}
2 & \text{if $i=j$} \\
-1 & \text{if $i-j= \pm 1$} \\
0 & \text{otherwise}
\end{cases}
%\quad
%\text{for $\kappa\geqq3$}
\]
for $\kappa\geqq3$ and
\[
a_{ij}=
\begin{cases}
2 & \text{if $i=j$} \\
-2 & \text{if $i \neq j$}
\end{cases}
%\quad
%\text{for $\kappa=2$.}
\]
for $\kappa=2$.
Then $\{\al_0,\al_1,\dots,\al_{\kappa-1},\fdwt_0\}$ forms a basis for $\h^*$.
Define $\fdwt_i \in \h^*$ ($i =1,2,\dots,\kappa-1$) by
$$\bra\fdwt_i,\al^\vee_j\ket = \delta_{ij},\quad \bra\fdwt_i,d\ket=0
\quad (j\in\Z/\kappa\Z).
%(j=0,1,\dots,\kappa-1)
$$
The weights $\fdwt_0,\fdwt_1,\dots,\fdwt_{\kappa-1}$ are called fundamental weights.
Put $\delta = \al_0+\al_1+\cdots+\al_{\kappa-1}$ (resp. $\delta^\vee=\al^\vee_0+\al^\vee_1+\cdots+\al^\vee_{\kappa-1}$),
which is called the null root (resp. the null coroot).

For $i \in \Z/\kappa\Z
%[0,\kappa-1]
$,
define the simple reflection $s_i \in GL(\h^*)$  by
\begin{align*}
&s_i(\zeta) = \zeta - \bra \zeta,\al^\vee_i \ket \al_i \quad (\zeta \in \h^*).
%\\
%&s_i(\eta^\vee) = \eta^\vee - \bra \al_i,\eta^\vee \ket \al^\vee_i \quad (\eta^\vee \in \h).
\end{align*}
%%
%It is called a simple reflection.
Define the {\it affine Weyl group} $W$ of type $A_{\kappa-1}^{(1)}$
% $W$ be
as the subgroup of $GL(\h^*)$
generated by simple reflections: % $\{s_i \mid i \in [0,\kappa-1]\}$:
$$W=\langle s_i\mid i\in\Z/\kappa\Z\rangle.$$
%The group $W$ is called 
%

The following is well-known:
%We remark that the group $W$ has the following relations:
%%%%%%%%%%%%%%%%%%%%%%%%%%%%%%%%%%%%%%%%%%%%%%%%%%%%%%%%%%%%%%%%%%%%%
\begin{prop}
The group $W$ has the following fundamental relations$:$
\begin{align}
&s_i^2=1, \label{eq:rel1} \\
&s_is_j=s_js_i \quad (i-j \neq 0,\pm1), \label{eq:rel2} \\
&s_is_{i+1}s_i=s_{i+1}s_is_{i+1}. \label{eq:rel3}
\end{align}
\end{prop}
%%%%%%%%%%%%%%%%%%%%%%%%%%%%%%%%%%%%%%%%%%%%%%%%%%%%%%%%%%%%%%%%%%%%
%
%%%%%%%%%%%%%%%%%%%%%%%%%%%%%%%%%%%%%
For $w  \in W$, 
we define the {\it length}  $\ell(w)$ of $w$
as the smallest $r$ for which an expression 
(or a word) 
$$w = s_{i_1}s_{i_2} \cdots s_{i_r} \in W\ \ (i_j \in \Z/\kappa\Z)$$
exists.
%and call the reduced expression.
An expression $w=s_{i_1}s_{i_2} \cdots s_{i_r}$ is said to be reduced if
 $\ell(w)=r$.
%By $\RE(w)$, we denote the set of reduced expressions for $w$.

Define the action of $W$ on $\h$ by
$$
s_i(h) = h - \bra \al_i,h \ket \al^\vee_i \quad (h \in \h).
$$

We put 
\begin{align*}
\Pi &= \{\al_0,\al_1,\dots,\al_{\kappa-1}\},\quad
\Pi^\vee = \{\al^\vee_0,\al^\vee_1,\dots,\al^\vee_{\kappa-1}\},\\
Q&=\left\{\sum_{i\in\Z/\kappa\Z}c_i \al_i\ \middle|\ c_i\in\Z\right\},\ \ 
Q_+=\left\{\sum_{i\in\Z/\kappa\Z}c_i \al_i\ \middle|\ c_i\in\Z_{\geqq 0}\right\}.
\end{align*}
The set $\Pi$ (resp. $\Pi^\vee$) is called
the set of  simple roots (resp. the set of simple coroots),
and $Q$ is called the root lattice.
%%%%%%%%%%%%%%%%
%%%%%%%%%%%%%%%%%%%%%%%%
Put
\begin{align*}
&R=W\Pi\subset\h^*, \quad R^\vee=W\Pi^\vee\subset\h. %\\
%&\tilde{R}=R\sqcup\Z\delta \quad (\tilde{R}^\vee=R^\vee\sqcup\Z\delta^\vee).
\end{align*}
Then %$\tilde{R}$ is the affine root system, %(resp. affine coroot system),
$R$ (resp. $R^\vee$) is the set of real roots (resp. coroots) 
and $R\sqcup\Z\delta$ is the affine root system. %(resp. imaginary roots).
Define the set $R_+$ of positive (real) roots and the set $R_-$ of negative (real) roots by
$$R_+ =R\cap Q_+= \left\{ \sum_{i=0}^{\kappa-1} c_i \al_i \in R \mid c_i \in \Z_{\geqq0}\right\}, \quad
 R_-= \left\{ \sum_{i=0}^{\kappa-1} c_i \al_i \in R \mid c_i \in \Z_{\leqq0}\right\}.$$
For $\beta=\sum_{i=0}^{\kappa-1} k_i\al_i\in R,$
define $\beta^\vee=\sum_{i=0}^{\kappa-1} k_i\al^\vee_i \in R^\vee.$
Then the correspondence $\beta\mapsto \beta^\vee$ gives a bijection
$R\to R^\vee$.
%Similarly,
Define the set of positive (resp. negative) coroots $R^\vee_+$ (resp. $R^\vee_-$)
as the image of $R_+$ (resp. $R_-$) by this bijection.

%%%%%%%%%%%%%%%%%%%%%%%%%%%%%%%%%%%%%%%%%%%%%%%%%%%%%%%%%%%
For $i,j\in\Z$ with $i<j$, we define
$$\al_{ij}=\sum_{i\leqq k\leqq j-1}\al_{\bar{k}},$$
where $\bar{k}=k\!\!\mod\kappa\Z\in\Z/\kappa\Z$.
%%%%%%%%%%%%%%%%%%%%%%%%%%%%%%%%%%%%%%%%%%%%%%%%%%%%%%%%%%%%%%%%
The followings are well-known:
\begin{align}
R_+&=\{\al_{ij}\mid i<j,\ j-i\notin \kappa\Z\}\label{eq:R+}
\\
&=\{\al_{ij}+k\delta\mid
1\leqq i<j\leqq\kappa,\ k\geqq 0\}\sqcup
\{-\al_{ij}+k\delta\mid
1\leqq i<j\leqq\kappa,\ k\geqq 1\},\label{eq:R+2}\\
R_-&=-R_+,\ \
%\{-\al_{ij}-k\delta\mid
%1\leqq i<j\leqq\kappa,\ k\geqq 0\}\sqcup
%\{\al_{ij}-k\delta\mid
%1\leqq i<j\leqq\kappa,\ k\geqq 1\},\label{eq:R-}\\
%&=\{-\al_{ij}\mid i<j,\ j-i\notin \kappa\Z\}\\
R=R_+\sqcup R_-.\notag
%\label{eq:R}
%=\{\al_{ij}+k\delta\mid
%1\leqq i\neq j\leqq \kappa,\ k\in\Z\}.\label{eq:R}
\end{align}
%where $\al_{ij}=\al_i+\al_{i+1}+\cdots+\al_{j-1}$ for $i<j$
%and $\al_{ij}=-\al_{ji}$ for $i>j$.
From the description of $R$ above, the following two lemmas  follow
easily  and they will be used later: 
%%%%%%%%%%%%%%%%%%%%%%%%%%%%%%%%%%%%%%%%%%%%%%%%%%%%%%%%%%%%%%%%%%%%%%%%%
\begin{lem}\label{lem:addnull}
%It holds that
%
If $\al\in R$, then $\al+k\delta\in R$ for all $k\in\Z$.
\end{lem}
%%%%%%%%%%%%%%%%%%%%%%%%%%%%%%%%%%%%%%%%%%%%%%%%%%%%%%%%%%%%%%%%%%%%%%%%%

%%%%%%%%%%%%%%%%%%%%%%%%%%%%%%%%%%%%%%%%%%%%%%%%%%%%%%%%%%%%%%%%%%%
\begin{lem}\label{lem:moddelta}
Let $\alpha \in R \sqcup \Z\delta$ and $\beta \in R$.
%Then $\bra \alpha, \beta \ket = 2$ if and only if $\alpha \equiv \beta \mod \delta$.
Then $\bra\alpha,\beta^\vee\ket=2$ if and only if $\alpha \equiv \beta \mod \delta$.
\end{lem}
%%%%%%%%%%%%%%%%%%%%%%%%%%%%%%%%%%%%%%%%%

%%%%%%%%%%%%%%%%%%%%%%%%%%%%%%%%%%%%%%%%%%%%%%%%%%%%%%%%%
\section{Hooks in cylindric diagrams}\label{sec:hooks}
%%%%%%%%%%%%%%%%%%%%%%%%%%%%%%%%%%%%%%%%%%%%%%%%%%
\subsection{Colored hook length}\label{ss:hook}
%%%%%%%%%%%%%%%%%%%%%%%%%%%%%%%%%%%%%%
% Then, for the reduced semi-infinite word 
%$w_{\cdiag,\lex}$, we have

%Moreover, by Proposition \ref{pr;hk=hkr}, we have
%In particular, t will be
%We will denote $R(w_{\cdiag,\lex})$ just by $R(w_\cdiag)$ below.

In this section, we will introduce colored hook length,
which is 
a key ingredient in this paper.

Fix $\kappa, m,\ell\in\Z_{\geqq1}$ with $\kappa=m+\ell$
and let $\cdiag$ be a cylindric diagram in $\cyli_{(m,-\ell)}$.

In the rest of this paper, we use the following notations:
$$\al(x)=\al_{\con(x)},\ \ s(x)=s_{\con(x)}\ \ \text{for } x\in \cdiag.$$
%%%%%%%%%%%%%%%%%%%%%%%%
\begin{df}\label{df:coloredhooklength}
For $x\in \cdiag$, put
\begin{align*}
\arm(x)&=\{x+(0,k)\in\cdiag\mid k\in\Z_{\geqq1}\},\\
\leg(x)&=\{x+(k,0)\in\cdiag\mid k\in\Z_{\geqq1}\},
\end{align*}
and define
%The "multiset" 
%\begin{df}
\begin{equation*}
\hk(x)=\al(x)+\sum_{y\in\arm(x)}\al(y)+\sum_{y\in\leg(x)}\al(y).
\end{equation*}
%where $\al(x)=\al_{\con(x)}$.
We call $\hk(x)$ the {\it colored hook length} at $x$.
\end{df}
%%%%%%%%%%%%%%%%%%%%%%%%%%%%%%%%%%%%%%%%%%%%%%%%%%%%%%%%%%%%%%%%%%%%%%%%%%%%%%%%%%%%%
% as before.

%%%%%%%%%%%%%%%%%%%%%%%%%%%%%%%%%%%%%%%%%%%%%%%%%%%%%%%%%%%%%%%%%%%%%%%%%%%%%%%%%
\begin{figure}[h!]
\begin{center}
\begin{tikzpicture}[scale=.45, rotate=270]

%%%%%%%%%%%%%%%%%%%%%%%%%%%%%%%%%%%%%%%%%%%%%%%%%%%%%%%%%%%%%%%%%%%%%%%%%%%%
\fill[pink] (0,-10) |- (1,5) |- (3,3) |- (4,1) -- (4,-10) -- cycle;

\coordinate (c) at (.5,.5);
\coordinate (x) at ($(1,-5)+(c)$); 
\fill[red] ($(x)+(c)$) rectangle +(-1,-1);
\fill[yellow] ($(x)+(c)$) rectangle +(-1,7);
\fill[green] ($(x)+(c)$) rectangle +(6,-1);
%\node at (x) {$x$};
%\node at ($(x)!.5!($(x)+(0,7.5)$)$) {$\arm(x)$};
%\node[rotate=-90] at ($(x)!.5!($(x)+(6.5,0)$)$) {$\leg(x)$};

\draw[xshift=-4cm, yshift=5cm, ultra thick] (2,3) -- (3,3) |- (4,1) -- (4,-5);
%\draw[xshift=-4cm, yshift=5cm, ultra thick] (0,-15) |- (1,5) |- (3,3) |- (4,1) -- (4,-5);
\draw[ultra thick] (0,-10) |- (1,5) |- (3,3) |- (4,1) -- (4,-5);
\draw[xshift=4cm, yshift=-5cm, ultra thick] (0,-5) |- (1,5) |- (3,3) |- (4,1) -- (4,-5);
\draw[xshift=8cm, yshift=-10cm, ultra thick] (0,0) |- (1,5) |- (2,3);

%%%%%%%%%%%%%%%%%%%%%%%%%%%%%%%%%%%%%%%%%%%%%%%%%%%%%%%%%%%%%%%%%%%%%%%%%%%%

\begin{scope}[xshift=-4cm, yshift=5cm]
\draw[dotted, ultra thick] (1,-14) -- +(0,16);
\draw (2,-15) grid (4,1);
\draw (2,1) grid (3,3);
%\draw (0,3) grid (1,5);
\end{scope}

\if0
\begin{scope}[xshift=-4cm, yshift=5cm]
\draw[dotted, ultra thick] (2,-16) -- +(0,-3); 
\draw (0,-15) grid (4,1);
\draw (0,1) grid (3,3);
\draw (0,3) grid (1,5);
\end{scope}
\fi

\draw[dotted, ultra thick] (2,-11) -- +(0,-3);
\draw (0,-10) grid (4,1);
\draw (0,1) grid (3,3);
\draw (0,3) grid (1,5);

\begin{scope}[xshift=4cm, yshift=-5cm]
\draw[dotted, ultra thick] (2,-6) -- +(0,-3);
\draw (0,-5) grid (4,1);
\draw (0,1) grid (3,3);
\draw (0,3) grid (1,5);
\end{scope}

\begin{scope}[xshift=8cm, yshift=-10cm]
%\draw (0,-5) grid (4,1);
\draw (0,0) grid (2,3);
\draw (0,3) grid (1,5);
\draw[dotted, ultra thick] (3,0) -- +(0,3);
\end{scope}

%\draw[dashed] (0,5) -- +(4,0) -- +(4,-5);
%\draw (0,5) to [out=60, in=120] +(4,0) to [out=-60, in=60] +(0,-5);
%\draw (0,5) sin (2,6) cos  (4,5);
%\node at (2,6) [fill=white,inner sep=5pt] {4};
%\draw (4,5) cos (5,2.5) sin (4,0);
%\node at (5,2.5) [fill=white,inner sep=5pt] {5};

%\draw[thick, ->] (3,4) -- +(4,-5);
%\node at (7,1.5) {$(4,-5)$};

\node at (x) {$x$};
\coordinate (a) at ($(x)!.5!($(x)+(0,7.5)$)$);
\fill[white] ($(a)+(-.4,-1.5)$) rectangle ($(a)-(-.4,-1.5)$);
\node at (a) {$\arm(x)$};
%\node at ($(x)!.5!($(x)+(0,7.5)$)$) {$\arm(x)$};
\coordinate (l) at ($(x)!.5!($(x)+(6.5,0)$)$);
\fill[white] ($(l)+(-1.5,-.4)$) rectangle ($(l)-(-1.5,-.4)$);
\node[rotate=-90] at (l) {$\leg(x)$};
%\node[rotate=-90] at ($(x)!.5!($(x)+(6.5,0)$)$) {$\leg(x)$};

\end{tikzpicture}
\end{center}
\caption{The sets $\arm(x)$ and $\leg(x)$ for $x$ in the cylindric diagram. 
%The colored hook length at $x$ is $\hk(x)%= \al_0 + 2\al_1 + 2\al_2 + 2\al_3 + 2\al_4 + 2\al_5 + \al_6 + \al_7 + \al_8 
%= \al_1 + \al_2 + \al_3 + \al_4 + \al_5 + \delta$
}\label{fig:coloredhook}
\end{figure}
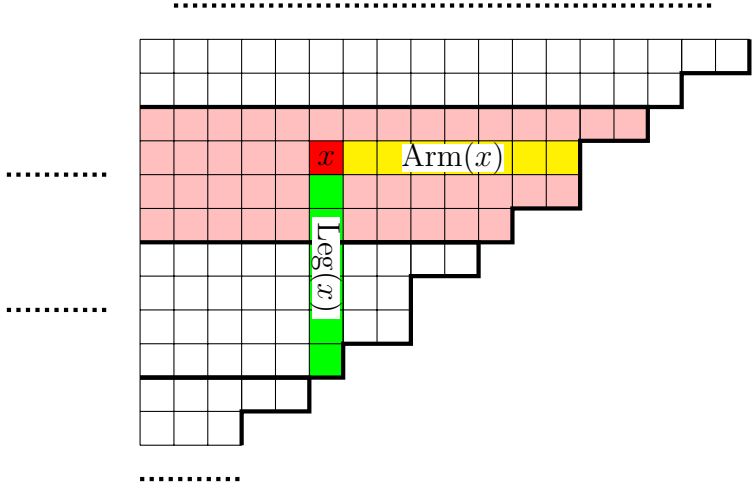

For $x\in \mathcal{C}_{(m,-\ell)}\setminus \cdiag$, we set $\hk(x)=0$ for convenience. 
%\end{df}
It is easy to see that for $x\in\cdiag$
$$\hk(x-(0,\ell))=\hk(x-(m,0))=\hk(x)+\delta$$
and
\begin{equation}\label{eq:hk=al}
\hk(x)=\al_{ij}\ \text{for some integers }i<j.
\end{equation}
% and thus $\hk(x)\in R_+$ by \eqref{eq:R+}.

%%%%%%%%%%%%%%%%%%%%%%%%%%%%%%%%%%
\begin{ex} (See Figure \ref{fig:coloredhook}.)\ 
Let $\peri=(4,-5)$. Then $\lm=(5,3,3,1) \in \ptn_\peri$.
For a cell $x=\pi (2,-4)\in \cyl\lm$, we have $\con (x)=3 +9\Z\in \Z/9\Z$.
The colored hook length at $x$ is 
\begin{align*}
\hk(x)
&=\al_{-6}+(\al_{-5}+\al_{-4}+\al_{-3}+\al_{-2}+\al_{-1}+\al_0+\al_1)\\
&+(\al_{-7}+\al_{-8}+\al_{-9}+\al_{-10}+\al_{-11}+\al_{-12})\\
&=\al_3+(\al_4+\al_5+\al_6+\al_7+\al_8+\al_0+\al_1)+(\al_2+\al_1+\al_0+\al_8+\al_7+\al_6)\\
&=\delta +\al_0+\al_1+\al_6+\al_7+\al_8,
\end{align*}
which can be expressed as 
$\hk(x)=\al_{-12, 2}$.
%(Note thta  $\con(8,-4)=-12$ and
%$\con(2,3)=1$ in the periodic diagram $\per\lm$.)
\end{ex}
%%%%%%%%%%%%%%%%%%%%%%%%%%%%%%%%%
%
%
%%%%%%%%%%%%%%%%%%%%%
\begin{rem}
(1) For $x\in\cdiag$, the ``multiset" $\mathrm{H}(x):=\{x\}\sqcup 
\arm(x)\sqcup \leg(x)$ is a cylindric 
analogue of  the hook at $x$.
%and $\hk(x)$ is thought as a colored version of hook length.

\smallskip\noindent
(2) A conjectural  hook formula 
concerning the number of standard tableaux on cylindric skew diagrams is proposed
in \cite{SuzToy2022}, where the hook length at $x\in \cdiag$ is given by 
 $|H(x)|=|\arm(x)|+|\leg(x)|+1$.
\end{rem}
%%%%%%%%%%%%%%%%%%%%%%%%%%%%%%

%Put $Q_+=\{\sum_{i\in\Z/{\kappa\Z}} c_i\al_i\mid c_i\in\Z_{\geqq 0}\}$.
For $\al\in Q_+$,
%$x \in \cdiag$,
define
\[ \nmult(\al)=\max\{k\in\Z \mid\al-k\delta \in Q_+\}. \]
%%%%%%%%%%%%%%%%%%%%%%%%%%%%%%%%%%%%%%%%%%%%%
\begin{lem}\label{lem:N=N}
%Let $\cdiag$ is a cylindric digagram in $\cyli_{(m,-\ell)}$.
For $x\in\cdiag$, it holds that
$$\nmult(\hk(x))=\max \{k\in\Z\mid
x+k(0,\ell)\in\cdiag\}.$$
\end{lem}
%%%%%%%%%%%%%%%%%%%%%%%%%%%%%%%%%%%%%%%%%%%%%%%%%%%%%%%%%%%%%%%%%
\begin{pf}
We put $N(x)=\max \{k\in\Z\mid
x+k (0,\ell)\in\cdiag\}$ and will show $N(\hk(x))=N(x)$.

Let $k\in\Z_{\geqq0}$.
Suppose that $x+k(0,\ell)\in\cdiag$.
Then, $\hk(x)-k\delta=\hk(x+k(0,\ell))\in Q_+$ and thus
$\nmult(\hk(x))\geqq N(x)$.

Suppose that $x-k(0,\ell)\notin\cdiag$.
Noting that $x-k(0,\ell)=x+k(m,0)$, we have
$$|\arm (x)|\leqq k\ell-1,\ \ \ |\leg(x)| \leqq km-1.$$
Thus we have $|\{x\}\cup\arm(x)\cup\leg(x)|\leqq k(\ell+m)-1$
and hence $\hk(x)-k\delta\notin Q_+$.
This means $ N(\hk(x))\leqq N(x)$.
\qed
\end{pf}
%%%%%%%%%%%%%%%%%%%%%%%%%%%%%%%%%%%%%%%%%%%%%%%%%%%%%%%%
Let $\btm=\{b_i\mid i\in\Z/\kappa\Z\}$ be the bottom set of $\cdiag$,
where $b_i$ is the minimum element of $\con^{-1}(i)$ as before.

%Let $Q$ denote the root lattice: $Q=\{\sum_{i\in\Z/{\kappa\Z}} c_i\al_i\mid c_i\in\Z\}$.
For $\al=\sum_{i\in\Z/{\kappa\Z}} c_i\al_i
\in Q_+
%R\sqcup \Z\delta
$, define its {\it support}
by
$$\Supp(\al)=\{b_i \mid c_i >0\ (i\in\Z/\kappa\Z)\}\subset \btm.$$
For example, we have $\Supp(\delta)=\btm$.
%For example, we have $\Supp(\delta)=\btm$.
%Note that $|\btmmax|=|\btmmin|$.
%
Let $x\in\cdiag$ with $\nmult(\hk(x))=0$. 
Then $\Supp(\hk(x))$ is a non-empty, proper and connected subset of $\btm$.
\if0
Moreover, we have
\begin{equation}\label{eq:hkx0}
\hk(x)=\sum_{y\in\Supp(\hk(x))}\al(y).
\end{equation}
\fi
%%%%%%%%%%%%%%%%%%%%%%%%%%%%%%%%%%%%%%%%%%%%%
\begin{lem}\label{lem:hkinR}
%Let $\cdiag$ is a cylindric digagram in $\cyli_{(m,-\ell)}$.
Let $x\in\cdiag$. Then $\hk(x)\in R_+$.
\end{lem}
%%%%%%%%%%%%%%%%%%%%%%%%%%%%%%%%%%%%%%%%%%%%%%%%%%%%%%%%%%%%%%%%%
\begin{pf}
By \eqref{eq:R+} and \eqref{eq:hk=al},
it is enough to show that $\hk(x)\notin \Z\delta$.

Put $k=\nmult(\hk(x))$ and $x_0=x+k(0,\ell)$.
Then $x_0\in\cdiag$ by Lemma \ref{lem:N=N} and $\nmult(\hk(x_0))=0$.
Since
%$\hk(x_0)=\sum_{y\in\Supp(\hk(x_0)}\al_{\con(y)}$ and moreover 
$\emptyset\neq\Supp(\hk(x_0))\subsetneqq\btm$,
we have $\hk(x_0)\notin \Z\delta$ and thus $\hk(x)=\hk(x_0)+k\delta\notin \Z\delta.$
\qed
\end{pf}

Let $\btm_{\mathrm{max}}$ (resp. $\btm_{\mathrm{min}}$) denote the set of maximal (resp. minimal)
elements in $\btm.$
Note that $|\btm_{\mathrm{max}}|=|\btm_{\mathrm{min}}|$.
One can easily see the following lemma. (See the figure below.)
%%%%%%%%%%%%%%%%%%%%%%%%%%%%%%%%%%%%%%%%%%%%%
\begin{lem}\label{lem:supportofhk}
Let $\al\in R_+$ with $\nmult(\al)=0$. Then
$\al=\hk(x)$ for some $x\in\cdiag$ if and only if
%$x\in \cdiag$ with $\nmult(x)=0$. Then
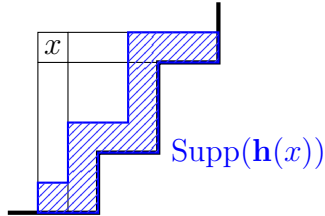
\begin{align*}
%\label{eq:max=min}
&|\Supp (\al)\cap \btmmax|+1=|\Supp (\al) \cap \btmmin|.
%&|\Supp (\delta-\hk(x))\cap \btmmax|=|\Supp (\delta-\hk(x)) \cap \btmmin|+1.
\end{align*}
\end{lem}
%%%%%%%%%%%%%%%%%%%%%%%%%%%%%%%%%%%%%%%%%%%%%%%%%%%%%%%%%%%%%%%%%
%\newcommand\im{{\mathrm{Im}}}
%%%%%%%%%%%%%%%%%%%%%%%%%%%%%%%%%%%%%%%%%%%%%%%%%%%%%%%%%%%

\begin{center}
\begin{tikzpicture}[scale=.4, rotate=-90]

%\draw[help lines] (1,1) grid (8,8);

\draw[ultra thick] (1,8) -| (3,6) -| (6,4) -| (8,1);

\coordinate (y) at (2,2) node at ($(y)+(.5,.5)$) {$x$};
\draw (y) rectangle +(1,1);
\draw ($(y)+(0,1)$) -- +(0,5);
\draw ($(y)+(1,1)$) -- +(0,5);
\draw ($(y)+(1,0)$) -- +(5,0);
\draw ($(y)+(1,1)$) -- +(5,0);
\draw[thick, blue] (8,2) -| (7,3) -| (5,5) -| (2,8) -| (3,6) -| (6,4) -| (8,2) -- cycle;
\fill[pattern=north east lines, pattern color=blue!80] (8,2) -| (7,3) -| (5,5) 
-| (2,8) -| (3,6) -| (6,4) -| (8,2) -- cycle;

\node[blue] at (6,9) {$\Supp (\hk(x))$};

\end{tikzpicture}
\end{center}

%%%%%%%%%%%%%%%%%%%%%%%%%%%%%%%%
\subsection{Predominant weights and hooks}\label{sec:predom}
%%%%%%%%%%%%%%%%%%%%%%%%%%%%%%%%%%%%%%%%%%%%%%%%%%%%%%%%%%%%%%
\begin{df}\label{df:predom}
We define $\predom_{\cdiag}\in \h^*$ by
\begin{gather}\label{eq:opo}
\predom_{\cdiag} = \sum_{i=0}^{\kappa-1} a_i \fdwt_i, \\
\hbox{where }
a_i = \begin{cases}
      1 & \hbox{if }\cmin_i\in \btm_{\mathrm{max}} \\
      -1 & \hbox{if }\cmin_i\in \btm_{\mathrm{min}} \\
      0 & \hbox{otherwise}.
      \end{cases} \notag
\end{gather}
\end{df}
%%%%%%%%%%%%%%%%%%%%%%%%%%%%%%%%%%%%%%%%%%%%%%%%%%%%%%%%%%%

Note that maximal and minimal elements are lined up alternatively in $\btm$.
This implies that the weight $\predom_\cdiag$ is predominant, namely,
$\bra \predom_\cdiag,\al^\vee \ket\geqq -1$ for all $\al^\vee\in R^\vee_+$.
%
\if0
Note that 
$$ \bra \predom_\cdiag,\al^\vee\ket =-1 \iff
%$ if and only if
|\Supp (\al)\cap \btmmax|-|\Supp (\al) \cap \btmmin|=-1.$$
for $\al\in R_+$ with $N(\al)=0$.
\fi
Define 
$$\DS(\predom_\cdiag)=\{\al\in R_+\mid \bra\predom_\cdiag,\al^\vee\ket=-1\}.$$
%%%%%%%%%%%%%%%%%%%%%%%%%%%%%%%%%%%%%%%%%%%%%%%%
\begin{thm}\label{thm:hkinR+}
The correspondence $x\mapsto \hk(x)$
gives a bijection $$\hk:\cdiag\to \DS(\predom_\cdiag).$$
%Moreover it holds that
%$$\im \hk=\{\al\in R_+\mid \bra\predom_\cdiag,\al^\vee\ket=-1\}.$$
\end{thm}
%%%%%%%%%%%%%%%%%%%%%%%%%%%%%%%%%%%%%%%%%%%%
\begin{pf}
First we will show that $\hk(\cdiag)=\DS(\predom_\cdiag)$.
Let $\al\in R_+$ and
put $\bar\al=\al-\nmult(\al)\delta$.
It follows from Lemma  \ref{lem:N=N}, 
$$\al\in \hk(\cdiag) \Leftrightarrow
\bar\al\in \hk(\cdiag).$$
On the other hand, as $\bra \predom_\cdiag,\delta^\vee\ket =0$, it holds that
$$\al\in \DS(\predom_\cdiag) \Leftrightarrow
\bar\al\in \DS(\predom_\cdiag).$$
Now we have $\hk(\cdiag)=\DS(\predom_\cdiag)$ by Lemma \ref{lem:supportofhk}.
%%%%%%%%%%%%%%%%%%%%%%%%%%%%%%%%%
\if0
Let $x\in\cdiag$.
We have $\hk(x)\in R_+$  by Lemma \ref{lem:hkinR}.
%Let $x\in\cdiag$ with $N(\hk(x))=0$. Then  we have
%$\hk(x)\in\DS(\predom_\cdiag)$ by Lemma \ref{lem:supportofhk}.
Let $x\in\cdiag$ and put $k=N(\hk(x))$. 
Put $x_0=x+k(0,\ell)$. Then  $N(\hk(x_0))=0$
and $\hk(x_0)\in \DS(\predom_\cdiag)$ by Lemma  \ref{lem:N=N}.
Noting $\bra \predom_\cdiag,\delta^\vee\ket =0$,
we have 
$$\bra \predom_\cdiag,\hk(x)^\vee\ket =\bra \predom_\cdiag,\hk(x_0)^\vee\ket
%\bra \predom_\cdiag,(\hk(x_0)+N\delta)^\vee\ket
=-1$$
By  Lemma \ref{lem:supportofhk}.
Thus $\hk(\cdiag)\subset \DS(\predom_\cdiag)$.

Let $\al\in \DS(\predom_\cdiag)$.
% and suppose that $\nmult(\al)=0$.
%Then 
%
%For a general $\al\in \DS(\predom_\cdiag)$, 
Put 
$\beta=\al-\nmult(\al)\delta$. Then $N(\beta)=0$ and
$\beta\in \DS(\predom_\cdiag)$. 
By Lemma \ref{lem:supportofhk}, $\beta=\hk(x)$ for some $x\in\cdiag$, 
for which we have
$\al=\hk(x+\nmult(\al)(0,\ell))$.
Therefore $\DS(\predom_\cdiag)\subset \hk(\cdiag)$.
\fi
%%%%%%%%%%%%%%%%%%%%%%%%%%%%%%%%%%%%%%%%%%%%%%

We will show the injectivity.
Suppose that $\hk(x)=\hk(y)$.
Then $N(\hk(x))=N(\hk(y))$. %which we denote by $N$ below.
Put $x_0=x+N(\hk(x))(0,\ell)$, $y_0=y+N(\hk(y))(0,\ell)$.
Then we have $N(\hk(x_0))=N(\hk(y_0))=0$ and thus $\hk(x_0)=\hk(y_0)$.
Now we have $\Supp(\hk(x_0))=\Supp(\hk(y_0))$
%by \eqref{eq:hkx0} 
and this imples  $x_0=y_0$ and hence $x=y$.
%$$\hk(x_0)=\sum_{z\in\Supp(\hk(x_0))}\al(z),\ \ \hk(y_0)=\sum_{z\in\Supp(\hk(y_0))}\al(z).$$
%Now $\hk(x_0)=\hk(y_0)$ implies $x_0=y_0$ and hence $x=y$.
%$\hk: \cdiag \to \DS(\predom_\cdiag)$.
%is  surjective.
%We put $\DS(\predom_\cdiag)=\{\al\in R_+\mid \bra\predom_\cdiag,\al^\vee\ket=-1\}$.
\qed
\end{pf}

%%%%%%%%%%%%%%%%%%%%%%%%%%%%%%%%%%%%%%%%%%%%%%%%%%%%%%%%%%%%%%%
\subsection{Weyl group elements and their inversion sets}%{Weyl group elements associated with standard tableaux and inversion sets}
\label{sec:inversion}
%%%%%%%%%%%%%%%%%%%%%%%%%%%%%%%%%%%%%%%%%%%%%%%%%%%%%%%%%%%%%%

The following proposition gives an alternative expression for $\hk(x)$.
%%%%%%%%%%%%%%%%%%%%%%%%%%%%%%%%%%%%%%%%%%%%%%%%%%%%%%%%%%%
\begin{prop}\label{prop:hk=hkr}
For any $x\in\cdiag$ and  $\lex\in \LE(\cdiag)$, it holds that 
\begin{equation}\label{eq;hk}
\hk(x)=s(\lex^{-1}(1))s(\lex^{-1}(2))\cdots s(\lex^{-1}(n-1)) \al(\lex^{-1}(n)),
%s(\lex^{-1}(1))s(\lex^{-1}(2))\cdots s(\lex^{-1}({k-1}))\al(x).
\end{equation}
where $n=\lex(x)$.
%In particular, $\hkr(x)$ is independent of $\lex$.
%where $k=\lex(x)$.
\end{prop}
%%%%%%%%%%%%%%%%%%%%%%%%%%%%%%%%%%%%%%%%%%%%%

The proof of Proposition~\ref{prop:hk=hkr} will be given in the next section.
In the rest of this section, we will see some consequences of the proposition.

Let $\lex\in\LE(\cdiag)$. For $n\in\Z_{\geqq1}$,
we define an element $w_{\cdiag,\lex}[n]$ of $W$ by
\begin{equation}
\label{eq:w[n]}
w_{\cdiag,\lex}[n]=s(\lex^{-1}(1))s(\lex^{-1}(2))\cdots s(\lex^{-1}(n)),
\end{equation}
and we set $w_{\cdiag,\lex}[0]=e$.
%%%%%%%%%%%%%%%%%%%%%%

%%%%%%%%%%%%%%%%%%%%%%%%%%%%%%%%%%%%%%%%%%%
\begin{ex}
Let $\lm=(5,4)$ and $\peri=(2,-3).$
For $\lex\in \LE(\cyl\lm)$ displayed in figure \ref{fig:w[6]}, we have
$w_{\cyl\lm,\lex}[6]=s_4s_2s_1s_3s_0s_2$.
\end{ex}
%%%%%%%%%%%%%%%%%%%%%%%%%%%%%%%%%%%%%%%%%%%%%%%%%%%%%

%%%%%%%%%%%%%%%%%%%%%%%%%%%%%%%%%%%%%%%%%%%%%%%%%%%%%%%%%%
\begin{figure}[h!]
\begin{center}
\begin{tikzpicture}[scale=.7,rotate=-90]

\foreach \m in {0,1,...,4}{
\draw[dashed] (3.5,\m+2.5) -- +(-3,-3);
\node at (4,\m+3) {$\m$};
}

\if0
\begin{scope}[black!60,xshift=-2cm,yshift=3cm]
\draw (1,-3) grid (2,5);
\draw (2,-3) grid (3,4);
\coordinate (c) at (.5,.5);
\foreach \n/\a/\b in {1/1/4, 2/2/3, 3/2/2, 4/1/3, 5/2/1, 6/1/2}
\node at ($(\a,\b)+(c)$) {$\n$};
\draw[dotted, ultra thick] (2,-3.5) -- +(0,-1);
\end{scope}
\fi

\begin{scope}[black!60,xshift=2cm,yshift=-3cm]
\draw (1,-3) grid (2,5);
\draw (2,-3) grid (3,4);
\coordinate (c) at (.5,.5);
\foreach \n/\a/\b in {1/1/4, 2/2/3, 3/2/2, 4/1/3, 5/2/1, 6/1/2, 7/1/1, 8/2/0, 9/1/0, 10/2/-1, 11/1/-1, 12/2/-2, 13/1/-2, 14/2/-3, 15/1/-3}
\node at ($(\a,\b)+(c)$) {$\n$};
\draw[dotted, ultra thick] (2,-3.5) -- +(0,-1);
\end{scope}

\draw (1,-6) grid (2,5);
\draw (2,-6) grid (3,4);
\coordinate (c) at (.5,.5);
\foreach \n/\a/\b in {1/1/4, 2/2/3, 3/2/2, 4/1/3, 5/2/1, 6/1/2, 7/1/1, 8/2/0, 9/1/0, 10/2/-1, 11/1/-1, 12/2/-2, 13/1/-2, 14/2/-3, 15/1/-3, 16/2/-4, 17/1/-4, 18/2/-5, 19/1/-5, 20/2/-6, 21/1/-6}{
\fill[white!0] ($(\a,\b)+(c)$) circle (.3);
\node at ($(\a,\b)+(c)$) {$\n$};
}
\draw[dotted, ultra thick] (2,-6.5) -- +(0,-1);
\end{tikzpicture}
\end{center}
\caption{}
\label{fig:w[6]}
\end{figure}

%%%%%%%%%%%%%%%%%%%%%%%%%%%%%%%%%%%%%%%%%%%%%%%%%%%%%%%%%%%
\begin{prop}\label{prop:reduced}
The expression \eqref{eq:w[n]}
%$$w_{\cdiag,\lex}[n]=s(\lex^{-1}(1))s(\lex^{-1}(2))\cdots s(\lex^{-1}(n))$$
is reduced. 
\end{prop}
%%%%%%%%%%%%%%%%%%%%%%%%%%%%%%%%%%%%%%%%%%%%%%%%%%%%%%%%
\begin{pf}
Put $p_k=\lex^{-1}(k)$ for $k\geqq 1$.
By Proposition \ref{prop:hk=hkr} and Theorem \ref{thm:hkinR+}, we have
\begin{align*}
%&\{\hk(x)\mid x\in\lex^{-1}([1,n])\}=\{\hkr(x)\mid x\in\lex^{-1}([1,n])\}\\
%&=\{s(\lex^{-1}(1))s(\lex^{-1}(2))\cdots s(\lex^{-1}(k-1))\al(\lex^{-1}(k))\mid k\in[1,n]\}
\hk(p_k)=s(p_1)s(p_2)\cdots s(p_{k-1})\al(p_k)
=w_{\cdiag,\lex}[k-1]\al(p_k)\in R_+
\end{align*}
for all $k\in[1,n]$.
%and $\hk(p_k)\in R_+$ 
This implies that 
$$\ell(w_{\cdiag,\lex}[k-1]s(p_k))>\ell(w_{\cdiag,\lex}[k-1])\ \ (k\in[1,n]).$$
Therefore we have $\ell(w_{\cdiag,\lex}[n])=n$ and thus 
the expression \eqref{eq:w[n]} is reduced.
% \ref{lem:hump1}. 
\qed
\end{pf}
%

%%%%%%%%%%%%%%%%%%%%%%%%%%%%%%%%%%%%%%%%%%%%%%%%%%%%%%%%%
 For $w \in W$,
the set $$R(w)= R_+ \cap w R_-$$ is called 
the {\it inversion set} of $w$.
It is known for any reduced expression $w=s_{i_1}s_{i_2}\cdots s_{i_\ell}$
that $\ell(w)=|R(w)|$ and
$$
R(w)=\{\al_{i_1}, s_{i_1}\al_{i_2},s_{i_1}s_{i_2}\al_{i_3},
\dots,s_{i_1}s_{i_2}\cdots s_{i_{\ell-1}} \al_{i_\ell}\}.$$

By \eqref{eq;hk}, \eqref{eq:w[n]} and Proposition \ref{prop:reduced},
we obtain the following proposition:

%%%%%%%%%%%%%%%%%%%%%%%%%%%%%%%%%%%%%%%%%%%%%%%%%%%%%%%%%%%
\begin{prop}\label{prop:Randhk}
Let $\lex\in\LE(\cdiag)$ and $n\in\Z_{\geqq1}$.
Then it holds that
$$R(w_{\cdiag,\lex}[n])=\{\hk(x)\mid x\in\lex^{-1}([1,n])\}.$$
In particular, it holds that $R(w_{\cdiag,\lex}[n])\subset \DS(\predom_\cdiag).$
\end{prop}
%%%%%%%%%%%%%%%%%%%%%%%%%%%%%%%%%%%%%%%%%%%

%%%%%%%%%%%%%%%%%%%%%%%%%%%%%%%%%%%%%%%%%%%%%
Define
$$R(w_{\cdiag,\lex})=\bigcup_{n\geqq 1}R(w_{\cdiag,\lex}[n]).$$
Then 
$$R(w_{\cdiag,\lex})=\left\{\hk(x)\mid x\in \cdiag\right\}=\DS(\predom_\cdiag).$$
In particular, $R(w_{\cdiag,\lex})$ is independent of $\lex$ and
we will denote it just by $R(w_\cdiag)$ in the rest.

\begin{rem}
The set $R(w_\cdiag)$ can be thought as the ``inversion set" associated with
the semi-infinite word
$$w_{\cdiag,\lex}:=s(\lex^{-1}(1))s(\lex^{-1}(2))\cdots\cdots.$$
\end{rem}

 %%%%%%%%%%%%%%%%%%%%%%%%%%%%%%%%%%%%%%%%%%%%%%%%%%%

\begin{df}
Let $\zeta\in P$ be an integral weight.

\smallskip\noindent
(1) An element $w$ of $W$ is 
said to be $\zeta$-{\it pluscule} if
%there exists a reduced expression
%$w=s_{i_1}s_{i_2}\cdots s_{i_n}$ such that
\begin{align*}
&\bra \zeta,\al^\vee\ket=-1\ \ \text{for all } \al\in R(w).
\end{align*}
%&\bra \zeta, s_{i_1}s_{i_2}\cdots s_{i_{k-1}}\al_{i_k}^\vee \ket= -1,
%\text{ for all }
%k=1,2,\dots, n,\\
%
%\smallskip
\noindent
(2) An element $w$ of $W$ is 
said to be $\zeta$-{\it minuscule} if
\begin{align*}
&\bra \zeta,\al^\vee\ket=1\ \ \text{for all } \al\in R(w^{-1}).
%&\bra \zeta, s_{i_1}s_{i_2}\cdots s_{i_{k-1}}\al_{i_k}^\vee \ket= -1,
%\text{ for all }
%k=1,2,\dots, n,\\
\end{align*}
%
\if0
or equivalently 
\begin{equation}\label{eq3}
s_{i_k} \cdots s_{i_1}\zeta = \zeta + \sum_{j=1}^k \al_{i_j}
\end{equation}
for all $k=1,2,\dots,n$.
%, where $p_j = \lex^{-1}(j)\ (j\in \Z_{\geqq1})$,
\fi
\end{df}
%%%%%%%%%%%%%%%%%%%%%%%%%%%%%%%%%%%%%%%%%%%%
\begin{df}
An element $w \in W$ is said to be \alert{fully commutative} if any reduced 
expression of $w$ can be obtained from any other by using only the relations \eqref{eq:rel2}.
\end{df}
%The weight $\predom_\cdiag$ is predominant, namely,
%$\bra \predom_\cdiag,\al^\vee \ket\geqq -1$ for all $\al^\vee\in R^\vee_+$.
% Moreover,
%$w_{\cdiag,\lex}$ is $\predom_\cdiag$-``pluscule'' in the sense of the following proposition:
%Proposition \cite{prop;pluscule}.
%, which we will see later.

\begin{rem}
(1) An element $w\in W$ is $\zeta$-pluscule if and only if 
$w$ is $(w^{-1}\zeta)$-minuscule.

\noindent
(2) It is known that if $w$ is $\zeta$-minuscule for some integral weight $\zeta$
then $w$ is fully commutative (\cite{Stem1996}).
\end{rem}
%%%%%%%%%%%%%%%%%%%%%%%%%%%%%%%%%%%%%%%%%%%%%%%
By Proposition \ref{prop:Randhk}, we have the following:
%Let $\cdiag$ be a cylindric diagram.
%%%%%%%%%%%%%%%%%%%%%%%%%%%%%%%%%%%%%%%%%%%%%%%%%%%%
\begin{prop}\label{prop;pluscule}
Let %$\cdiag$ be a cylindric diagram and
 $\lex\in \LE(\cdiag)$ and $n\in\Z_{\geqq1}$.
 % a standard tableau  of $\cdiag$.
%Putting $p_j = \lex^{-1}(j)\ (j\in \Z_{\geqq1})$,
%For any $x\in\cdiag$, 
Then $w_{\cdiag,\lex}[n]$ is $\predom_\cdiag$-pluscule
and fully commutative.
%\begin{equation}\label{eq2}
%\bra \predom_\cdiag, 
%s(p_1) \cdots s(p_{k-1}) \al^\vee(p_k) 
%\hkr(x) \ket = -1
%\end{equation}
%for all $k \geqq1$,
\end{prop}
%%%%%%%%%%%%%%%%%%%%%%%%%%%%%%%%%%%%%%%%%

%%%%%%%%%%%%%%%%%%%%%%%%%%%%%%%%%%%%%%%%%%%%%%%%%%%%%%%
%

%%%%%%%%%%%%%%%%%%%%%%%%%%%%%%%%%%%%%%%%%%%%%%%%%%%%%%%%%%%%%%%
\subsection{Proof of Proposition~\ref{prop:hk=hkr}}
%%%%%%%%%%%%%%%%%%%%%%%%%%%%%%%%%%%%%%%%%%%%%%%%%%%%%%%%%%%%%%
%%%%%%%%%%%%%%%%%%%%%%%%%%%%%%%%%%%%%%%%%%%%%%%%%%%%%%
For $\lex\in \LE(\cdiag)$ and  $x\in\cdiag$, we put
\begin{align}\label{eq;hkr}
\hkr(x)
&=s(\lex^{-1}(1))s(\lex^{-1}(2))\cdots s(\lex^{-1}(n-1)) \al(\lex^{-1}(n)),
%\\
%&=w_{\cdiag,\lex}[n-1]\al(x)\notag,
\end{align}
where $n=\lex(x)$.
%We have $R(w_{\cdiag,\lex})=\{\hkr(x)\mid x\in \cdiag\}$.
%%%%%%%%%%%%%%%%%%%%%%%%%%%%%%%%%%%%%%%%%%%%%%%%%%%%%
%We put
%$p_j=\lex^{-1}(j)$ and
%$$\hkr(p_k)=s(p_1)s(p_2)\cdots s(p_{k-1})\al(p_k).$$
%Then $\cdiag=\{p_j\mid j\in\Z_{\geqq1}\}$
%and $\eta(p_k)=w_{\cdiag,\lex}[k-1]\al(p_k)$.

For $x\in \cdiag$,
%We proceed by induction on $x\in\cdiag$.
%If $x$ is minimal in $\cdiag$ then $\hk(x)=\al(x)=\hkr(x)$.
%Suppose that $x$ is not minimal and
 put $\xs=x+(1,0),\ \xe=x+(0,1),\ \xse=x+(1,1)\in \cyli_\peri$.
%Then at least one of $z_1$ and  $z_2$ belongs to $\cdiag$.
We will use the following lemma later:
%%%%%%%%%%%%%%%%%%%%%%%%%%%%%%%%
\begin{lem}\label{lem;rechkr}
Let $x\in \cdiag$.

\smallskip\noindent
$\mathrm{(1)}$ If $x\notin \btm$, then $\xs,\xe,\xse\in\cdiag$ and
\begin{align}\label{eq;hkrrec1}
\hkr(x)&=%\begin{cases}
\hkr(\xs)+\hkr(\xe)-\hkr(\xse).
%& \text{ if }x\notin\btm
%\end{cases}
\end{align}

\smallskip\noindent
$\mathrm{(2)}$ If $x\in \btm$, then $\xse\notin \cdiag$ and
\begin{empheq}[left = {\hkr(x)= \empheqlbrace\,}]{alignat=3}
%\begin{align}\label{eq;hkrrec2}
%\hkr(x)=\begin{cases}
&\al(x)+\hkr(\xs)+\hkr(\xe) &\text{ if } \xs,\xe\in\cdiag\\
&\al(x)+\hkr(\xs) &\text{ if } \xs\in\cdiag, \xe\notin\cdiag\\
&\al(x)+\hkr(\xe) &\text{ if } \xe\in\cdiag, \xs\notin\cdiag\\
&\al(x) &\text{ if } \xs,\xe\notin\cdiag.
%\end{cases}
%\end{align}
\end{empheq}
%where we set $\hkr(y)=0$ if $y\notin\cdiag$.
\end{lem}
%%%%%%%%%%%%%%%%%%%%%%%%%%
\begin{pf}
%Let us prove \eqref{eq;hkrrec}.
We put $p_k=\lex^{-1}(k)$ $(k\in\Z_{\geqq 1})$.

\smallskip\noindent
(1) 
 Let $x=p_j$, $\xse=p_i$, $\xe=p_{k_1}$ and $\xs=p_{k_2}$.
Then $j>k_1,k_2>i$ and  we may assume that $k_2<k_1$.
Put $\con(x)=r$. 
Then $\con(\xe)=r-1$, $\con(\xs)=r+1$.
We have
$$\hkr(x)=w_1s(p_i)w_2s(p_{k_1})w_3s(p_{k_2})w_4\al(p_j)=w_1s_r w_2 s_{r+1} w_3 s_{r-1}w_4 \al_r,$$
where $w_1=s(p_1)\cdots s(p_{i-1})$, $w_2=s(p_{i+1})\cdots s(p_{k_1-1})$,
$w_3=s(p_{k_1+1})\cdots s(p_{k_2-1})$  and $w_4=s(p_{k_2+1})\cdots s(p_{j-1})$.

Note that $\con(p_d)-r\neq 0,\pm1$
for all $d\in [i+1, j-1]\setminus\{k_1,k_2\}$.
Actually, if $\con(p_d)-r= 0,\pm1$ then 
$p_d$ is comparable with $p_j$ and $p_i$, and hence
$p_j>p_d>p_i$. But such $d$ must be $k_1$ or $k_2$.
%Using Lemma \ref{prop:content}, we have
%reason of $p_h$ ($k_2<h<j$) are incomparable with $p_j$
Now we have
\begin{align*}
\hkr(x)&=w_1s_r w_2 s_{r+1} w_3 s_{r-1} w_4 \al_r =w_1s_r w_2 s_{r+1} w_3 s_{r-1} \al_r \\
&=w_1s_r w_2 s_{r+1} w_3 (\al_{r-1}+\al_r)=\hkr(\xs)+w_1s_r w_2 s_{r+1} w_3 \al_r\\
&=\hkr(\xe)+w_1s_r w_2 s_{r+1} \al_r=\hkr(\xe)+w_1s_r w_2 (\al_r+\al_{r+1})\\
&=\hkr(\xs)+\hkr(\xe)+w_1s_r w_2 \al_r =\hkr(\xs)+\hkr(\xe)+w_1s_r\al_r\\
&=\hkr(\xs)+\hkr(\xe)-\hkr(\xse).
\end{align*}

\smallskip\noindent
(2) Suppose that $\xs,\xe\notin \cdiag$, or equivalently,
suppose that $x$ is  minimal element in $\btm$.
Let $x=p_j$.
Then $p_d\ (d\in [1,j-1])$
%,p_2,\dots,p_{k-1}$ 
is not comparable with $p_j$.
%(Accturally, if $p_j$ is comparable  with $p_k$,
%then $p_k<p_j$ by minimality of $p_k$. This implies $k=\lex(p_k)<\lex(p_j)=j$.)
Hence 
%by Proposition \ref{prop:content} (1), we have 
%$s(p_1) \cdots s(p_{k-1}) \alpha(p_k) = \alpha(p_k)$,
%and obtain $\pring_k=\bra  \predom_{\cdiag},\alpha^\vee(p_k) \ket = -1$ by (\ref{eq:opo}).
%and hence \eqref{eq:opo} holds.
$$\hkr(x)=s(p_1) \cdots s(p_{j-1})\alpha(p_j)
=\alpha(p_j)$$
The other cases are reduced to the case where $x$ is  minimal in $\btm$,
via a similar argument as in the proof of the statement (1),

\qed
\end{pf}
%%%%%%%%%%%%%%%%%%%%%%%%%%%%%%%%%%%%%%%%%%%%%%%%%%%%%%%%%%%%%%%%%%%%%%%%

%%%%
\begin{pf}[Proposition~\ref{prop:hk=hkr}]
Let $x\in \cdiag$.
%We proceed by induction on $x\in\cdiag$.
%If $x$ is minimal in $\cdiag$ then $\hk(x)=\al(x)=\hkr(x)$.
%Suppose that $x$ is not minimal and
Put $\xs=x+(1,0),\ \xe=x+(0,1),\ \xse=x+(1,1)$.
%Then at least one of $z_1$ and  $z_2$ belongs to $\cdiag$.
It is easy to see the following:
\begin{align}\label{eq;hkrec} 
\hk(x)=\begin{cases}
\hk(\xs)+\hk(\xe)-\hk(\xse)& \text{ if }x\notin\btm\\
\al(x)+\hk(\xs)+\hk(\xe)& \text{ if }x\in\btm \text{ and } \xs,\xe\in \cdiag\\
\al(x)+\hk(\xs)& \text{ if }x\in\btm \text{ and } \xs\in \cdiag,\ \xe\notin\cdiag\\
\al(x)+\hk(\xe)& \text{ if }x\in\btm \text{ and } \xe\in \cdiag,\ \xs\notin\cdiag\\
\al(x)& \text{ if }x\in\btm \text{ and } \xs,\xe\notin \cdiag
\end{cases}
\end{align}
%In particular,  we have
%$\hk(x)=\al(x)$ for a minimal element $x$ in $\cdiag$, and 
%the recurrence relation \eqref{eq;hkrec} determines $\hk(x)$ for any $x\in\cdiag$
%uniquely.
%if $x$ is minimal in $\cdiag$
%where $z_1=x+(1,0),\ z_2=x+(0,1),\ y=x+(1,1)$.
On the other hand, we have shown that
%Hence the statement will be shown if we prove that 
$\hkr(x)$ satisfies the same 
recurrence relations in Lemma \ref{lem;rechkr}.\qed
\end{pf}

\newcommand\pring{f}

%%%%%%%%%%%%%%%%%%%%%%%%%%%%%%%%%%%%%%%%%
\section{Poset structure of cylindric diagrams}\label{sec:1stresult}
%%%%%%%%%%%%%%%%%%%%%%%%%%%%%%%%%%%%%%%%

%%%%%%%%%%%%%%%%%%%%%%%%%%%%%%%%%%%%%%%%%%%%%%%%%%%%%%
\subsection{Partial orders on the inversion set}\label{sec:inversion set}
%%%%%%%%%%%%%%%%%%%%%%%%%%%%%%%%%%%%%%%%%%%%%%%%%%%%%%

Recall that $Q$ denote the root lattice: $Q=\bigoplus_{i\in\Z/\kappa\Z}\Z \al_i$.
%%%%%%%%%%%%%%%%%%%%%%%%%%%%%%%%%%
\begin{df}
Define the partial order $\oord$ on $Q$ by
\begin{equation*}
\alpha \oord \beta \iff
\beta-\alpha \in Q_+=\bigoplus_{i\in\Z/\kappa\Z}\Z_{\geqq0}\al_i
%=\left\{\sum_{i\in\Z/\kappa\Z}k_i \al_i\mid k_i\in\Z_{\geqq0}
%\ (i\in\Z/\kappa\Z)\right\}.
%\text{ is a sum of elements in }\Pi
\end{equation*}
The order $\oord$ is called the {\it ordinary order}.
\end{df}
%%%%%%%%%%%%%%%%%%%%%%%%%%
The restriction of the ordinary order defines a poset structure on
 $R(w_\cdiag)$.

Let $\cdiag$ be a cylindric diagram in $\cyli_\peri$ with $|\peri|=\kappa$. 
We have introduced a poset structure on $\cdiag$
and also  have seen that the map $\hk$ gives a bijection 
between $\cdiag$ and $R(w_\cdiag)$.
Remark that this is not a poset isomorphism as seen in the following example:
%Note that $(R(w_\cdiag),\oord)$  is not equivalent to $(\cdiag,\leq)$
%%%%%%%%%%%%%%%%%%%%%%%%%%%%%
\begin{ex}\label{ex:ord}
%For example,
Let $\lm=(4,2),\ \peri=(2,-2)$ and consider the cylindric diagram $\cyl\lm$ in
$\cyli_\peri$.
Then $x=\pi(1,2)$ and $y=\pi(2,1)$ are incomparable in $\cyl\lm$.
On the other hand, $\hk(x)=\delta+\al_3$ and $\hk(y)=\al_0+\al_2+\al_3$,
and hence $\hk(x)-\hk(y)=\al_1+\al_3$. This implies $\hk(y)\oord\hk(x)$. 
\end{ex}
%%%%%%%%%%%%%%%%%
%\end{rem}

We will introduce a modified ordinary order $\neword$, 
%which is
%will turn out to be 
for which we will have
$(\cdiag,\leq)\cong (R(w_\cdiag),\neword)$.

Let $\btm=\{b_i\mid i\in\Z/\kappa\Z\}$
be the bottom set of $\cdiag$,
where $b_i$ is the element such that $\con(b_i)=i$.
Let $\btmmax$ (resp. $\btmmin$) denote the set of maximal (resp. minimal)
elements in $\btm$.

%%%%%%%%%%%%%%%%%%%%%%%%
\if0
For $i \in \Z/\kappa\Z$,
we denote by $\tilde{i}$ its representative contained in $\{0,1,\dots,\kappa-1\}$.
%
Let $i,j \in \Z/\kappa\Z$ with $i \neq j$, and define
$$
((i,j)) =
\begin{cases}
\{ k \in \Z/\kappa\Z \mid \tilde{i} < \tilde{k} < \tilde{j}\} & \text{if $\tilde{i}<\tilde{j}$} \\
\{ k \in \Z/\kappa\Z \mid \tilde{k} < \tilde{j} \ \text{or}\  \tilde{i} < \tilde{k}\} & \text{if $\tilde{i}>\tilde{j}$},
\end{cases}
$$
and put
\[ \al[i,j] = \al_i + \sum_{k \in ((i,j))} \al_k + \al_j. \]
%Note that $\al[i,j]$ may be the null root.
%which is an element of $R \sqcup \Z\delta$.
\fi

%%%%%%%%%%%%%%%%%%%%%%%%%%%%%%%%%%%%%%%%%%%%%%%%%%%%%%%%%%%%%%%%%%%%%%%
\begin{df}\label{df:Deltaord}
Define
\[ \Pi_\cdiag = \Pi_\cdiag^0 \sqcup \Pi_\cdiag^\mathrm{arm} \sqcup \Pi_\cdiag^\mathrm{leg}. \]
Here,
\begin{align*}\label{eq:delta}
\Pi_{\cdiag}^0 
&=\{\al(x) \mid x\in\btm\setminus (\btmmax\sqcup \btmmin)\},\\
%\{ \alpha_i \in \Pi \mid \bra \predom_{\cdiag},\alpha^\vee_i \ket = 0 \}, \\
%
\Pi_{\cdiag}^\mathrm{arm} 
&=\left\{\al(x)+\sum_{y\in \arm(x)}\al(y)\ \middle|\  x\in\btmmax\right\},\\
%\{ \al[i,j] \middle\mid \bra\predom_\cdiag,\al^\vee_i\ket=1,\ \bra\predom_\cdiag,\al^\vee_j\ket=-1,\ 
%\bra\predom_\cdiag,\al^\vee_k\ket=0\ \forall k \in ((i,j))\}, \\
%
\Pi_{\cdiag}^\mathrm{leg} 
&=\left\{\al(x)+\sum_{y\in \leg(x)}\al(y)\ \middle|\ x\in\btmmax\right\}.
%\{ \al[i,j] \mid \bra\predom_\cdiag,\al^\vee_i\ket=-1,\ \bra\predom_\cdiag,\al^\vee_j\ket=1,\ 
%\bra\predom_\cdiag,\al^\vee_k\ket=0\ \forall k \in ((i,j))\}.
\end{align*}
\end{df}
%%%%%%%%%%%%%%%%%%%%%%%%%%%%%%%%%%%%%%%%%%%
%It is easy to see that $\bra\zeta_\cdiag^\circ,\al^\vee\ket=0$ for
%$\al\in \Pi_\cdiag$.
Note that %$\al[i,j]$ may be the null root and 
$\Pi_\cdiag \subset R_+ \sqcup \Z_{\geqq0}\delta$.

\begin{ex}
For the cylindric diagram described in Fig. \ref{content}, we have
\begin{align*}
&\Pi_\cdiag^0=\{\al_3,\al_5,\al_7\},\\
&\Pi_\cdiag^\mathrm{arm} =\{\al_6+\al_7+\al_8, \al_2+\al_3+\al_4, \al_0+\al_1\},\\
&\Pi_\cdiag^\mathrm{leg} =\{\al_4+\al_5+\al_6, \al_1+\al_2, \al_0+\al_8\}.
\end{align*}
\end{ex}
\begin{ex}
Let $\lm=(n)$ and $\peri=(1,-n+1)$.
Then, for the corresponding cylindric diagram $\cyl\lm$, we have
\begin{align*}
&\Pi_{\cyl\lm}^0=\{\al_1,\al_2,\dots,\al_{n-2}\},\\
&\Pi_{\cyl\lm}^\mathrm{arm} =\{\delta\},\\
&\Pi_{\cyl\lm}^\mathrm{leg} =\{\al_0+\al_{n-1}\}.
\end{align*}
\end{ex}
%%%%%%%%%%%%%%%%%%%%%%%%%%%%%%%%%%%%%%%%
\begin{df}\label{df:neword}
Define the partial order $\neword$ on $R(w_{\cdiag})$ by
%as the transitive closure of the relations
\begin{equation}\label{eq:neword}
\alpha \neword \beta \iff
\beta-\alpha \in\sum_{\gam\in\Pi_\cdiag}\Z_{\geqq0}\gam
=\left\{
\sum_{\gamma\in \Pi_\cdiag}
k_\gamma \gamma \ \middle| \ k_\gamma\in\Z_{\geqq0}
\ (\forall \gamma \in\Pi_\cdiag)\right\}.
%\text{ is a sum of elements in }\Pi_\cdiag
\end{equation}
\end{df}
%%%%%%%%%%%%%%%%%%%%%%%%%%%%%%%%%%%%%%

\if0
It is easy to see the following:
%%%%%%%%%%%%%%%%%%%%%%%%%%%%%%%%%%%%%%%%%%
\begin{lem}\label{lem:ordersinR}
Let $\al,\beta\in R(w_\cdiag)$.
Then
$$
%\al\tcord\beta\ 
%\implies\ 
\al\neword \beta \ 
\implies\
\al \oord\beta.$$
\end{lem}
%%%%%%%%%%%%%%%%%%%%%%%%%%%%%%%%
\fi

%%%%%%%%%%%%%%%%%%%%%%%%%%%%%%%%%%%%%%%%%%%%%%%%%%%%%%%%%%%%%%%%%%%%%%%%%%

%%%%%%%%%%%%%%%%%%%%%%%%%%%%%%%%%%%%%%%%%%%%%%%%%
\begin{prop}\label{prop:cdiagtotcord}
Let $x,y\in\cdiag$. Then 
$$ x\le y \implies \hk(x)\neword \hk(y).$$
In other words, the bijection 
%$\hk:\cdiag \to R(w_\cdiag)$ is order preserving
%isomorphism
$$ \hk:(\cdiag,\leq) \to (R(w_\cdiag),\neword)$$
is order preserving.
%$\hpord$ and $\tcord$on $R(w_\cdiag)$ coincide.
\end{prop}
%%%%%%%%%%%%%%%%%%%%%%%%%%%%%%%%%%%%%%%%%%%%%%%%%%%%%%%%%%%%%%%%%%%

\begin{pf}
%[Proposition \ref{prop:sankaku}]
 We assume that $y$ covers $x$, and
% We suppose $x \gtrdot y$ and
 will show that $\hk(y)-\hk(x)\in\Pi_\cdiag$
 by induction on $y$ concerning the poset structure on $\cdiag$.

Put $y^S=y+(1,0),\ y^E=y+(0,1),\ y^{SE}=y+(1,1)$. Then
 $x=y^S$ or $x=y^E$.
 
 When $y\in\btm$, it follows from \eqref{eq;hkrec}
that  $\hk(y)- \hk(x)\in \Pi_\cdiag$.

Suppose that $y\notin \btm$.
Note that $y^{SE}\in \cdiag$.
Since $y^E$ covers $y^{SE}$ and $y>y^E$,
we have
$\hk(y^E)-\hk(y^{SE})\in \Pi_\cdiag$ by induction hypothesis.
By the recursion relation \eqref{eq;hkrec}, we have 
$$\hk(y)-\hk(y^S)=\hk(y^E)-\hk(y^{SE})\in\Pi_\cdiag.$$
Similar argument implies 
$\hk(y)-\hk(y^E)\in \Pi_\cdiag$.
In both cases, we have $\hk(y)-\hk(x)\in\Pi_\cdiag$.
Therefore, the statement
%\eqref{eq:hkleq} 
is proved.
%$x>y$ implies $\hk(x)\tcord \hk(y)$.
\qed
%
\if0
%%%%%%%%%%%%%%%%%%%%%%%%%%%%%%%%%%%%%
In order to show 
\begin{equation}
\label{eq:hkleq2}
\hk(x)\tcord\hk(y)\implies 
x\leq y,
\end{equation}
it is enough to show
\begin{equation}
\label{eq:hkleq3}
 \al\tcord\beta\implies\alpha \hpord \beta.
\end{equation}
for all $\al,\beta\in R(w_\cdiag)$
as we know $(R(w_\cdiag), \hpord )\cong (\cdiag,\leq)$.

Suppose that
$\al \tcord \beta$ is a covering relation.
Note that $\beta-\al$ is an element of
$\Pi_\cdiag\subset  R \sqcup \Z\delta$.
%
Assume that $\bra\al,\beta^\vee\ket= 0$.
Then $\bra\beta-\al,\beta^\vee\ket=2$ and
Lemma \ref{lem:moddelta} implies $\al\equiv 0\mod \delta$.
This contradicts the fact that $\al$ is a real root.
Therefore $\bra\al,\beta^\vee\ket\neq 0$.
%Suppose $\gam(b)-\gam(a) \in \Pi_\cdiag^0$.
%Then 
Thus, $\al$ and $\beta$ is comparable  with respect to $\hpord$. %and hence
% $\gam^{-1}(\al)$ and $\gam^{-1}(\beta)$ is comparable in $(\Z_{\geq1},\hpord)$.
%Since $\hk$ is order preserving, 
Now, $\al\tcord\beta$ implies $\al\hpord\beta$.
\fi
%%%%%%%%%%%%%%%%%%%%%%
%\qed
\end{pf}

It is easy to see that 
$$\al\neword\beta\implies\al\oord \beta$$
for any $\al, \beta\in R(w_\cdiag)$. Thus we have the following:
%%%%%%%%%%%%%%%%%%%%%%%%%%%%%%%%
\begin{cor}\label{cor:ordersinR}
Let $x,y\in\cdiag$. Then
$$x<y\implies\hk(x)\oord \hk(y).$$
\end{cor}
%%%%%%%%%%%%%%%%%%%%%%%%%%%%%%%%%

\newcommand\tcoord{\unlhd^\mathrm{tc}}
%%%%%%%%%%%%%%%%%%%%%%%%%%%%%%%%%%%%%%%%%%%%%%%%%%%%%%%%%%%%%%%%%%%%%%%%%

\if0
%%%%%%%%%%%%%%%%%%%
\begin{rem}
(1) 
%In the definition of the order $\tcord$,
%there is some ambiguity and 
%one can adopt for example  $\Pi$ in stead of $\Pi_\cdiag$:
Define the partial order $\tcoord$ on $R(w_{\cdiag})$ as 
the transitive closure of the relations
\begin{equation}\label{eq:tcoord}
\alpha \tcoord \beta \ \text{whenever}\ \beta-\alpha \in \Pi.
\end{equation}
The same argument in the proof of Proposition \ref{prop:cdiagtotcord}
%and Proposition \ref{prop:orderpres},
 implies  the poset isomorphism
 $(\cdiag,\leq)\cong (R(w_\cdiag), \tcoord)$ and thus
 % isomorphism
\begin{equation}\label{eq:tc=tc}
(R(w_\cdiag), \tcoord) = (R(w_\cdiag), \tcord).
\end{equation}
%$$\alpha \tcoord \beta\ \Longleftrightarrow 
%\ \alpha \tcord \beta.$$
%\smallskip
\noindent
(2) It is known that the partial order $\tcoord_Q$ on $Q$ defined as 
the transitive closure of the same relations \eqref{eq:tcoord} coincides with
the ordinary order $\oord$.
On the other hand, as we have seen in Example \ref{ex:ord},
$$(R(w_\cdiag), \oord)\ncong (R(w_\cdiag), \tcord)$$
in general.
Thus the  restriction of $\tcoord_Q$ to $R(w_\cdiag)$ is not equivalent to $\tcoord$.
\end{rem}
%%%%%%%%%%%%%%%

\fi
%%%%%%%%%%%%%%%%%%%%%%%%%%%%%%%%%%%%%%%%%%%%%%%%%%%%%%%%%%%%%%%%%%%%%%%%%%%%%%%%

%%%%%%%%%%%%%%%%%%%%%%%%%%%%%%%%%%%%%%%%%%%%%%%%%%%%%%%%
\subsection{Poset isomorphism}
%%%%%%%%%%%%%%%%%%%%%%%%%%%%%%%%%%%%%%%%%%%%%%%%%%%%%%%
Our next goal is to prove that the order preserving map
$$\hk: (\cdiag,\leq)\to (R(w_\cdiag),\neword)$$ is actually a poset isomorphism.
%%%%%%%%%%%%%%%%%%%%%%%%%%%%%%%%%%%%%%%%%%%%%%%%
We start with some preparations.

%%%%%%%%%%%%%%%%%%%%%%%%%%%%%%%%%%%%%%%%%%%%

As before, we denote by 
$\Supp(\al)$ the support of
%=\{b_i \mid \bra  (i\in\Z/\kappa\Z)\}\subset \btm,$
 $\al\in Q$.
The following lemma is almost obvious from Definition \ref{df:Deltaord}.
%%%%%%%%%%%%%%%%%%%%%%%%%%%%%%%%%%%%%%%%%%%%%%%%%%
\begin{lem}\label{lem:supportofPi}
Let $\al\in\Pi_\cdiag$.
Then 
%it is easy to see that
\begin{equation}\label{eq:max=min}
|\Supp (\al)\cap \btmmax|=|\Supp (\al) \cap \btmmin|=
\begin{cases}
0&\ (\al\in\Pi_\cdiag^0)\\
1&\ (\al\in\Pi_\cdiag^\mathrm{arm}\sqcup\Pi_\cdiag^\mathrm{leg}).
\end{cases}
\end{equation}
\end{lem}
%%%%%%%%%%%%%%%%%%%%%%%%%%%%%%%%%%%%%%%%%%%%

%It is easy to see that

%By using Lemma \ref{lem:addnull},
%we can see that the following lemma:
It is easy to see the next lemma:
%%%%%%%%%%%%%%%%%%%%%%%%%%%%%%%%%%%%%%%%%%%%%%%%%%%%%%%%%%%%
\begin{lem}\label{lem200}
$\mathrm{(1)}$ Let $\al\in R_+$. Then
$\nmult(\al)=\max\{\nmult\in\Z \mid\al-\nmult\delta \in R_+\}.$

\smallskip\noindent
$\mathrm{(2)}$ Let $x,y \in \cdiag$.
If $x<y$ then $\nmult(\hk(x)) \leqq \nmult(\hk(y))$.
\end{lem}
%%%%%%%%%%%%%%%%%%%%%%%%%%%%%%%%%%%%%%%%%%%%%%%%%%
\begin{pf}
(1) Follows from Lemma \ref{lem:addnull}.
%\eqref{eq:R+2}.
%\end{pf}
%\if0
%\begin{pf}

%\smallskip
\noindent
(2) Suppose $x<y$.
By Corollary \ref{cor:ordersinR}, we have
$\nmult(\hk(x))\delta\oord \hk(x) \oord \hk(y)$.
%Since $\hk(x)-\nmult(\hk(x))\delta\in R_+$, it holds that 
%$\nmult(\hk(x))\delta\oord \hk(x)$.
%Suppose $x<y$. Then by Corollary \ref{cor:ordersinR}, 
%we have $\hk(x)\oord\hk(y)$ and hence $\nmult(\hk(x))\delta\oord \hk(y)$.

As $\hk(y)-\nmult(\hk(x))\delta$ is in $R$ by Lemma \ref{lem:addnull},
it must be a positive root.
This means $\nmult(\hk(x))\leqq \nmult(\hk(y))$.\qed
\end{pf}
%\fi

%%%%%%%%%%%%%%%%%%%%%%%%%%%%%%%%%%%%%%%%%%%%%%%%%%%%%%%%%%%%%%%%%
\begin{lem}\label{lem201}
 Let $x,y \in \cdiag$ such that $\nmult(\hk(x))=\nmult(\hk(y))=0$.
Then 
$$x<y \iff \hk(x) \oord \hk(y).$$
In particular, if $x$ and $y$ are incomparable,
then $\hk(x)$ and $\hk(y)$ are also incomparable with respect to $\oord$.
\end{lem}
%%%%%%%%%%%%%%%%%%%%%%%%%%%%%%%%%%%%%%%%%%%%%%%%%%%%%%%%%%%%%%%%%%%%
\begin{pf}

\smallskip\noindent
By Corollary \ref{cor:ordersinR}, we have
% and Proposition \ref{prop:cdiagtotcord}, the implication
$$x<y\implies\hk(x) \oord \hk(y).$$
We shall prove the opposite implication.
Suppose $\hk(x) \oord \hk(y)$. 
Then noting that $0\oordneq \hk(x), \hk(y)\oordneq \delta$, %by \eqref{eq:R+}.
we have $\Supp(\hk(x))\subset\Supp(\hk(y))\subsetneq \btm$.
%Thus $\Supp(\hk(x))$ detrmines $\hk(x)$ for $x\in\cdiag$ with $\nmult(x)=0$
%and 
%$$x<y\Longleftrightarrow \Supp(\hk(x))\subset\Supp(\hk(y))$$ 
%Note that $Supp(\hk(x))$
%It follows that
%from \eqref{eq:rimhook} that
%Now it can be seen that
%$\Supp(\hk(x))\subset \Supp(\hk(y))$
This implies $x<y$.
\qed
\end{pf}

%%%%%%%%%%%%%%%%%%%%%%%%%%%%%%%%%%%%%%%%%%%%%%%%%%%%%%%%%%%%%%%
\begin{lem}\label{lem202}
Let $x,y \in \cdiag$.
Suppose that  $x$ and $y$ are incomparable in $\cdiag$.
Then $\nmult(\hk(y))-\nmult(\hk((x))=1,0$ or $-1$. Moreover the followings hold:

\noindent
$\mathrm{(1)}$ If $\nmult(\hk(y))-\nmult(\hk(x))=1$, then 
$$\hk(y)-\delta \oord \hk(x)\oord \hk(y).$$
\noindent
$\mathrm{(2)}$ If $\nmult(\hk(y))-\nmult(\hk(x))=-1$, then 
$$\hk(x)-\delta \oord \hk(y)\oord \hk(x).$$
\noindent
$\mathrm{(3)}$ If $\nmult(\hk(y))-\nmult(\hk(x))=0$, then $\hk(x)$ and $\hk(y)$ are incomparable 
with respect to $\oord$.
\end{lem}
%%%%%%%%%%%%%%%%%%%%%%%%%%%%%%%%%%%%%%%%%%%%%%%%%%%%%%%%%%%%%%%%%%%%%%
%
\begin{pf}
In this proof, we denote $N(\hk(x))$ by $N(x)$ for $x\in\cdiag$.
Put 
$$x_k=x+(\nmult(x)-k)(0,\ell),\ \ y_k=y+(\nmult(y)-k)(0,\ell)$$ 
for $k\in\Z_{\geqq0}$.
Then $\nmult(x_k)=\nmult(y_k)=k$.
Putting $n=\nmult(x)$, one can see that 
$$\cdiag\setminus \left(
\{z\in\cdiag\mid z\geqq x\}\sqcup \{z\in\cdiag\mid z \leqq x\}\right)
=[x_{n-1}-(1,1),x_{n+1}+(1,1)].
$$
As
 $x$ and $y$ are incomparable, 
 $y$ belongs to this interval and
hence 
\begin{equation}\label{eq:ininterval}
x_{n-1}<y<x_{n+1}
\end{equation}
and
$n-1\leqq \nmult(y)\leqq n+1$ by Lemma \ref{lem200}.
Namely, we have $\nmult(y)-\nmult(x)=-1,0$ or $1$.

%Now we suppose that  $\hk(x)<\hk(y)$.

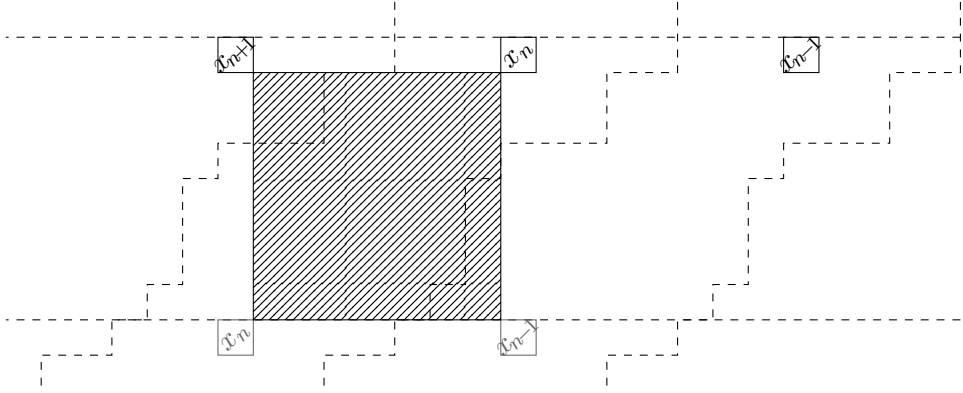
\begin{figure}[h]
\begin{center}
\begin{tikzpicture}[scale=.47, rotate=-90]

%\draw[help lines] (1,1) grid (9,9);

\draw[dashed, yshift=2cm] (0,7) -- (1,7) -| (2,5) -| (4,2) -| (5,1) -| (8,0) -| (9,-1) -| ($(2,5)+(8,-8)$) -- ($(3,5)+(8,-8)$);

\draw[dashed, yshift=-6cm] (0,7) -- (1,7) -| (2,5) -| (4,2) -| (5,1) -| (8,0) -| (9,-1) -| ($(2,5)+(8,-8)$) -- ($(3,5)+(8,-8)$);
\draw[dashed, yshift=-14cm] (0,7) -- (1,7) -| (2,5) -| (4,2) -| (5,1) -| (8,0) -| (9,-1) -| ($(2,5)+(8,-8)$) -- ($(3,5)+(8,-8)$);

\draw[dashed] (1,9) -- (1,-18);
\draw[dashed] (9,9) -- (9,-18);

\coordinate (x) at (1,-4) node[rotate=45] at ($(x)+(.5,.5)$) {\small$x_n$};
\draw (x) rectangle +(1,1);
\node[black!60, rotate=45] at ($(x)+(.5,.5)+(8,-8)$) {\small$x_n$};
\draw[black!60] ($(x)+(8,-8)$) rectangle +(1,1);

\node[rotate=45] at ($(x)+(.5,.5)+(0,-8)$) {\small$x_{\! n \! + \! 1}$};
\draw ($(x)+(0,-8)$) rectangle +(1,1);

\node[black!60, rotate=45] at ($(x)+(.5,.5)+(8,0)$) {\small$x_{\! n \! - \! 1}$};
\draw[black!60] ($(x)+(8,0)$) rectangle +(1,1);

\node[rotate=45] at ($(x)+(.5,.5)+(0,8)$) {\small$x_{\! n \! - \! 1}$};
\draw ($(x)+(0,8)$) rectangle +(1,1);

\filldraw[pattern=north east lines] ($(x)+(1,0)$) rectangle ($(x)+(8,-7)$);

\end{tikzpicture}
\end{center}
\caption{The cells in the shadow are incomparable with $x=x_n$.}
\end{figure}

\noindent
(1) Suppose that $\nmult(y)-\nmult(x)=1$.
In this case,  
$$\hk(y)-\hk(x)
= \hk(y_0)+ \delta - \hk(x_0).$$
By definition, $\hk(x_0)$ and $\hk(y_0)$ are positive roots.
By Lemma \ref{lem:addnull},  $\hk(x_0)-\delta$ 
is also a root and it is not positive.
Therefore $\delta-\hk(x_0)\in R_+$
and $\hk(y)-\hk(x)$ is a sum of two positive roots.
This implies $\hk(x)\oord \hk(y)$.
%Now the statement has been shown.
Combining with \eqref{eq:ininterval}, we have
$x-\delta\oord y\oord x$.

\noindent
(2) Follows from (1).
%Suppose that $\nmult(y)-\nmult(x)=-1$.
%Then by the same argument as above, we have
%$y-\delta\oord x\oord y$.

\noindent
(3) Suppose that $\nmult(y)-\nmult(x)=0$.
Note that $x_0$ and $y_0$ are incomparable this case, and
it follows from Lemma \ref{lem201} that
 $\hk(y_0)$ and $\hk(x_0)$ are also incomparable with respect to $\oord$.
Now we have
$$\hk(y)-\hk(x)
= \hk(y_0)+\nmult(y)\delta - \hk(x_0)-\nmult(x)\delta
= \hk(y_0)-\hk(x_0).$$
 and hence 
$\hk(x)$ and $\hk(y)$ are incomparable with respect to $\oord$.
\qed
\end{pf}

\if0
If $\nmult(y)-\nmult(x)=0$, then it must be $x_0<y_0$,
and hence $x<y$, it contradicts the incomparably $x$ and $y$.
%
%If $\nmult(y)-\nmult(x)=1$,
%then $\hk(y_0)$ and $\delta-\hk(x)$ are elements of $R_+$.
%
Suppose $\nmult(y)-\nmult(x) \geqq2$.
We remark that $y_{\nmult(y)+2}=y+(-m,-\ell)$.
Moreover, we remark that $x$ and $y$ are incomparable if and only if
$y \in [x+(-1,-1),x+(1-m,1-\ell)]$.
Then $x$ and $y$ are comparable, it contradicts.
%
Therefore, it must be $\nmult(y)=\nmult(x)+1$.
%\end{pf}
\fi

%%%%%%%%%%%%%%%%%%%%%%%%%%%%%%%%%%%%%%%%%%%%%%%%%%%%%%%%%%%
\begin{thm}\label{th:cylord}
The map 
$$\hk: (\cdiag,\leq)\to (R(w_\cdiag),\neword)$$ is a poset isomorphism.
%%%%%%%%%%%%%%%%%%%%%%%%%%%%%%%%%%%%%%%%%%%%%%%%
%Let $x,y\in\cdiag$. Then
%$$x\leq y\ \Longleftrightarrow\  \hk(x)\neword \hk(y).$$
%Let $\al,\beta\in R(w_\cdiag)$. Then
%$$\al\neword \beta\Longleftrightarrow \al\tcord\beta$$
\end{thm}
%%%%%%%%%%%%%%%%%%%%%%%%%%%%%%%%%%%%%%%%%%%%%%%%%%%%%%%%%%%
%%%%%%%%%%%%%%%%%%%%%%%%%%%%%%%%%%%%%%%%%%%%%%%
\begin{pf}
%[Proposition \ref{prop:tcord}]
%It is easy to see that
By  Proposition \ref{prop:cdiagtotcord},
%Lemma \ref{lem:ordersinR},
we have
$$x\leq y\implies\hk(x) \neword\hk(y),$$
$$y\leq x\implies\hk(y) \neword\hk(x).$$
Thus the statement follows if we prove that
$$x\text{ and }y\text{ are incomparable}\implies
\hk(x)\text{ and }\hk(y)\text{ are incomparable with respect to }\neword.$$
%We assume that $x$ and $y$ is imcomparable, and then  will prove that
%$\hk(x)$ and $\hk(y)$ is imcomparable 
%with respect to $\neword$.
% then the statement will be proved.
%Take
% $x,y \in \cdiag$ such that $\al=\hk(x)$ and $\beta=\hk(y)$.
%By the assumption and the equivalence
%$(\cdiag,\leq)\cong(R(w_\cdiag),\tcord)$,
%the cells $x$ and $y$ are incomparable in $\cdiag$.
%
Suppose that $x$ and $y$ are incomparable.
Then, putting $n=N(\hk(x))$, we have
$N(\hk(y))=n+1,n$ or $n-1$ by Lemma \ref{lem202}.
%, and it is enough to show that 
%$N(y)-N(x)$ is actually $0$.

First we assume that $N(y)=n$.
Then Lemma \ref{lem202} implies that $\hk(x)$ and $\hk(y)$ must be incomparable.

Next, assume that $N(y)=n+1$. Then
\[ \hk(y)-\hk(x) = \hk(y_0) + \delta - \hk(x_0), \]
where $x_0=x+n(0,\ell)$ and $y_0=y+(n+1)(0,\ell)$.
% with  $n=\nmult(x)$.
By Lemma \ref{lem:supportofhk}, we have
\begin{align*}
&|\Supp(\hk(y_0))\cap\btmmax|+1=|\Supp(\hk(y_0))\cap\btmmin|,\\
&|\Supp(\delta-\hk(x_0))\cap\btmmax|=|\Supp(\delta-\hk(x_0))\cap\btmmin|+1.
\end{align*}
They are not compatible with  Lemma \ref{lem:supportofPi} and thus 
%neither $\hk(y_0)$ or $\delta-\hk(x_0)$
%belong to %be a sum of elements in
%  $\Z_{\geqq 0}\Pi_\cdiag$.
we have  
\begin{equation}\label{eq:hknotinpi}
\hk(y_0)\notin \Z_{\geqq 0}\Pi_\cdiag,\ \ 
\delta-\hk(x_0)\notin \Z_{\geqq 0}\Pi_\cdiag
\end{equation}
%%%%%%%%%%%%%%%%%%%%%%%%%%%%%%%%%%%%%%%%%
\begin{figure}[h]
\begin{center}
\begin{tikzpicture}[scale=.2, rotate=-90]

\draw[thick, dashed] (-2,1) -- (20,1);

%\draw[help lines] (1,1) grid (20,20);
\draw[ultra thick] (-1,20) -| (3,18) -| (4,17) -| (6,15) -| (10,12) -| (13,9) -| (14,8) -| (16,7) -| (17,5) -| (18,1) -| (20,1);

\draw[blue!60, thick] (1,20) -- (-1,20);
\draw[blue!60, thick] (1,19) -- (-1,19);
\fill[pattern=north east lines, pattern color=blue!30] (1,20) rectangle (-1,19);
\draw[blue, thick] (1,20) -| (3,18) -| (4,17) -| (5,16) -| (3,17) -| (2,19) -- (1,19);
\fill[pattern=north east lines, pattern color=blue!80] (1,20) -| (3,18) -| (4,17) -| (5,16) -| (3,17) -| (2,19) -| (1,20) -- cycle;

\draw[thick, dashed] (1,-18) -- (1,20);

\coordinate (x0) at (5,3) node at ($(x0)+(-.8,.5)$) {$x_0$};
\filldraw (x0) rectangle +(1,1);	% x0
\draw ($(x0)+(1,0)$) -- +(12,0);	% leg of x0
\draw ($(x0)+(1,1)$) -- +(12,0);
\draw ($(x0)+(0,1)$) -- +(0,13);	% arm of x0
\draw ($(x0)+(1,1)$) -- +(0,13);
%\draw[blue, thick, xshift=19cm, yshift=-19cm] (1,20) -| (3,18) -| (4,17) -| (5,16) -| (3,17) -| (2,19) -| (-2,22) -| (-1,20) -- (1,20) -- cycle;
%\fill[xshift=19cm, yshift=-19cm, pattern=north east lines, pattern color=blue!80] (1,20) -| (3,18) -| (4,17) -| (5,16) -| (3,17) -| (2,19) -| (-2,22) -| (-1,20) -- (1,20) -- cycle;

\coordinate (y0) at (8,7) node at ($(y0)+(-.8,.5)$) {$y_0$};
\filldraw (y0) rectangle +(1,1);
\draw ($(y0)+(1,0)$) -- +(6,0);	% leg of y0
\draw ($(y0)+(1,1)$) -- +(6,0);
\draw ($(y0)+(0,1)$) -- +(0,7);	% arm of y0
\draw ($(y0)+(1,1)$) -- +(0,7);
\draw[red, thick] (8,15) -| (10,12) -| (13,9) -| (14,8) -| (16,7) -| (13,7) -| (13,8) -| (12,11) -| (9,14) -| (8,15) -- cycle;
\fill[pattern=north east lines, pattern color=red!80] (8,15) -| (10,12) -| (13,9) -| (14,8) -| (16,7) -| (13,7) -| (13,8) -| (12,11) -| (9,14) -| (8,15) -- cycle;
\node[anchor=west, red] at (14,10) {$\hk(y_0)$};

\begin{scope}[xshift=19cm, yshift=-19cm]

\draw[thick, dashed] (1,1) -- (1,20);

\draw[ultra thick]  (1,20) -| (3,18) -| (4,17) -| (6,15) -- (8,15);

\coordinate (x0) at (5,3) node [black!60] at ($(x0)+(-.8,.5)$) {$x_0$};
\filldraw[black!60] (x0) rectangle +(1,1);
\draw[black!50] ($(x0)+(1,0)$) -- +(2,0);	% leg of x0
\draw[black!50] ($(x0)+(1,1)$) -- +(2,0);
\draw[black!50] ($(x0)+(0,1)$) -- +(0,13);	% arm of x0
\draw[black!50] ($(x0)+(1,1)$) -- +(0,13);

\node[anchor=west, blue] at (4,19) {$\delta-\hk(x_0)$};
\end{scope}

%\draw[blue, thick, xshift=19cm, yshift=-19cm] (1,20) -| (3,18) -| (4,17) -| (5,16) -| (3,17) -| (2,19) -| (-2,22) -| (-1,20) -- (1,20) -- cycle;
\draw[blue, thick, xshift=19cm, yshift=-19cm] (1,19) -| (-2,22) -| (-1,20) -- (1,20);
\fill[xshift=19cm, yshift=-19cm, pattern=north east lines, pattern color=blue!80] (1,19) -| (-2,22) -| (-1,20) -| (1,19) -- cycle;
\draw[blue!60, thick, xshift=19cm, yshift=-19cm] (1,20) -| (3,18) -| (4,17) -| (5,16) -| (3,17) -| (2,19) -- (1,19);
\fill[xshift=19cm, yshift=-19cm, pattern=north east lines, pattern color=blue!30] (1,20) -| (3,18) -| (4,17) -| (5,16) -| (3,17) -| (2,19) -| (1,20) -- cycle;
%\fill[xshift=19cm, yshift=-19cm, pattern=north east lines, pattern color=blue!80] (1,20) -| (3,18) -| (4,17) -| (5,16) -| (3,17) -| (2,19) -| (-2,22) -| (-1,20) -- (1,20) -- cycle;

\end{tikzpicture}
\end{center}
\caption{}\label{fig:rimhook}
\end{figure}

We need to show that $\hk(y_0)+(\delta-\hk(x_0))\notin \Z_{\geqq 0}\Pi_\cdiag.$
%By \eqref{eq:ininterval} and ,
%Since $x$ and $y$ are incomparable,
It follows from Lemma \ref{lem202} that
\begin{align*}
&0\oord \hk(y_0)\oord\hk(x_0)\oord\delta
\end{align*}
and thus we have
\begin{align*}
&0\oord 
\hk(y_0)+(\delta-\hk(x_0))
%=
%\delta-(\hk(x_0)-\hk(y_0))
\oord \delta,
\end{align*}
and
$$
\Supp \left( \hk(y_0)+(\delta-\hk(x_0)) \right)
=\Supp( \hk(y_0))\sqcup \Supp(\delta-\hk(x_0)).
$$
%\end{align*}
%
 By \eqref{eq:ininterval}, it holds that 
$y=y_{n+1}<x_{n+1}$. 
Thus
 $y_0<x_0$ and moreover
 $x_0$ and $y_0$ are not located in the same row or column.
Hence
\begin{equation}\label{eq:armlegnotin}
x_0^\mathrm{arm},x_0^\mathrm{leg}\notin
\Supp \left( \hk(y_0)+(\delta-\hk(x_0)) \right),
\end{equation}
%
%=\Supp(\hk(y_0))cup \Supp(\delta-\hk(x_0))
where
$x_0^\mathrm{arm}$ (resp. $x_0^\mathrm{leg}$) is the 
minimal element in $\{x_0\}\cup \arm(x_0)$ (resp. $\{x_0\}\cup \leg(x_0)$).
%$\{x_0^\mathrm{arm}\}=\arm(x_0)\cap \btm_\mathrm{min}$
%(resp. $\{x_0^\mathrm{leg}\} =\leg(x_0)\cap \btm_\mathrm{min}$).

%
%$\{x_0^\mathrm{arm}\}=\arm(x_0)\cap \btm_\mathrm{min}$ and 
%$\{x_0^\mathrm{leg}\} =\leg(x_0)\cap \btm_\mathrm{min}$.

Suppose that 
$$\hk(y_0)+(\delta-\hk(x_0))=\sum_{i=1}^r\beta_i$$ 
with 
$\beta_1,\dots,\beta_r\in \Pi_\cdiag$.
Then
$0\oord \beta_i\oord \delta$ ($i=1,\dots,r)$,
$0\oord \sum_{i=1}^r\beta_i\oord \delta$ and 
$$\Supp \left(\sum_{i=1}^r\beta_i\right)=\bigsqcup_{i=1}^r\Supp (\beta_i).$$ 
%and hence $\Supp(\beta_i)$ is included in 
%$\Supp \left( \hk(y_0)+(\delta-\hk(x_0)) \right)$.
Note that
each $\Supp (\beta_i)$ is an interval in $\cdiag$.
Combining with \eqref{eq:armlegnotin}, this implies that
$$\Supp(\beta_i)\subset
\Supp \left( \hk(y_0) \right)\ \text{ or }
%a subset of 
\ \Supp(\beta_i) \subset 
\Supp \left(\delta-\hk(x_0)) \right).$$
Thus there exist
$i_1,\dots,i_s$ for which we have
$\hk(y_0)=\beta_{i_1}+\cdots+\beta_{i_s}$, but this contradics 
\eqref{eq:hknotinpi}.
%the fact that
%$\hk(y_0)\notin \Z_{\geqq0}\Pi_\cdiag$.
%
Therefore $\hk(y)-\hk(x)= \hk(y_0)+(\delta-\hk(x_0)) $ 
cannot be a sum of elements in $\Pi_\cdiag$,
and thus  $\hk(x)$ and $\hk(y)$ are incomparable with respect to $\neword$.

The same argument implies that
 $\hk(x)$ and $\hk(y)$ are incomparable also in the case where  $N(\hk(y))=n-1$.
 \qed
\end{pf}
%

%We can give  another description for the order $\neword$.
%order on $R(w_\cdiag)$, which will interpolate
%$\hpord$ and $\neword$.
%%%%%%%%%%%%%%%%%%%%%%%%%%%%%%%%%

\begin{prop}\label{pr:tcord}
Let $\al,\beta\in R(w_\cdiag)$ with $\al\neword\beta$.
Then there exists a sequence
$$\al=\gamma_1,\gamma_2,\dots,\gamma_k=\beta$$
in $R(w_\cdiag)$
such that $\gamma_{i+1}-\gamma_i\in \Pi_\cdiag\ (i=1,\dots,k-1)$.

In other words,
the partial order $\neword$ on $R(w_\cdiag)$ coincides with
the transitive closure of the relations
\begin{equation}\label{eq:ord3}
\alpha \neword\hspace{-1.2mm} \raisebox{0.3mm}{$\cdot$} \,  \beta \ \text{whenever}\ \beta-\alpha \in \Pi_{\cdiag}.
\end{equation}
\end{prop}
\begin{pf}
Let $\tcord$ denote the transitive closure
of the relations above.
It follows from the same argument in the proof of Proposition \ref{prop:cdiagtotcord} 
that
$$x\le y\implies\hk(x)\tcord \hk(y)$$
for any $x,y\in\cdiag$.
It is clear that 
$$\hk(x)\tcord\hk(y)\implies
\hk(x)\neword\hk(y).$$
Combining with Theorem \ref{th:cylord}, the statement follows.\qed
\end{pf}

\subsection{Heaps}
%%%%%%%%%%%%%%%%%%%%%%%%%%%%%%%%%%%%%%%%%%%%%%%%%%%%%%%%%%
\if0
Let $\cdiag$ be a cylindric diagram.
%Heaps are a class of posets introduced by Stembridge associated with 
%(finite) $d$-complete posets.
%reduced expressions of elements in Coxeter groups \cite{Stem1996}.
Stembridge introduced posets called heaps asssociated with reduced expressions of 
elements in Coxeter groups \cite{Stem1996}.
In this section, we consider a semi-infinite analogue of heaps. 
%We have introduced a poset structure on $\cdiag$.
%We will introduce several more posets associated with $\cdiag$.
\fi
Let $\cdiag$ be a cylindric diagram.
%%%Let $\lex\in \LE(\cdiag)$
Recall that standard tableaux on $\cdiag$ have been defined as  order preserving bijection 
from $(\cdiag,\leq)$ to $(\Z_{\geqq1},\leqq)$.
Through the bijection $\lex$,  %Let
the set $\Z_{\geqq1}$ inherits a partial order from $\cdiag$, which 
we will investigate in this section.
%, which will be denoted by $\preceq_\lex$.%as before.
%Then $s_{i_k}=s(\lex^{-1}(k))$ and
%Let $\lex\in\LE(\cdiag)$. $$w_{\cdiag,\lex}[n]=s_{i_1}s_{i_2}\cdots\cdots s_{i_n}$$
%Define a partial order $\preceq_\lex$ on $\Z_{\geqq1}$
%%%%%%%%%%%%%%%%%%%%%%%%%%%%%%%%%%%%%%%%%%%%%%%%%%%%%%%%%%%%%%%%%%%%%%%%%%%%%

%In this section, we give a description of the posets $(\Z_{\geqq1}, \preceq_\lex)$,\begin{df}

%\hbox{$a \preceq_\lex b$ whenever $a \leqq b$ 
%and either $s_{i_a} s_{i_b} \neq s_{i_b} s_{i_a}$ or $i_a = i_b$.}
%The poset $(\Z_{\geqq1},\preceq_\lex)$ is called the \alert{heap} of $w_{\cdiag,\lex}$.  %(\cite{Stem2001}).
%\end{df}
%%%%%%%%%%%%%%%%%%%%%%%%%%%%%%%%%%%%%%%%%%%%%%%%%%%%%%%%%%%%%%%%%%%%%%

%%%%%%%%%%%%%%%%%%%%%%%%%%%%%%%%%%%%%%%%%%%%%%%%%%%%%%%%%%%%%%%%%%%%%%%%%%%%
\begin{df}
Let $\lex \in \LE(\cdiag)$.
%\noindent
Define a partial order $\preceq_\lex$ on $\Z_{\geqq1}$
as the transitive closure of the relations
\begin{align*}
&\hbox{$a \preceq_\lex b$ whenever $a \leqq b$ and}
\hbox{ either
$s_{i_a}s_{i_b}=s_{i_b}s_{i_a}$
or 
 $i_a =i_b$.}
% $s(\lex^{-1}(a)) s(\lex^{-1}(b)) \neq s(\lex^{-1}(b)) s(\lex^{-1}(a))$ or 
% $\con(\lex^{-1}(a)) =\con(\lex^{-1}(b))$.}
\end{align*}
where $i_k=\con(\lex^{-1}(k))$ for $k\in\Z$.
The poset $(\Z_{\geqq1},\preceq_\lex)$ is called the \alert{heap} of $w_{\cdiag,\lex}$.  
%(\cite{Stem2001}).
\end{df}
%%%%%%%%%%%%%%%%%%%
\begin{prop}\label{prop:cdiagandheap}
Let $\cdiag$ be a cylindric diagram and 
$\lex$ a standard tableau on $\cdiag$.
%\smallskip\noindent
%$\mathrm{(1)}$ 
Then, the map $\lex:\cdiag\to\Z_{\geqq 1}$ gives a poset isomorphism  
%between 
$$(\cdiag,\le)\cong (\Z_{\geqq1},\preceq_\lex).$$
%
%\smallskip
\if0
\noindent
$\mathrm{(2)}$ The map $\hk:\cdiag \to R(w_\cdiag)$ gives a   poset isomorphism
$$ (\cdiag,\leq) \cong (R(w_\cdiag),\hpord_\lex).$$
%for all $\lex\in\LE(\cdiag)$. 
In other words, the partial order $\hpord_\lex$ and $\neword$ on
$R(w_\cdiag)$ coincide.
\fi
\end{prop}
%%%%%%%%%%%%%%%%%%%%%%%%%

\begin{pf}
%The statements can be shown in a similar way as the finite case treated in \cite{Stem2001}.
%We just give a proof for (1) here.
%(1) 
Let $x,y\in \cdiag$.
 Suppose that $x < y$ is a covering relation in $\cdiag$.
Then $y=x-(1,0)$ or $y=x-(0,1)$ and
it is easy to see that 
 $\lex(x)<\lex(y)$ and $s(x)s(y) \neq s(y)s(x)$.
Hence $\lex(x) \preceq_\lex \lex(y)$.

Conversely,
suppose that $\lex(x) \prec_\lex \lex(y)$ is a covering relation in $\Z_{\geqq1}$.
%By definition, we have $a<b$ and either $s_{i_a} s_{i_b} \neq s_{i_b} s_{i_a}$ or $i_a = i_b$.
%Put $\lex^{-1}(a) = x$ and $\lex^{-1}(b) =y$.
%Suppose that $x \nless y$.
Then $s(x)s(y) \neq s(y)s(x)$ or $\con(x)=\con(y)$,
and hence  $\con(x)-\con(y) \neq 0,\pm1$.
By Proposition \ref{prop:content} (1),
$x$ and $y$ are comparable.
Since $\lex$ is order preserving, 
%we have $x \ngtr y$.
%Therefore, 
we must have $x<y$, and hence $\hk$ is a poset isomorphism.
%and hence, $(\cdiag,\le)$ and $(\Z_{\geqq1},\preceq_\lex)$ are isomorphic.
%The statement (2) can be shown in a similar way as the finite case treated in \cite{Stem2001},
%or as Proposition \ref{prop:cdiagandheap2}, which we will see soon.
\qed
\end{pf}

%%%%%%%%%%%%%%%%%%%%%%%%%%%%%%%%%%%%%%%%%%%%%%%%%%%%%%%%%%%%%%%%%%%%%%
The posets $(\Z_{\geqq1}, \preceq_\lex)$ are 
thought as semi-infinite analogue of heaps introduced by
Stembridge \cite{Stem1996}.
Stembridge also introduced the heap order on the inversion sets.
We treat a slightly modified version of heap order by Nakada \cite{Nak2008}.
% by Nakada. 
%In classical case. 
%Several orders in $R(w_\cdiag[n])$ which interpolate $(R(w_\cdiag[n])$
%and $\cdiag$ have been proposed (\cite{Stem1996}, Nakada).
%%%%%%%%%%%%%%%%%%%%%%%%%%%%%%%%%%%%%%%%%%%%%%%%%%%%%%%%%%%%%%%%%%%%%%%
\begin{df}
Define  a partial order $\hpord$ on $R(w_{\cdiag})$
as the transitive closure of the relations
$$
\hbox{$\alpha \hpord
\beta$ whenever $\alpha \oord \beta$ and 
$\bra\alpha,\beta^\vee\ket\neq 0$}.
%\hbox{$\gam(a) \hpord_\lex \gam(b)$ whenever $a \leqq b$ and $\bra\gam(a),\gam(b)^\vee\ket\neq 0$}.
$$
%The order $\hpord_{\lex}$ is called the \alert{heap order} on $R(w_{\cdiag})$.
\end{df}
%%%%%%%%%%%%%%%%%%%%%%%%%%%%%%%%%%%%%%%%%%%%%%%%%%%%%%%%%%%%%%%%%%%%%%%%

%%%%%%%%%%%%%%%%%%%%%%%%%%%%%%%%%%%%%%%%%%%%%%%%%%%%%%%%%%%%%%%%%%%%%%%%
\if0
%The partial order in $R(w_\cdiag)$ defined in the following is also called the heap order:
%%%%%%%%%%%%%%%%%%%%%%%%%%%%%%%%%%%%%%%%%%%%%%%%%%%%%%%%%%%%%%%%%%%%%%%
\begin{df}
%\noindent
%(2)
 Define  a partial order $\hpord_{\lex}$ on $R(w_{\cdiag})$
as the transitive closure of the relations
$$
\hbox{$\alpha \hpord_\lex 
\beta$ whenever $\lex( \hk^{-1}(\alpha)) \leqq \lex (\hk^{-1}(\beta))$ and 
$\bra\alpha,\beta^\vee\ket\neq 0$}.
%\hbox{$\gam(a) \hpord_\lex \gam(b)$ whenever $a \leqq b$ and $\bra\gam(a),\gam(b)^\vee\ket\neq 0$}.
$$
%The order $\hpord_{\lex}$ is called the \alert{heap order} on $R(w_{\cdiag})$.
\end{df}
%%%%%%%%%%%%%%%%%%%%%%%%%%%%%%%%%%%%%%%%%%%%%%%%%%%%%%%%%%%%%%%%%%%%%%%%
\fi

%%%%%%%%%%%%%%%%%%%%%%%%%%%%%%%%%%%%%%%%%%%%%%%%%%%%%%%%%%%%%%%%%

%The following order in $R(w_\cdiag)$ is proposed by Nakada \cite{Nakada}.
% by Nakada. 
%In classical case. 
%Several orders in $R(w_\cdiag[n])$ which interpolate $(R(w_\cdiag[n])$
%and $\cdiag$ have been proposed (\cite{Stem1996}, Nakada).
%%%%%%%%%%%%%%%%%%%
\begin{prop}\label{prop:cdiagandheap2}
The map $\hk:\cdiag \to R(w_\cdiag)$ gives a   poset isomorphism
$$ (\cdiag,\leq) \cong (R(w_\cdiag),\hpord).$$
%for all $\lex\in\LE(\cdiag)$. 
In other words, the partial order $\hpord$ and $\neword$ on
$R(w_\cdiag)$ coincide.
%\fi
\end{prop}
%%%%%%%%%%%%%%%%%%%%%%%%%
\begin{pf}
%%%%%%%%%%%%%%%%%%%%%%%%%%%%%%%%%%%%%%%%%%
 Let $x,y\in \cdiag$.
Suppose that $x < y$ is a covering relation in $\cdiag$.
Then $\hk(x)\oord\hk(y)$ and
$\hk(y)-\hk(x)\in\Pi_\cdiag\subset R\sqcup\Z\delta$.
We have
$$%\bra\hk(y)-\hk(x),\hk(y)^\vee\ket=
\bra\hk(y)-\hk(x),\hk(y)^\vee\ket
=2-\bra\hk(x),\hk(y)^\vee\ket.$$
If $\bra\hk(y),\hk(x)^\vee\ket=0$ then
$\hk(y)-\hk(x) \equiv \hk(y)\mod \Z \delta$ by Lemma \ref{lem:moddelta}, 
and thus $\hk(x)=k\delta$ for some $k\in\Z$.
This is a contradiction.
Therefore  $\bra\hk(x),\hk(y)^\vee\ket\neq 0$, from which it follows that
 $\hk(x)\hpord \hk(y)$.

Next, suppose that $\hk(x)\hpord \hk(y)$ is a covering relation.
%Then $\hk(x)\oord\hk(y)$ and $\bra\hk(x),\hk(y)^\vee\ket\neq 0$.
%
Put $x_0=x+N(\hk(x))(0,\ell)$ and $y_0=x+N(\hk(y))(0,\ell)$.
Then  $\hk(x_0)=\hk(x)-N(\hk(x))\delta,\ \hk(y_0)=\hk(y)-N(\hk(y))\delta$
and 
\begin{equation}\label{eq:nonzeroprod}
\bra\hk(x_0),\hk(y_0)^\vee\ket=\bra\hk(x),\hk(y)^\vee\ket\neq 0
\end{equation}
by assumption.

%It is enough to show that $x$ and $y$ are comparable.
We assume that $x$ and $y$ are incomparable. %and wll deduce a contradiction.
Then as $\hk(x)\oord\hk(y)$, 
we have $N(\hk(y))=N(\hk(x))+1$ and 
\begin{equation}\label{eq:x0y0}
\hk(y_0)\oord\hk(x_0)\oord \delta
\end{equation}
by Lemma \ref{lem202}.
Moreover, by \eqref{eq:armlegnotin} in the proof of Theorem \ref{th:cylord},
we have
\begin{equation}\label{eq:y0notin}
y_0\notin \arm(x_0)\cup \leg(x_0).
\end{equation}
(See also Figure \ref{fig:rimhook}.)

Recall that positive roots $\hk(x_0)$ and $\hk(y_0)$ can be expressed as 
$\hk(x_0)=\al_{ij}$ and $\hk(y_0)=\al_{kl}$
for some
$i,j,k,l\in \Z$ with $i<j$, $k<l$.
By \eqref{eq:x0y0} and \eqref{eq:y0notin},
the indices $i,j,k$ and $l$ can be chosen in such a way  that they satisfy
 $j-i\leqq \kappa-1$ and $i<k<l<j$.
Thus we have
\begin{align*}
\bra\hk(x_0),\hk(y_0)^\vee\ket&=
\bra\al_{ij},\al_{kl}^\vee\ket=\bra\al_{k-1\, l+1},\al_{kl}^\vee\ket\\
%\bra\al_i+\al_{i+1}+\cdots+\al_{j-1},(\al_k+\al_{k+1}++\cdots+\al_{l-1})^\vee\ket\\
%&=\bra\al_{k-1}+\al_{k}+\cdots+\al_{l-1}+\al_l,(\al_k+\al_{k+1}+\cdots+\al_{l-1})^\vee\ket\\
%&=\bra\al_{}+\al_{i+1}+\cdots+\al_{k-1},(\al_k+\al_{k+1}+\cdots+\al_{l-1})^\vee\ket
%\\
&=\bra\al_{k-1},\al_{kl}^\vee\ket+\bra\al_{k},\al_{kl}^\vee\ket
+\sum_{d=k+1}^{l-1}\bra\al_{d},\al_{kl}^\vee\ket
+\bra\al_{l-1},\al_{kl}^\vee\ket+\bra\al_{l},\al_{kl}^\vee\ket\\
&=-1+1+0+1-1
=0
\end{align*}
This contradicts \eqref{eq:nonzeroprod}.
Therefore $x$ and $y$ are comparable, and thus $x<y$ as $\hk(x)<\hk(y)$.
%%%%%%%%%%%%%%%%%%%%%%%%%%%%%%%%%%%%%%%%%%%%%%%%%%%%%%%%%%%%%%%%%%%%%%%%
\qed
\end{pf}

%%%%%%%%%%%%%%%%%%%%%%%%%%%%%%%%%%%%%%%%%%%%%%%%%%%%%%%%%%%%%%%
\section{
Poset structure of the set of order ideals}
\label{sec:2ndresult}
%%%%%%%%%%%%%%%%%%%%%%%%%%%%%%%%%%%%%%%%%%%%%%%%%%%%%%
\subsection{Standard tableaux on cylindric skew diagrams}

For a poset $P$, let $\I(P)$ denote 
the set of proper order ideals and regard $\I(P)$
as a poset with the inclusion relation.

 Let $\omega\in \Z_{\geqq 1}\times \Z_{\leqq -1}$ and 
fix a cylindric  diagram $\cdiag$ in $\cyli_\peri$.
%%%%%%%%%%%%%%%%%%%%%%%%%%%%%%%%%%%%%%%%%%%%%%%%%%%%%%%%%%%%%%%%%%%%%%%%%%%%%%%%
%\if0
%Recall that proper order ideal of $\cdiag$ is cylindric skew diagerams.
In this section, we will investigate the poset structure of 
the set $\I(\cdiag)$ of order ideals of $\cdiag$, in other words,
cylindric skew diagrams included in $\cdiag$.

%and 
%of proper order ideals (cylindric skew diagrams) of $\cdiag$
%can be seen as a poset with the inclusion relation.
%in terms of associated word $w_\cdiag[n]$.

Recall that any cylindric skew diagram $\ideal\in \I(\cdiag)$ is a finite set
and $\I(\cdiag)=\bigsqcup_{n=0}^\infty \I_n(\cdiag)$,
where 
$$\I_n(\cdiag)=\{\ideal\in\I(\cdiag)\mid |\ideal|=n\}.$$
% as before.
%%%%%%%%%%%%%%%%%%%%%%%%%%%%%%%%%%%%%%%%%%%%%%%%%%%%
%\begin{lem}
%Let $\cdiag$ be a cylindric diagram.
%%
%Then for each $n\in\Z_{\ge 1}$, the correspondence
%$$\lex\mapsto \lex^{-1}([1,n])$$
%gives a surjective map $\LE(\cdiag)\to \I(\cdiag)[n]$.
%\end{lem}

For $\ideal\in\I_n(\cdiag)$ and $\lex\in\LE(\ideal)$,
define a word $w_{\ideal,\lex}$ by
%$w_{\cdiag/\ofil,\lex}\in W$ and a semi-infinite word  $w_{\ofil,\lex} $ by
\begin{align}\label{eq;wzeta}
&w_{\ideal,\lex} = s(\lex^{-1}(1)) s(\lex^{-1}(2)) \cdots s(\lex^{-1}(n)).
\end{align}
%where $n=|\xi|$.
We sometimes regard $w_{\ideal,\lex}$ as a Weyl group element.

 %Moreover,  $w_{\ideal,\lex}$ is independent of the choice of $\lex$
%satisfying \eqref{eq:lexcondition}.
% and $w_{\ofil}$ is reduced and fully commutative semi-infinite word.
%%%%%%%%%%%%%%%%%%%%%%%%%%%%%%%%%%%%%%%
\begin{prop}
\label{prop:w_ideal}
%Let $\lex\in\LE(\ideal)$.
The word $w_{\ideal,\lex}$ is reduced.
As an element of  Weyl group, $w_{\ideal,\lex}$ is fully commutative
and independent of $\lex$.
\end{prop}
%%%%%%%%%%%%%%%%%%%%%%
\begin{pf}
It follows from Lemma \ref{lem:standideal} that the standard tableau $\lex$ on $\ideal$ can be extended to a standard tableau 
$\tilde\lex$ on $\cdiag$, for which 
we have $w_{\cdiag,\tilde\lex}[n]=w_{\ideal,\lex}$. 
By Proposition \ref{prop:reduced} and Proposition \ref{prop;pluscule},
the right hand side of \eqref{eq;wzeta} is a reduced expression and
 $w_{\ideal,\lex}$ is a fully commutative element of $W$.
It follows from Proposition \ref{prop:Randhk}
that $$R(w_{\ideal,\lex})=\{\hk(x)\mid x\in\ideal\}.$$ 
Hence the set $R(w_{\ideal,\lex})$ is independent of $\lex$
and so is $w_{\ideal,\lex}$.
%Since the inversion set $R(w)$ determines $w$ as a Weyl group element,
%t holds
%This implies  that $w_{\ideal,\lex}$ is independent of $\lex$.
\qed
\end{pf}
%%%%%%%%%%%%%%%%%%%%%%%%%%%%%%%%%%%%%%%%%%%%%%%
We denote by $w_\ideal$ the Weyl group element determined by the word
$w_{\ideal,\lex}$ for a/any standard tableau $\lex\in\LE(\ideal)$.
%%%%%%%%%%%%%%%%%%%%%%%%%%%%%%%%%%%%%%%%%%%%%%%%%%%%%%%%%
\begin{lem}[See\ {\cite[Theorem\ 3.2]{Stem1996}}]
\label{lem:lex to red}
%Let $\cdiag$ be a cylindric diagram and
%$\ideal\in\I_n(\cdiag)$.
 %with $|\ideal|=n$. 
%$\lex \in \LE(\cdiag)$.
%Then the map
The map
$$
\lex \mapsto w_{\ideal,\lex}= s(\lex^{-1}(1)) s(\lex^{-1}(2)) \cdots s(\lex^{-1}(n))$$
gives a bijection from 
$\LE(\ideal)$ to the set of reduced expressions for $w_{\ideal}$.
\end{lem}
%%%%%%%%%%%%%%%%%%%%%%%%%%%%%%%%%%%%%

%%%%%%%%%%%%%%%%%%%%%%%%%%%%%%%%%%%%%%%%%%%%%%%%%%%%%%%%%%%%%%%%%%%%%%
\begin{pf}
First, we prove that the correspondence is injective.
For $\lex_1,\lex_2 \in \LE(\ideal)$,
consider two words $w_{\ideal,\lex_1}
=s(p_1)s(p_2)\cdots s(p_n)$ and
$w_{\ideal,\lex_2}=s(q_1)s(q_2)\cdots s(q_n)$,
where $p_k=\lex_1^{-1}(k)$ and $q_k=\lex_2^{-1}(k)$.
%which are reduced.
%
Assume that $w_{\ideal,\lex_1}=w_{\ideal,\lex_2}$ as words.
Then $\con(p_1)=\con(q_1)$ and it holds that
$p_1$ and $q_1$ are minimal elements of $\ideal$.
Hence we have $p_1=q_1$.
Inductively, we have $p_k=q_k$ for any $k \in [1,n]$ by similar argument.

\if0
We will show that $p_k=q_k$ for any $k \in [1,n]$ by induction on $k$.
%
If $k=1$,
then $p_1$ and $q_1$ are minimal elements of $\ideal$.
Then
the two cells $p_1$ and $q_1$ are comparable,
and hence $p_1=q_1$.
%
Suppose that $p_j=q_j$ for any $j \in [1,k-1]$.
Then $p_k$ and $q_k$ are minimal elements of $\ideal \setminus \{p_1,\dots,p_{k-1}\}$,
and hence $p_k=q_k$ by the same argument.
\fi

Next, we prove that the map is surjective.
Take $\lex\in \LE(\ideal)$ and
put $p_j=\lex^{-1}(j)$ $(j\in[1,n])$. 
Then 
$w_{\ideal,\lex}= s(p_1) s(p_2) \cdots s(p_n)$,
which is a reduced expression of $w_\ideal$.
%and 

%Hence it is enough to show that 

Suppose that $s(p_k)s(p_{k+1})=s(p_{k+1})s(p_{k})$.
%  commute for some $k \in [1,n-1]$.
%
%Let $v=s_{i_1} s_{i_2} \cdots s(p_{k-2}) \cdot ts \cdot  s(p_{k+1})  \cdots s_{i_n}$ be a reduced expression for $w_{\cdiag,\lex}[n]$.
%Since $\lex$ is order preserving, $p_k$ and $p_{k+1}$
%are incomparable if $p_k$ is not covered by $p_{k+1}$.
Then $\con(p_k) - \con(p_{k+1}) \neq \pm1$,
and thus $p_k$ is not covered by $p_{k+1}$. %by Proposition \ref{prop:content},
%and hence 
This means that $p_k$ and $p_{k+1}$ are incomparable.
% as $\lex$ is order preserving.
%It can be seen that $p_k$ and $p_{k+1}$ must be incomparable.
%Actually,
%then there exists $q$ such that $p_k<q<p_{k+1}$.
Define the map $\lex^{(k)}:\ideal\to [1,n]$ by
$$
\lex^{(k)}(p_j)=
\begin{cases}
{k+1} & \text{if $j=k$}, \\
{k} & \text{if $j=k+1$}, \\
{j} & \text{otherwise}. \\
\end{cases}
$$
Then $\lex^{(k)}\in\LE(\ideal)$ 
%is a standard tableau of $\ideal$ 
and
$w_{\ideal,\lex^{(k)}}=s(p_1)s(p_2)\cdots s(p_{k+1})s(p_k)\cdots s(p_n)$.
Now fully-commutativity of $w_\ideal$ implies the surjectivity.
\qed
\end{pf}

%%%%%%%%%%%%%%%%%%%%%%%%%%%%%%%%%%%%%%%%%%%%%%%%%%%%%%%%%%
%\fi

%%%%%%%%%%%%%%%%%%%%%%%%%%%%%%%%%%%%%%%%%%%%%%%%%%%%%%%%
\subsection{Bruhat intervals}
%%%%%%%%%%%%%%%%%%%%%%%%%%%%%%%%%%%%%%%%%%%%%%%%%%%%%%%%
%%%%%%%%%%%%%%%%%%%%%%%%%
For $v,w \in W$,
we write $v \precdot w$ if $\ell(w)=\ell(v)+1$ and $w=v s_i$ for some simple reflection $s_i$.
Write $v \rble w$ if there is a sequence
$v=w_0 \precdot w_1 \precdot \cdots \precdot w_n= w$.
It is clear that the relation $\rbleq$ is a partial order of $W$,
and it is called the \alert{weak right Bruhat order}.
%We define the \alert{weak left Bruhat order} $\lbleq$ in the same way.
%%%%%%%%%%%%%%%%%%%%%%%%%%%
%

For $w \in W$, we define
$$[e,w] =\{ x\in W\mid e\rbleq x\rbleq w\}.$$
%where $\rbleq$ denotes the weak right Bruhat order (\S 3).
Note that when $\ell(w)=n$, we have
\begin{equation}
[e,w] = \left\{ s_{i_1} s_{i_2}\cdots s_{i_k} \in W \ \middle|\ 
\begin{tabular}{l}
$0\leqq k \leqq n \text{ and there exist }i_{k+1},\dots, i_n\text{ such that }$
\\
$s_{i_1} 
\cdots s_{i_k} s_{i_{k+1}}
 \cdots s_{i_n} \text{ is a reduced expression for }w$
\end{tabular}
\right\}.
\end{equation}
%where $\rleq$ denotes the \alert{weak right Bruhat order} in $W$.
%
%%%%%%%%%%%%%%%%%%%%%%%%%%%%%%%%%%%%%
%%
%
Let $\cdiag$ be a cylindric diagram.
For $\lex\in\LE(\cdiag)$, we define 
$$[e,w_{\cdiag,\lex})=\bigcup_{n=1}^\infty[e,w_{\cdiag,\lex}[n]].$$
We will see that the ``semi-infinite Bruhat interval" $[e, w_{\cdiag,\lex})$
is actually independent of $\lex\in\LE(\cdiag)$.
\newcommand\lexs{{\mathfrak{s}}}
%%%%%%%%%%%%%%%%%%%%%%%%%%%%%%%%%%%%%%%%%%%%%%%%%%%%
\begin{lem}\label{lem4}
Let $\lex_1$ and ${\lex_2}$ be two standard tableaux  on $\cdiag$.
%and 
%$w_\lex$ and $w_\delta$ be the expressions corresponding to $\lex$ and $\delta$.
Then for each $n \geqq 1$,
there exist $r\geqq n$ 
%reduced semi-infinite word
and $\lexs\in \LE(\cdiag)$ for which it holds that
$w_{\cdiag,\lexs}[r]=w_{\cdiag, \lex_1}[r]$ as elements of $W$ and 
$w_{\cdiag,\lexs}[n] = w_{\cdiag, \lex_2}[n]$ as words.
\end{lem}
%%%%%%%%%%%%%%%%%%%%%%%%%%%%%%%%%%%%%%%%%%%%%%%%%%%%%%%%%%%%%%%%%
\begin{pf}
Choose  $r\geqq n$ such that $\lex_2^{-1}[1,n]\subset \lex_1^{-1}[1,r]$.
Put $\ideal_1=\lex_1^{-1}[1,r]$ and $\ideal_2=\lex_2^{-1}[1,n]$.
Note that 
%$\ideal_1$ and $\ideal_2$ are cylindric diagrams and 
$\ideal_1\setminus \ideal_2$ is an order ideal of the cylindric diagram $\cdiag\setminus
\ideal_2$.
Take $\lex\in\LE(\cdiag\setminus \ideal_2)$ such that 
$\lex^{-1}[1,r-n]=\ideal_1\setminus\ideal_2$
(Lemma \ref{lem:standideal}).
Define a map $\lexs:\cdiag\to \Z_{\geqq1}$ by
\begin{equation*}
\lexs(p)=
\begin{cases}
\lex(p)+n \ \ &(p\in \cdiag\setminus\ideal_2)\\
\lex_2(p) \ \ &(p\in \ideal_2)
\end{cases}
\end{equation*}
Then we have $\lexs\in\LE(\cdiag)$, 
%$\lexs^{-1}[1,m]=\ideal_1$ and $\lexs^{-1}[1,n]=\ideal_2$,
which satisfies the desired conditions by Proposition \ref{prop:w_ideal}. 
\qed
\end{pf}
%%%%%%%%%%%%%%%%%%%%%%%%%%%%%%
%%%%%%%%%%%%%%%%%%%%%%%%%%%%%%%%%%%%%%%%%%%%%%%%%%%%%%%%%%%%%%%%%%%%%%%%%%%%%
\begin{prop}\label{prop:cile}
Let $\lex_1$ and $\lex_2$ be two standard tableaux of $\cdiag$.
Then
$$
[e,w_{\cdiag,\lex_1}) = [e,w_{\cdiag,\lex_2}) \ \text{as subsets of $W$}.
$$ 
\end{prop}
%%%%%%%%%%%%%%%%%%%%%%%%%%%%%%%%%%%%%%%%%%%%%%%%%%%%%%%%%%%
\begin{pf}%[Proposition \ref{prop:cile}]
%We write $w_{\cdiag,\lex_1}=w_1$ and $w_{\cdiag,\lex_2}=w_2$.
%
Let $n \geqq 1$.
By Lemma \ref{lem4},
there exist $r\geqq n$ and $\lexs\in\LE(\cdiag)$
such that $w_{\cdiag,\lexs}[r]=w_{\cdiag,\lex_1}[r]$  and 
$w_{\cdiag,\lexs}[n] = w_{\cdiag,\lex_2}[n]$.
Now we have
$$[e,w_{\cdiag,\lex_2}[n]]=[e,w_{\cdiag,\lexs}[n]]\subset [e,w_{\cdiag,\lexs}[r]]
\subset [e,w_{\cdiag,\lex_1}[r]].$$
Hence we obtain
$$
[e,w_{\cdiag,\lex_2}) = \bigcup_{n=1}^\infty[e,w_{\cdiag,\lex_2}[n]] \subset [e,w_{\cdiag,\lex_1}).
$$
Similarly,
we obtain $[e,w_{\cdiag,\lex_1}) \subset [e,w_{\cdiag,\lex_2})$,
and hence $[e,w_{\cdiag,\lex_1}) = [e,w_{\cdiag,\lex_2})$.
\qed
\end{pf}
%%%%%%%%%%%%%%%%%%%%%%%%%%%%%%%%%%%%%%%%%%%%%%%%%%%%%%%%%%%%

%%%%%%%%%%%%%%%%%%%%%%%%%%%%%%%
We denote $[e,w_{\cdiag,\lex})$ just by  $[e,w_{\cdiag})$
in the rest.
We have $$[e,w_\cdiag)=\bigcup_{\ideal\in\I(\cdiag)}[e,w_\ideal]$$
by the following lemma:
%
\if0
Note that $w_\ideal\in [e,w_\cdiag)$ for $\ideal\in\I(\cdiag)$
since there exist $\lex\in\LE(\cdiag)$ and $n$ such that $w_\ideal=w_{\cdiag,\lex}[n]
\in [e,w_\cdiag)$.
\fi
%%%%%%%%%%%%%%%%%%%%%%%%%%%%%%%%%%%%%%%%%%%%%%%%%%%%%%%%%%%%%%%%%%%%%%%%%%
\begin{lem}\label{lem:intervalcondition}
Let $v\in W$. Then 
$v\in[e,w_\cdiag)$ if and only if $v=w_\ideal$ for some $\ideal\in\I(\cdiag)$.
\end{lem}
%%%%%%%%%%%%%%%%%%%%%%%%%%%%%%%%%%%%%%%%%%%%%%%%%%%%%%%%%%%%%%%%%%%%%%%%%%%
\begin{pf}
Let $v\in [e,w_{\cdiag})$.
Then $v\in [e,w_{\cdiag,\lex}[n]]$
for some
$\lex\in \LE(\cdiag)$ and $n$.
By Lemma \ref{lem:lex to red},
there exist  $\lex'\in \LE(\cdiag)$ and $k$ such that
$v=w_{\cdiag,\lex'}[k]$.
Putting $\ideal={\lex'}^{-1}[1,k]$,
we have $v=w_{\ideal}$.

%
\if0
Then $v\rbleq w_{\cdiag,\lex}[n]=w_{\ideal'}$ with $\ideal'=\lex^{-1}[1,n]$.
There exist $\lexs\in\LE(\ideal)$ and $r\leqq n$ such that
$v=w_{\ideal',\lexs}[r]$ by Lemma \ref{lem:lex to red}.
Therefore $v=w_{\ideal}$ with $\ideal=\lexs^{-1}[1,r]\in \I(\cdiag)$
\fi

Let $\ideal\in\I(\cdiag)$.
Then
there exist $\lex\in\LE(\cdiag)$ and $n$ such that $w_\ideal=w_{\cdiag,\lex}[n]$.
Therefore 
$w_\ideal \in [e,w_\cdiag)$.
\qed
\end{pf}
%$\mathrm{(1)}$ Let . Then 
%
%\smallskip\noindent
%$\mathrm{(2)}$ Let $w\in [e,w_\cdiag)$. Then $R(w)$ is an order ideal of $R(w_\cdiag)$.
%

The following theorem can be
seen as a semi-infinite version of the results established 
in \cite{Stem1996} (see also \cite{NarOka2019, Proc1999ME}).
%See also \cite{SuzToy2022} for cylindric analogue of the skew hook formulas in \cite{NarOka2019}

%%%%%%%%%%%%%%%%%%%%%%%%%%%%%%%%%%%%%%%%%%%%%%%%%%%%%%%%%%%%%%%%%%%%%%%%%%%
\begin{thm}\label{th:ideal}
Let $\cdiag$ be a cylindric Young diagram in $\cyli_\peri$.
%and $\lex \in \LE(\cdiag)$.

\smallskip\noindent
$\mathrm{(1)}$ The map 
$$\Phi: (\I(\cdiag),\subset) \to ([e,w_{\cdiag}),\rbleq)$$
given by $\Phi(\ideal) = w_{\ideal}$ is a poset isomorphism.

\smallskip\noindent
$\mathrm{(2)}$ The map 
$$\Psi: ([e,w_{\cdiag}),\rbleq)\to (\I(R(w_\cdiag)),\subset)$$
given by $\Psi(w) = R(w)$ is a poset isomorphism.
\end{thm}
%%%%%%%%%%%%%%%%%%%%%%%%%%%%%%%%%%%%%%%%%%%%%%%%%%%%%%%%%%%%%%%%%%%%%%%%%%%%%%%%
\begin{pf}
We will show (1) and (2) togather.
Note that the poset isomorphism 
$\hk:
\cdiag\to R(w_\cdiag)$ induces a poset isomorphism
$\I(\cdiag)\to \I(R(w_\cdiag))$, 
under which $\ideal\in \I(\cdiag)$ corresponds to
$$\{\hk(x)\mid x\in \ideal\}=R(w_\ideal)=\Psi\circ \Phi(\ideal).$$
%
%
%It follows from Lemma \ref{lem3} that 
%First we show that $\im \Phi= [e,w_\cdiag)$.
\if0
Let $v\in [e,w_{\cdiag})$ and take
 $\lex\in \LE(\cdiag)$ and $n$ such that
$v\in [e,w_{\cdiag,\lex}[n]]$.
Then $v\rbleq w_{\cdiag,\lex}[n]=w_\ideal$ with $\ideal=\lex^{-1}[1,n]$
and there exist $\lexs\in\LE(\ideal)$ and $r\leqq n$ such that
$v=w_{\ideal,\lexs}[r]$ by Lemma \ref{lem:lex to red}.
Therefore $v=w_{\ideal'}$ with $\ideal'=\lexs^{-1}[1,r]$ and the map $\Phi$ is surjective.
\fi
%
Hence  $\Psi\circ\Phi$ is bijective and thus  $\Phi$ is injective.
As
 $\Phi$ is surjective %$\im \Phi= [e,w_\cdiag)$
 by Lemma \ref{lem:intervalcondition},
$\Phi$ is bijective.
Thus $\Psi$ is also bijective.

We will show that $\Phi$ and $\Psi$ are order preserving.

Suppose that
$\ideal'$ covers $\ideal$, or equivalently that
$\ideal'=\ideal\sqcup\{x\}$ for a maximal element $x$ of $\ideal'$.
Then there exists
$\lex\in \LE(\ideal')$ satisfying $\lex^{-1}(n)=x$, for which we have
$$w_{\ideal'}=s(\lex^{-1}(1))s(\lex^{-1}(2))\cdots
s(\lex^{-1}(n-1))s(\lex^{-1}(n))=w_\ideal s(x),$$
This implies that $w_{\ideal'}$ covers $w_{\ideal}$. Hence $\Phi$ is order preserving.

It is easy to see that $v\rbleq w$ implies $R(v)\subset R(w)$.
Hence $\Psi$ is order preserving.

As we know that $(\Psi\circ\Phi)^{-1}$ is order preserving, it holds that
$\Phi^{-1}$ and 
$\Psi^{-1}$ are also 
order preserving.
\qed
\end{pf}
%%%%%%%%%%%%%%%%%%%%%%%%%%%%%%%%%%%%%%%%%%%%%%%%%%%%%%%%%

\begin{prop}
Let $\cdiag$ be a cylindric diagram. Then
$$[e,w_\cdiag)=\{w\in W\mid w\text{ is }\predom_\cdiag\text{-pluscule}\}$$
\end{prop}
%%%%%%%%%%%%%%%%%%%%%%%%%%%%%%%%%%%%%
\begin{pf}
It follows from Proposition \ref{prop;pluscule} that
any element of $[e,w_\cdiag)$ is $\predom_\cdiag$-pluscule.

Let $w\in W$ be $\predom_\cdiag$-pluscule and $w=s_{i_1}s_{i_2}\cdots s_{i_n}$ 
its reduced expression.
We will show that $w\in [e,w_\cdiag)$ by induction on $n=\ell(w)$.
By induction hypothesis, $v:=s_{i_1}s_{i_2}\cdots s_{i_{n-1}}$ belongs to 
$ [e,w_\cdiag)$, and thus $v=w_\ideal$ for some $\ideal\in\I(\cdiag)$.

Let  $x$ be the minimum element of $\con^{-1}(i_n)\cap (\cdiag\setminus \ideal)$
and put $\ideal'=\ideal\sqcup\{x\}$.
Take $\lex\in\LE(\ideal')$ such that $\lex(n)=x$.
 Then
 $w=s(\lex^{-1}(1))s(\lex^{-1}(2))\cdots s(\lex^{-1}(n))$.
Since $w$ is $\predom_\cdiag$-pluscule, %, it is fully commutative.
%by \cite[Proposition 2.3]{Stem2001},
%Hence it holds that 
if $i_n=i_k$ then there exist $j_+,j_-\in[k,n]$
such that $j_+=i_n+1$ and $j_-=i_n-1$ by \cite[Proposition 2.3]{Stem2001}.
This implies that
the subset $\ideal'$
%$=\{\lex^{-1}(k)\mid k\in[1,n]\}$ of $\cdiag$ 
satisfies the condition (v) in Proposition \ref{prop:skew}.
Therefore $\ideal'$ is a cylindric skew diagram in $\cdiag$
and $w=w_{\ideal\sqcup\{x\}}$.
Therefore $w\in[e,w_\cdiag)$.
\qed
\end{pf}
%%%%%%%%%%%%%%%%%%%%%%%%%%%%%%%%%%%%%%%%%%%%%%%%%%%%%%%%%%%%%%%%%%%%%%

%%%%%%%%%%%%%%%%%%%%%%%%%%%%%
\subsection{Skew diagrams and classical case}
%%%%%%%%%%%%%%%%%%%%%%%%%%%
\newcommand\sideal{{\eta}}
Let $\cdiag$ be a cylindric diagram in $\cyli_\peri$.
%%%%%%%%%%%%%%%%%%%%%%%%%%%%%%%%%%%%%%%%%%%%%%%%%%%%%%%%%%%%%%%
Let $\ideal\in\I_n(\cdiag)$ and take
 $\lex\in\LE(\cdiag)$ such that $\ideal=\lex^{-1}[1,n]$. 
Then we have
$w_{\cdiag,\lex}[n]=w_{\ideal}$ and 
$\hk(\ideal)=R(w_\ideal)$.
Thus the next theorem follows easily from 
Theorem \ref{th:cylord}:
%%%%%%%%%%%%%%%%%%%%%%%%%%%%%%%%%%%%%%%%%%%%%%%%%%%%%%%%%%%%
\begin{thm}\label{thm:skew1}
Let $\ideal\in\I_n(\cdiag)$.

\smallskip\noindent
$\mathrm{(1)}$ The map 
$\hk: (\ideal,\leq)\to (R(w_\ideal),\cylord)$
is a poset isomorphism.

\smallskip\noindent
$\mathrm{(2)}$
For $\lex\in\LE(\ideal)$, the map 
$\lex:(\ideal,\leq)\to ([1,n],\hpord_\lex)$ is a poset isomorphism.
\end{thm}
%%%%%%%%%%%%%%%%%%%%%%%%%%%%%%%%%%%%%%%%%%%%%%%%%%%%%%%%

Note that $\I(\ideal)=\{\sideal\in\I(\cdiag)\mid \sideal\subset \ideal\}$.  
Theorem \ref{th:ideal} implies the following:
%%%%%%%%%%%%%%%%%%%%%%%%%%%%%%%%%%%%%%%%%%%%%%%%%
\begin{thm}\label{thm:skew2}
Let $\ideal\in\I(\cdiag)$.

\smallskip\noindent
$\mathrm{(1)}$ 
The map $\Phi:(\I(\ideal),\subset)\to ([e,w_\ideal],\rbleq)$
given by $\Phi(\sideal)=w_\sideal$
is a poset isomorphism.

\smallskip\noindent
$\mathrm{(2)}$ 
The map $\Psi:([e,w_\ideal],\rbleq)\to (\I(R(w_\ideal)),\subset)$
given by $\Psi(w)=R(w)$
is a poset isomorphism.
\end{thm}
%%%%%%%%%%%%%%%%%%%%%%%%%%%%%%%%%%%%%%%%%%%%%%%%%%%%%%%%

%%%%%%%%%%%%%%%%%%%%%%%%%%%%%%%%%%%%%%%%%%%%%%%%%
\newcommand\cldiag{{\sem\lm/\sem\phi}}
\newcommand\clskew{{[\lm/\mu]}}
In the rest, we will see that 
 description for %classical (finite) 
non-cylindric diagrams
can be deduced  from the results above.
Let $m\in\Z_{\geqq1}$ and 
let  $\lm=(\lm_1,\dots,\lm_m),\ \mu=(\mu_1,\dots,\mu_m)$ be partitions such that 
$\lm_i\geqq\mu_i \geqq0 \ (i\in[1,m])$.
Under the notation in Section \ref{SS:cylindric},
the associated classical skew Young diagram is represented as the subset $\sem\lm/\sem\mu$
of $\Z^2$:
$$\sem\lm/\sem\mu=\left\{(a,b)\in\Z^2\mid
a\in[1,m],\ b\in[\mu_a+1,\lm_a]\right\}.$$
Note that the classical normal Young diagram associated with $\lm$ is
a special skew diagram $\sem\lm/\sem\phi$ with $\phi=(0,0,\dots,0)$.
%$$\sem\lm/\sem\phi=\left\{(a,b)\in\Z^2\mid
%a\in[1,m],\ b\in[1,\lm_a]\right\}.$$
\if0

\begin{align*}
\cldiag&=\left\{ (a,b)\mid a\in[1,m],\ b\in[1,\lm_a]\right\}\\
\clskew&=\left\{ (a,b)\mid a\in[1,m],\ b\in[\mu_a+1,\lm_a]\right\}
\end{align*}
\fi

%

%
To connect classical diagrams and cylindric diagrams, 
we take  $\ell\in\Z_{\geq1}$ such that 
$$\ell\geqq \lm_1-\mu_m.$$
Then the partitions 
$\lm,\mu$ are $\ell$-restricted, and moreover
it is easy to see that the skew diagram $\sem\lm/\sem\mu$ is isomorphic
to the cylindric skew  diagram $\cyl\lm/\cyl\mu=\pi(\sem\lm/\sem\mu)$ 
as a poset.
Under this identification $\sem\lm/\sem\mu=\cyl\lm/\cyl\mu$,
%it follows from 
Theorem \ref{thm:skew1} and Lemma \ref{lem201} %that
for the order ideal $\cyl\lm/\cyl\mu$ of the cylindric diagram $\cyl\lm$ imply the followings:
$$([1,n],\hpord_\lex)\cong (\sem\lm/\sem\mu,\leq)\cong (R(w_{\sem\lm/\sem\mu}),\cylord)
=(R(w_{\sem\lm/\sem\mu}),\oord)$$
for each $\lex\in \LE(\sem\lm/\sem\mu)=\LE(\cyl\lm/\cyl\mu)$,
and it follows from Theorem \ref{thm:skew2} that
$$(\I(\sem\lm/\sem\mu),\subset)\cong ([e,w_{\sem\lm/\sem\mu}], \rbleq)\cong 
(\I(R(w_{\sem\lm/\sem\mu}),\subset).$$
Remark that by redefining the content as 
$$\con(a,b)=b-a+m-\mu_m,$$
we have $\con (\sem\lm/\sem\mu)\subset [1,\kappa-1]$, and
$$w_{\sem\lm/\sem\mu}\in \bar W,\  R(w_{\sem\lm/\sem\mu})\subset \bar R,$$
where 
$\bar W$ and $\bar R$  denote the Weyl group and  the root system
of type $A_{\kappa-1}$ respectively.
%Note also that the order $\cylord$ is defined by
%$$\al\cylord\beta\Leftrightarrow \beta-\al\in\Pi_{\cyl\lm}$$
%and $\Pi_{\cyl\lm}$ is a subset of the set $\bar \Pi$ of simple roots in $\bar R$.

%%%%%%%%%%%%%%%%%%%%%%%%%%%%%%%%%%%%%%%%%%%%%%%%%%%%%%%%%%%%%%

%%%%%%%%%%%%%%%%%%%%%%%%%%%%%%%%%%%%%%%%%%%%%%%%%%%%
\newcommand\cohk{{\mathbf{coh}}}
%\begin{rem}
We will see the relation between the results above and preceding works.
Let $n \in \Z_{\geqq1}$ and $\lm$ be a partition of $n$.
% and let $\cldiag$ denote the associated classical Young diagram, where
%$\phi=(0,0,\dots,0)$.
Fix $\lex\in \LE(\cldiag)$ and
put
$$w_\lm:=w_{\cldiag}
=s(\lex^{-1}(n))s(\lex^{-1}(n-1))\cdots s(\lex^{-1}(1)).
$$
The element $w_\lm$ is independent of $\lex$
and it is called the Grassmannian permutation associated with $\lm$.
%$\bar W$, where 
%,
%is independent of $\lex$ and is called the Grassmanian permutation associated to $\lm$.

It  has been  shown  in \cite{Stem2001, NarOka2019}
that the  map 
$$\cohk:\cldiag \to R(w_\lm^{-1})$$
given by 
\begin{equation}\label{eq:defcohk}
\cohk(x)=s(\lex^{-1}(n))s(\lex^{-1}(n-1))\cdots 
s(\lex^{-1}(k+1))\al(\lex^{-1}(k)),
\end{equation}
where $k=\lex^{-1}(x)$,
leads an {\it dual isomorphism} of posets:
\begin{equation}\label{eq:cohkisom}
\cohk : (\cldiag,\leq)\to (R(w_\lm^{-1}),\oord),
\end{equation}
where $\oord$ is the ordinary order as before.

On the other hand, as a classical version of Theorem \ref{thm:skew1},
we have a poset isomorphism 
\begin{equation}\label{eq:hkisom}
\hk: (\cldiag,\leq)\to (R(w_\lm),\neword).
\end{equation}
Now define the map $\iota: R\to R$  
 by $\iota (\al )=-w_\lm^{-1}\al$.
Then it follows immediately from the expression \eqref{eq;hk} and \eqref{eq:defcohk}
that $\iota \circ \hk(x)=\cohk (x)$ for all $x\in \cldiag$.
%$\iota$ rekates the two isomorphisms \eqref{eq:cohkisom} and \eqref{eq:hkisom}. 
Therefore  we have the following:
%%%%%%%%%%%%%%%%%%%%%%%%%%%%%%%%%%%%%
\begin{prop}
The restriction of $\iota$ gives a dual poset isomorphism
$$\iota:(R(w_\lm),\cylord)\to (R(w_\lm^{-1}),\oord)$$
and  moreover $\iota \circ \hk=\cohk$.
In other words, the following diagram of poset isomorphisms commutes $:$
\begin{equation}
\xymatrix{
(\cldiag,\leq)\ar[r]^-{\hk}\ar[d]_-{\cohk}&
(R(w_{\lm}),\cylord)\ar[dl]^-{\iota}
%\ar@{}@<-1.5ex>[dl]|{\circlearrowright}
\\
(R(w_{\lm}^{-1}),\oord)^{\mathrm{op}}
}
\end{equation}
where $(R(w_{\lm}^{-1}),\oord)^{\mathrm{op}}$ denotes the poset obtained 
from  $(R(w_{\lm}^{-1}),\oord)$ by reversing the order.
\end{prop}
%\end{rem}

%%%%%%%%%%%%%%%%%%%%%%%%%%%%%%%%%%%%%

\end{document}